\documentclass[12pt,english,leqno]{article}
\usepackage[T1]{fontenc}
\usepackage[latin1]{inputenc}
\usepackage{babel}
\usepackage{amssymb}

\makeatletter

\providecommand{\LyX}{L\kern-.1667em\lower.25em\hbox{Y}\kern-.125emX\@}

\usepackage{amsfonts,amsmath}
\def \nz{\normalsize}

 \newcounter{secnum}

 \makeatother
 
\begin{document}

 \title{ACCURACY OF DIFFUSION APPROXIMATIONS FOR HIGH FREQUENCY MARKOV DATA. \thanks{
 This research was supported by grant 436RUS113/845/0-1
from the Deutsche Forschungsgemeinschaft and by grant 05-01-04004
from the Russian Foundation of Fundamental Researches. The first
author worked on the paper during a visit at the Laboratory of
Probability Theory and Random Models of the University Paris VI in
2005/2006. He is grateful for the hospitality during his stay.
 }\\
 }

 \author{Valentin KONAKOV \\
  \nz Central Economics Mathematical Institute, Academy of Sciences\\
  \nz Nahimovskii av. 47, 117418 Moscow, Russia\\ \nz E mail:
  kv24@mail.ru\\
  Enno MAMMEN \\
  \nz Department of Economics, University of Mannheim \\
  \nz L 7, 3-5, 68131 Mannheim, Germany\\
  \nz E mail: emammen@rumms.uni-mannheim.de}

 \date{\today{}}

 \maketitle
 \noindent We consider triangular arrays of Markov chains that converge
weakly to a diffusion process. Edgeworth type expansions of third
order for transition densities are proved. This is done for time
horizons that converge to $0$. For this purpose we represent the
transition density as a functional of densities of sums of i.i.d.
variables. This will be done by application of the parametrix
method. Then we apply Edgeworth expansions to the densities. The
resulting series gives our Edgeworth-type expansion for the
transition density of Markov chains. The research is motivated by
applications to high frequency data that are available on a very
fine grid but are approximated by a diffusion model on a more
rough grid. \vskip .2in
 \noindent \textsl{1991 MSC:} primary
 62G07, secondary 60G60


 \noindent \textsl{Short title:} Edgeworth expansions for Markov chains.

 \newpage

\section{Introduction.}

In this paper we study triangular arrays of Markov chains
$X_{k,h}$ ($0\leq k\leq n$) that converge weakly to a diffusion
process $Y_s$ ($0 \leq s \leq T$) for $n\rightarrow\infty$. Here
$h=T/n$ denotes the discretization step. We allow that $T$ depends
on $n$. In particular, we consider the case that $T \to 0$  for
$n\rightarrow\infty$. Our main result will give Edgeworth type
expansions for the transition densities. The order of the
expansions is $o(h^{-1-\delta}),\delta>0$. This is done for time
horizons that may converge to $0$. The research is motivated by
applications to high frequency data that are available on a very
fine grid but are approximated by a diffusion model on a more
rough grid. The work of this paper generalizes the results in
Konakov and Mammen (2005) in two directions. The time horizon is
allowed to converge to $0$ and also cases are treated with
nonhomogenous diffusion limits.

The theory of Edgeworth expansions is well developed for sums of
independent random variables. For more general models approaches
have been used where the expansion is reduced to models with sums
of independent random variables. This is also the basic idea of
our approach. We will make use of the parametrix method. In this
approach the transition density is represented as a nested sum of
functionals of densities of sums of independent variables.
Plugging Edgeworth expansions into this representation will result
in an expansion for the transition density.

Weak convergence of the distribution of scaled discrete time
Markov processes to diffusions has been extensively studied in the
literature {[}see Skorohod (1965) and Stroock and Varadhan
(1979){]}. Local limit theorems for Markov chains were given in
Konakov and Molchanov (1984) and Konakov and Mammen
(2000,2001,2002). In Konakov and Mammen (2000) it was shown that
the transition density of a Markov chain converges with rate
$O(n^{-1/2})$ to the transition density in the diffusion model.
For the proof there an analytical approach was chosen that made
essential use of the parametrix method. This method permits to
obtain tractable representations of transition densities of
diffusions that are based on Gaussian densities, see Lemma 1
below. Similar representations hold for discrete time Markov
chains $X_{k,h}$, see Lemma 3 below. For a short exposition of the
parametrix method, see Section 3 and Konakov and Mammen (2000).
The parametrix method for Markov chains developped in Konakov and
Mammen (2000) is exposed in Section 4. Applications to Markov
random walks are given in Konakov and Mammen (2001). In Konakov
and Mammen (2002) the approach is used to give Edgeworth-type
expansions for Euler schemes for differential equations. Related
treatments of Euler schemes can be found in [Bally and Talay
(1996a,b), Protter and Talay (1997), Jacod and Protter (1998),
Jacod (2004), Jacod, Kurtz, M\'{e}l\'{e}ard and Protter (2005) and
Guyon (2006)]. Standard references for the parametrix method are
the books by Friedman (1964) and Lady\u{z}enskaja, Solonnikov and
Ural\'{}ceva (1968) on parabolic PDEs [see also McKean and Singer
(1967)].

The paper is organized as follows. In the next section we will
present our model for the Markov chain and state our main result
that gives an Edgeworth-type expansion for Markov chains. In
Section 3 we will give a short introduction into the parametrix
method for diffusions. In Section 4 we will recall the parametrix
approach developed in Konakov and Mammen (2000) for Markov chains.
Technical discussions, auxiliary results and proofs are given in
Sections 5-7.

\section{Results.}
  Let $n\geq 2$ , $T=T(n)\leq 1$ \ and $h=T/n.$ Suppose that $q\left(
t,x,\cdot \right) ,$ $\left( t,x\right) \in \lbrack 0,1]\times \Bbb{R}^{d}$
is a given family of densities on $\Bbb{R}^{d}$, $\chi _{\nu }(t,x)$ is the $%
\nu -th$ cumulant corresponding to the density $q(t,x,\cdot )$ and $m$ is a
function from $[0,1]\times \Bbb{R}^{d}$ into $\Bbb{R}^{d}.$ We shall impose
the following conditions
\begin{description}
 \item[(A1)] $\ \int_{\Bbb{R}^{d}}yq\left(
t,x,y\right) dy=0,\;0\leq t\leq 1,\ x\in \Bbb{R}^{d}.$
 \item[(A2)] There exists positive constants $\sigma
_{\ast }$ and $\sigma ^{\ast }$ such that the covariance matrix $\
\sigma \left( t,x\right) =\int_{\Bbb{R}^{d}}yy^{T}q\left(
t,x,y\right) dy$ \ satisfies
\[
\sigma _{\ast }\leq \theta ^{T}\sigma \left( t,x\right) \theta \leq \sigma
^{\ast },
\]
for all $\left\| \theta \right\| =1$ and $t\in \lbrack
0,1],$ $x\in \Bbb{R}^{d}.$

\item[(A3)] There exists a positive integer $S^{\prime }$ and a
real nonnegative function $\psi \left( y\right) ,$ $y\in
\Bbb{R}^{d}$
satisfying $\sup_{y\in \Bbb{R}^{d}}\psi \left( y\right) <\infty $ \ and $%
\int_{\Bbb{R}^{d}}\left\| y\right\| ^{S}\psi \left( y\right) dy<\infty $ \
with $S=(S^{\prime }+2)d+4$ \ such that
\[
\left| D_{y}^{\nu }q\left( t,x,y\right) \right| \leq \psi \left( y\right)
,\;t\in \lbrack 0,1],\;x,y\in \Bbb{R}^{d}\;\left| \nu \right| =0,1,2,3,4
\]
\[
\left| D_{x}^{\nu }q\left( t,x,y\right) \right| \leq \psi \left( y\right)
,\;t\in \lbrack 0,1],\;x,y\in \Bbb{R}^{d}\;\left| \nu \right| =0,1,2.
\]
\end{description}
It follows from \textbf{(A2), }Lemma 5 and (\ref{eq:12e}),(\ref
{eq:12f}) below (applied with $h=1$ ) that the following condition
holds
\begin{description}
 \item[(A3${^{\prime }}$)]  For all $x,y\in R,$ $%
h>0,0\leq t,t+jh\leq 1,j\geq j_{0},$ with\ $\ j_{0}$ \ does not
depending on $x,t$
\[
\left| D_{x}^{\nu }q^{(j)}\left( t,x,y\right) \right| \leq Cj^{-d/2}\psi
\left( j^{-1/2}y\right) ,\left| \nu \right| =0,1,2,3
\]
for a constant $C<\infty .$\ \ Here $q^{(j)}(t,x,y)$ denotes the $j-$ fold \
convolution of \ $q$ \ for fixed $x$ as a function of $\ y$.
\[
q^{(j)}(t,x,y)=\int q^{(j-1)}(t,x,u)q(t+(j-1)h,x,y-u)du,
\]
$q^{(1)}(t,x,y)=q(t,x,y).$
\end{description}

It follows also from (A3) that for \ $%
1\leq j\leq j_{0}$%
\[
\int \left\| y\right\| ^{S}q^{(j)}(t,x,y)dy\leq C(j_{0\;}).
\]
\begin{description}
 \item[(B1)] The functions $m\left( t,x\right) $ and $\sigma
\left( t,x\right) $ and their first and second derivatives w.r.t.
$t$ and their derivatives up to the order six w.r.t. $x$ are
continuous and bounded uniformly in $t$ and $x.$ All these
functions are Lipschitz continuous with respect to $x$ with a
Lipschitz constant that does not depend on $t.$ The function $\chi
_{\nu }(t,x)$, $\left| \nu \right| =3,4,$ is Lipschitz continuous
with respect to $\ t$ \ with a Lipschitz constant that does not
depend on $x.$ A sufficient condition for this is the following
inequality
\[
\int_{\Bbb{R}^{d}}(1+\left\| z\right\| ^{4})\left| q\left( t,x,z\right)
-q\left( t^{\prime },x,z\right) \right| dz\leq C\left| t-t^{\prime }\right|
,0\leq t,t^{\prime }\leq 1,
\]
with a constant that does not depend on $\ x\in \Bbb{R}^{d}.$ Furthemore, $%
D_{x}^{\nu }\sigma \left( t,x\right) $ exists for $\left| \nu \right| =6$
and is Holder continuous w.r.t. $x$ with positive exponent and a constant
that does not depend on $t.$

 \item[(B2)] There exist\textbf{\ }\ $\varkappa <\frac{1}{5}$ $\
$such that $\liminf  $ $_{n\rightarrow \infty }T(n)n^{\varkappa
}>0$ (remind that we consider $T(n)\leq 1$ ).
\end{description}

Consider a family of Markov processes in $\Bbb{R}^{d}$ of the
following form
\begin{equation}
X_{k+1,h}=X_{k,h}+m\left( kh,X_{k,h}\right) h+\sqrt{h}\xi
_{k+1,h},\;X_{0,h}=x\in \Bbb{R}^{d},\;k=0,...,n-1,  \label{eq:001}
\end{equation}
where $\left( \xi _{i,h}\right) _{i=1,...,n}$ is an innovation
sequence satisfying the Markov assumption: the conditional
distribution of $\xi
_{k+1,h}$ given $X_{k,h}=x_{k},...,X_{0,h}=x_{0}$ depends only on $%
X_{k,h}=x_{k}$ and has conditional density $q\left( kh,x_{k},\cdot \right) .$
The conditional covariance matrix corresponding to this density is $\sigma
(kh,x_{k}).$ The transition densities of $\left( X_{i,h}\right) _{i=1,...,n}$
are denoted by $p_{h}\left( 0,kh,x,\cdot \right) .$

We will consider the process (\ref{eq:001}) as an approximation to
the
following stochastic differential equation in $\Bbb{R}^{d}:$%
\[
dY_{s}=m\left( s,Y_{s}\right) ds+\Lambda \left( s,Y_{s}\right)
dW_{s},\;Y_{0}=x\in \Bbb{R}^{d},\;s\in \lbrack 0,T],
\]
where $\left( W_{s}\right) _{s\geq 0}$ is the standard Wiener process and $%
\Lambda $ is a symmetic positive definite $d\times d$ matrix such that $%
\Lambda \left( s,y\right) \Lambda \left( s,y\right) ^{T}=\sigma \left(
s,y\right) .$ The conditional density of $Y_{t},$ given $Y_{0}=x$ is denoted
by $p\left( 0,t,x,\cdot \right) .\,\,$We shall consider \ the following
differential operators $L$ and $\widetilde{L}$ :
\[
Lf(s,t,x,y)=\frac{1}{2}\sum_{i,j=1}^{d}\sigma _{ij}(s,x)\frac{\partial
^{2}f(s,t,x,y)}{\partial x_{i}\partial x_{j}}+\sum_{i=1}^{d}m_{i}(s,x)\frac{%
\partial f(s,t,x,y)}{\partial x_{i}},
\]
\begin{equation}
\tilde{L}f(s,t,x,y)=\frac{1}{2}\sum_{i,j=1}^{d}\sigma _{ij}(s,y)\frac{%
\partial ^{2}f(s,t,x,y)}{\partial x_{i}\partial x_{j}}%
+\sum_{i=1}^{d}m_{i}(s,y)\frac{\partial f(s,t,x,y)}{\partial x_{i}}.
\label{eq:001a}
\end{equation}
To formulate our main result we need also\ the following operators
\[
L^{\prime }f(s,t,x,y)=\frac{1}{2}\sum_{i,j=1}^{d}\frac{\partial \sigma
_{ij}(s,x)}{\partial s}\frac{\partial ^{2}f\left( s,t,x,y\right) }{\partial
x_{i}\partial x_{j}}+\sum_{i=1}^{d}\frac{\partial m_{i}(s,x)}{\partial s}%
\frac{\partial f\left( s,t,x,y\right) }{\partial x_{i}}
\]
\begin{equation}
\widetilde{L}^{\prime }f(s,t,v,z)=\frac{1}{2}\sum_{i,j=1}^{d}\frac{\partial
\sigma _{ij}(s,y)}{\partial s}\frac{\partial ^{2}f\left( s,t,x,y\right) }{%
\partial x_{i}\partial x_{j}}+\sum_{i=1}^{d}\frac{\partial m_{i}(s,y)}{%
\partial s}\frac{\partial f\left( s,t,x,y\right) }{\partial x_{i}}.
\label{eq:001b}
\end{equation}
and the convolution type binary operation $\otimes :$%
\[
f\otimes g\left( s,t,x,y\right) =\int_{s}^{t}du\int_{R^{d}}f\left(
s,u,x,z\right) g\left( u,t,z,y\right) dz.
\]

Konakov and Mammen (2000) obtained a nonuniform rate of
convergence for the difference $p_{h}\left( 0,T,x,\cdot \right)
-p\left( 0,T,x,\cdot \right) $ as $n\rightarrow \infty $ in the
case $T\asymp 1.$ Konakov (2006) proved an analogous result for
the case $T=o\left( 1\right).$ Edgeworth type expansions for the
case $T\asymp 1$ and homogenous diffusions were obtained in
Konakov and Mammen (2005). The goal of the present paper is to
obtain the
Edgeworth type expansions for nonhomegenious case and for both cases $%
T\asymp 1$ or $T=o\left( 1\right) .$ The following theorem contains our main
result. It gives Edgeworth type expansions for $p_{h}.$ For the statement of
\ the theorem we introduce the following differential operators

\begin{eqnarray*}
\mathcal{F}_{1}[f](s,t,x,y)&=&\sum_{\left| \nu \right|
=3}\frac{\chi _{\nu }(s,x)}{\nu !}D_{x}^{\nu }f(s,t,x,y),
\\
\mathcal{F}_{2}[f](s,t,x,y)&=&\sum_{\left| \nu \right|
=4}\frac{\chi _{\nu }(s,y)}{\nu !}D_{x}^{\nu }f(s,t,x,y).
\end{eqnarray*} Furthermore, we introduce two terms corresponding
to the classical Edgeworth expansion
\begin{eqnarray}
\widetilde{\pi }_{1}(s,t,x,y)&=&(t-s)\sum_{\left| \nu \right| =3}\frac{%
\overline{\chi }_{\nu }(s,t,y)}{\nu !}D_{x}^{\nu }\widetilde{p}(s,t,x,y),
\label{eq:pi1}
\\
\label{eq:pi2}
\widetilde{\pi }_{2}(s,t,x,y)&=&(t-s)\sum_{\left| \nu \right| =4}\frac{%
\overline{\chi }_{\nu }(s,t,y)}{\nu !}D_{x}^{\nu }\widetilde{p}(s,t,x,y)
\\ && \nonumber
+\frac{1}{2}(t-s)^{2}\left\{ \sum_{\left| \nu \right| =3}\frac{\overline{%
\chi }_{\nu }(s,t,y)}{\nu !}D_{x}^{\nu }\right\}
^{2}\widetilde{p}(s,t,x,y),
\end{eqnarray}
where
\[
\overline{\chi }_{\nu }(s,t,y)=\frac{1}{t-s}\int_{s}^{t}\chi _{\nu }(u,y)du,
\]
and where $\chi _{\nu }(t,x)$ is the $\nu -th$ cumulant of the
density of the
innovations $q(t,x,\cdot )$. The gaussian transition densities $\widetilde{%
p}(s,t,x,y)$\ are defined in (\ref{eq:002}). Note, that in the
homogenous case $\chi _{\nu }(u,y)\equiv \chi _{\nu }(y)$ and
$\overline{\chi }_{\nu }(s,t,y)\equiv \chi _{\nu }(y),$ where
$\chi _{\nu }(y)$ is the $\nu -th$ cumulant of \ the density
$q(y,\cdot ).$

\textbf{Theorem 1. } \emph{Assume (A1)-(A3), (B1),(B2). Then there exists a constant $%
\delta >0$ such that the following expansion holds:}
\[
\sup_{x,y\in R^{d}}\left[ T^{d/2}\left( 1+\left\| \frac{y-x}{\sqrt{T}}%
\right\| ^{S^{\prime }}\right) \times \left|
p_{h}(0,T,x,y)-p(0,T,x,y)\right. \right.
\]
\[
\left. \left. -h^{1/2}\pi _{1}(0,T,x,y)-h\pi _{2}(0,T,x,y)\right| \right]
=O(h^{1+\delta }),
\]
\emph{where $S^\prime$ is defined in Assumption (A3) and where}
\begin{eqnarray*}
\pi _{1}(0,T,x,y)&=&p\otimes \mathcal{F}_{1}[p])(0,T,x,y),
\\
\pi _{2}(0,T,x,y)&=&(p\otimes \mathcal{F}_{2}[p])(0,T,x,y)+p\otimes \mathcal{F}%
_{1}[p\otimes \mathcal{F}_{1}[p]](0,T,x,y)
\\ &&
+\frac{1}{2}p\otimes (L_{\star }^{2}-L^{2})p(0,T,x,y)-\frac{1}{2}p\otimes
(L^{\prime }-\widetilde{L}^{\prime })p(0,T,x,y).
\end{eqnarray*}
\emph{Here $p(s,t,x,y)$ is a transition density of the diffusion
$Y(t)$ and
the operator L}$_{\ast }\emph{\ }$\emph{is defined analogously to $%
\widetilde{L}$ but with the coefficients {}``frozen'' at the point $x.$ The
norm $\left\| \bullet \right\| $ is the usual Euclidean norm.}

\bigskip \textbf{Discussion and remarks}

\textbf{1.} It can be easily shown that
\begin{eqnarray*}
\left| h^{1/2}\pi _{1}(0,T,x,y)\right| &\leq&
C_{1}n^{-1/2}T^{-d/2}\exp \left[ -C_{2}\left\|
\frac{y-x}{\sqrt{T}}\right\| ^{2}\right] ,
\\
\left| h\pi _{2}(0,T,x,y)\right| &\leq& C_{1}n^{-1}T^{-d/2}\exp
\left[ -C_{2}\left\| \frac{y-x}{\sqrt{T}}\right\| ^{2}\right] ,
\end{eqnarray*}
and that by definition of $h$%
\[
\left| O(h^{1+\delta })T^{-d/2}\left[ 1+\left\| \frac{y-x}{\sqrt{T}}\right\|
^{S^{\prime }}\right] ^{-1}\right| \leq C_{1}n^{-1-\delta }T^{1+\delta -d/2}%
\left[ 1+\left\| \frac{y-x}{\sqrt{T}}\right\| ^{S^{\prime }}\right] ^{-1}
\]
with some positive constants $C_{1}$ and $C_{2}.$

\textbf{2.} If the innovation density $q(t,x,\cdot )$ does not
depend on $x$
then $L_{\ast }=L,L^{\prime }=\widetilde{L}^{\prime }$ and $p(s,t,x,y)=%
\widetilde{p}(s,t,x,y)$ \ where $\widetilde{p}$ \ is defined in (\ref{eq:002}%
) with $\sigma (s,t,y)=\sigma (s,t)=\int_{s}^{t}\sigma (u)du$ \ and \ $%
m(s,t,y)=m(s,t)=\int_{s}^{t}m(u)du.$ \ This gives
\begin{eqnarray*}
&&\pi _{1}(0,T,x,y)=\int_{0}^{T}ds\int
\widetilde{p}(0,s,x,v)\sum_{\left| \nu \right| =3}\frac{\chi _{\nu
}(s)}{\nu !}D_{v}^{\nu }\widetilde{p}(s,T,v,y)dv
\\
&&\qquad=-\sum_{\left| \nu \right| =3}\int_{0}^{T}\frac{\chi _{\nu }(s)}{\nu !}%
dsD_{y}^{\nu }\int \widetilde{p}(0,s,x,v)\widetilde{p}(s,T,v,y)dv
\\
&&\qquad=-\sum_{\left| \nu \right| =3}\frac{T}{\nu
!}\overline{\chi }_{\nu }(0,T)D_{y}^{\nu
}\widetilde{p}(0,T,x,y)=\widetilde{\pi }_{1}(0,T,x,y),
\end{eqnarray*}
\begin{eqnarray*}
&&\widetilde{p}\otimes \mathcal{F}_{1}[\widetilde{p}](s,T,z,y)=\int_{s}^{T}du%
\int \widetilde{p}(s,u,z,w)\sum_{\left| \nu \right| =3}\frac{\chi _{\nu }(u)%
}{\nu !}D_{w}^{\nu }\widetilde{p}(u,T,w,y)dw
\\
&&\qquad=-\sum_{\left| \nu \right| =3}D_{y}^{\nu }\int_{s}^{T}\frac{\chi _{\nu }(u)}{%
\nu !}\widetilde{p}(s,T,z,y)=(T-s)\sum_{\left| \nu \right| =3}\frac{%
\overline{\chi }_{\nu }(s,T)}{\nu !}D_{z}^{\nu }\widetilde{p}(s,T,z,y),
\end{eqnarray*}
\begin{eqnarray*}
&&\mathcal{F}_{1}[\widetilde{p}\otimes \mathcal{F}_{1}[\widetilde{p}%
](s,T,z,y)=(T-s)\sum_{\left| \nu \right| =3}\frac{\chi _{\nu }(s)}{\nu !}%
D_{z}^{\nu }\left[ \sum_{\left| \nu ^{\prime }\right| =3}\frac{\overline{%
\chi }_{\nu ^{\prime }}(s,T)}{\nu ^{\prime }!}D_{z}^{\nu ^{\prime }}%
\widetilde{p}(s,T,z,y)\right]
\\
&&\qquad=(T-s)\sum_{\left| \nu \right| =3,\left| \nu ^{\prime
}\right| =3}\frac{\chi _{\nu }(s)}{\nu !}\frac{\overline{\chi
}_{\nu ^{\prime }}(s,T)}{\nu ^{\prime }!}D_{z}^{\nu +\nu ^{\prime
}}\widetilde{p}(s,T,z,y),
\\
&&\widetilde{p}\otimes \mathcal{F}_{2}[\widetilde{p}](0,T,x,y)+\widetilde{p}%
\otimes \mathcal{F}_{1}[\widetilde{p}\otimes \mathcal{F}_{1}[\widetilde{p}%
](0,T,x,y)=T\sum_{\left| \nu \right| =4}\frac{\overline{\chi }_{\nu }(0,T)}{%
\nu !}D_{x}^{\nu }\widetilde{p}(0,T,x,y)
\\ &&
\qquad\qquad+\int_{0}^{T}ds\int
\widetilde{p}(0,s,x,z)(T-s)\sum_{\left| \nu \right|
=3,\left| \nu ^{\prime }\right| =3}\frac{\chi _{\nu }(s)}{\nu !}\frac{%
\overline{\chi }_{\nu ^{\prime }}(s,T)}{\nu ^{\prime }!}D_{y}^{\nu
+\nu ^{\prime }}\widetilde{p}(s,T,z,y)
\\
&&\qquad=T\sum_{\left| \nu \right| =4}\frac{\overline{\chi }_{\nu }(0,T)}{\nu !}%
D_{x}^{\nu }\widetilde{p}(0,T,x,y)\\ && \qquad\qquad+\sum_{\left|
\nu \right| =3,\left| \nu
^{\prime }\right| =3}\frac{1}{\nu !}\frac{1}{\nu ^{\prime }!}%
\int_{0}^{T}\chi _{\nu }(s)\left( \int_{s}^{T}\chi _{\nu ^{\prime
}}(u)du\right) dsD_{x}^{\nu +\nu ^{\prime }}\widetilde{p}(s,T,x,y).
\end{eqnarray*}
For $\nu =\nu ^{\prime }$ \ we have
\[
\int_{0}^{T}\chi _{\nu }(s)\left( \int_{s}^{T}\chi _{\nu ^{\prime
}}(u)du\right) ds=\frac{1}{2}\int_{0}^{T}\int_{0}^{T}\chi _{\nu }(s)\chi
_{\nu }(u)dsdu=\frac{T^{2}}{2}\overline{\chi }_{\nu }(0,T)\overline{\chi }%
_{\nu }(0,T).
\]
For $\ \nu \neq \nu ^{\prime }$ we consider
\begin{eqnarray*}
&&\int_{0}^{T}\chi _{\nu }(s)\left( \int_{s}^{T}\chi _{\nu
^{\prime }}(u)du\right) ds+\int_{0}^{T}\chi _{\nu ^{\prime
}}(s)\left( \int_{s}^{T}\chi _{\nu }(u)du\right) ds
\\
&& \qquad =\int_{0}^{T}\int_{s}^{T}\left[ \chi _{\nu }(s)\chi
_{\nu ^{\prime
}}(u)+\chi _{\nu ^{\prime }}(s)\chi _{\nu }(u)\right] dsdu\\
&& \qquad =\frac{1}{2}%
\int_{0}^{T}\int_{0}^{T}\left[ \chi _{\nu }(s)\chi _{\nu ^{\prime
}}(u)+\chi _{\nu ^{\prime }}(s)\chi _{\nu }(u)\right] dsdu
\\ &&\qquad =\frac{T^{2}}{2}\overline{\chi }_{\nu }(0,T)\overline{\chi }_{\nu
^{\prime
}}(0,T)+\frac{T^{2}}{2}\overline{\chi }_{\nu ^{\prime }}(0,T)\overline{\chi }%
_{\nu }(0,T).
\end{eqnarray*}
From the last equations we obtain
\begin{eqnarray*}
&&
\widetilde{p}\otimes \mathcal{F}_{2}[\widetilde{p}](0,T,x,y)+\widetilde{p}%
\otimes \mathcal{F}_{1}[\widetilde{p}\otimes \mathcal{F}_{1}[\widetilde{p}%
](0,T,x,y)
\\&&\qquad
=T\sum_{\left| \nu \right| =4}\frac{\overline{\chi }_{\nu }(0,T)}{\nu !}%
D_{x}^{\nu }\widetilde{p}(0,T,x,y)+\frac{T^{2}}{2}\left\{ \sum_{\left| \nu
\right| =3}\frac{\overline{\chi }_{\nu }(0,T)}{\nu !}D_{x}^{\nu }\right\}
\widetilde{p}(0,T,x,y)\\
&& \qquad =\widetilde{\pi }_{2}(0,T,x,y).
\end{eqnarray*}
Thus for this case we get the first two terms of the classical
Edgeworth
expansion \ $h^{1/2}\widetilde{\pi }_{1}(0,T,x,y)+h\widetilde{\pi }%
_{2}(0,T,x,y)$ \ for the sums of independent non identically distributed
random vectors.

\textbf{3.} If \ $\chi _{\nu }(t,x)=0$ \ for $\left| \nu \right|
=3$ \ and
for \ $t\in \lbrack 0,T]\times R^{d}$ \ then it holds that $\ \mathcal{F}%
_{1}\equiv 0.$ \ The Theorem 1 holds with
\begin{eqnarray*}
\pi _{1}(0,T,x,y)&=&0,\\
\pi _{2}(0,T,x,y)&=&(p\otimes \mathcal{F}_{2}[p])(0,T,x,y)
\\&&
+\frac{1}{2}p\otimes (L_{\star }^{2}-L^{2})p(0,T,x,y)-\frac{1}{2}p\otimes
(L^{\prime }-\widetilde{L}^{\prime })p(0,T,x,y).
\end{eqnarray*}
If in addition $\chi _{\nu }(t,x)=0$ \ for $\left| \nu \right| =4$ \ then
the first four moments of the innovations \ coincide with the first four
moments of a normal distribution with zero mean and covariance matrix $%
\sigma (t,x).$ In this case we have $\mathcal{F}_{2}=0$ and
\begin{eqnarray*}
\pi _{1}(0,T,x,y)&=&0,
\\
\pi _{2}(0,T,x,y)&=&\frac{1}{2}p\otimes (L_{\star }^{2}-L^{2})p(0,T,x,y)-\frac{%
1}{2}p\otimes (L^{\prime }-\widetilde{L}^{\prime })p(0,T,x,y)
\end{eqnarray*} and the first two terms of the Edgeworth expansion do not
depend on the innovation density. In particular the Edgeworth
expansion for the homogeneous Euler scheme holds with the same
$\pi _{1}$ and $\pi _{2}$ \ as in the last two equations. For the
homogenous case
\begin{eqnarray*}
\pi _{1}(0,T,x,y)&=&0,
\\
\pi _{2}(0,T,x,y)&=&\frac{1}{2}p\otimes (L_{\star
}^{2}-L^{2})p(0,T,x,y).
\end{eqnarray*}
This result for $T=[0,1]$ under a weaker condition on the
diffusion matrix was obtained by Bally and Talay (1996).

\section{\label{sec:Pardiff} Parametrix method for diffusions.}

For any $s\in \lbrack 0,T],$ $x,y\in \Bbb{R}^{d}$ we consider an
additional family of ''frozen'' diffusion processes
\[
d\widetilde{Y}_{t}=m\left( t,y\right) dt+\Lambda \left( t,y\right) dW_{t},\;%
\widetilde{Y}_{s}=x,\;s\leq t\leq T.
\]
Let $\widetilde{p}^{y}\left( s,t,x,\cdot \right) $ be the conditional
density of $\widetilde{Y}_{t},$ given $\widetilde{Y}_{s}=x.$ In the sequel
for any $z$ we shall denote $\widetilde{p}\left( s,t,x,z\right) =\widetilde{p%
}^{z}\left( s,t,x,z\right) ,$ where the variable $z$ acts here twise: as the
arument of the density and as defining quantity of the process $\widetilde{Y}%
_{t}.$

The transition densities $\widetilde{p}$ can be computed explicitly

\begin{eqnarray}
\widetilde{p}\left( s,t,x,y\right) &=&\left( 2\pi \right) ^{-d/2}\left( \det
\sigma \left( s,t,y\right) \right) ^{-1/2}  \nonumber \\
&&\times \exp \left( -\frac{1}{2}\left( y-x-m\left( s,t,y\right) \right)
^{T}\sigma ^{-1}\left( s,t,y\right) \left( y-x-m\left( s,t,y\right) \right)
\right) ,  \label{eq:002}
\end{eqnarray}
where
\[
\sigma \left( s,t,y\right) =\int_{s}^{t}\sigma \left( u,y\right)
du,\;\;\;\;\;m\left( s,t,y\right) =\int_{s}^{t}m\left( u,y\right) du.
\]
Note that the differential operators $\ L$ and $\tilde{L}$ corresponds to
the infinitesimal operators of $Y$ or of the frozen process $\widetilde{Y}%
_{s,x,y}$, respectively, i.e.
\begin{eqnarray*}
Lf(s,t,x,y)&=&\lim_{h\rightarrow 0}h^{-1}\{E[f(s,t,Y(s+h),y)\mid
Y(s)=x]-f(s,t,x,y)\},
\\
\tilde{L}f(s,t,x,y)&=&\lim_{h\rightarrow 0}h^{-1}\{E[f(s,t,\tilde{Y}%
_{s,x,y}(s+h),y)]-f(s,t,x,y)\}.
\end{eqnarray*}
We put
\[
H=(L-\tilde{L})\tilde{p}.
\]
Then
\begin{eqnarray*}
H\left( s,t,x,y\right) &=&\frac{1}{2}\sum_{i,j=1}^{d}\left( \sigma
_{ij}\left( s,x\right) -\sigma _{ij}\left( s,y\right) \right) \frac{\partial
^{2}\widetilde{p}\left( s,t,x,y\right) }{\partial x_{i}\partial x_{j}} \\
&&+\sum_{i,j=1}^{d}\left( m_{i}\left( s,x\right) -m_{i}\left( s,y\right)
\right) \frac{\partial \widetilde{p}\left( s,t,x,y\right) }{\partial x_{i}}.
\end{eqnarray*}
In the following lemmas $k$ - fold convolution of $H$ is denoted by $%
H^{\left( k\right) }.$ \ The following results are taken from Konakov and
Mammen (2000).

\textbf{Lemma 1.} \emph{\label{km-1}Let $0\leq s<t\leq T$. It holds
\[
p(s,t,x,y)=\sum_{r=0}^{\infty}\widetilde{p}\otimes H^{(r)}(s,t,x,y).
\]
}

\textbf{Lemma 2. \label{km-2}} \emph{Let $0\leq s<t\leq T$. There
are constants $C$ and $C_{1}$ such that}
\begin{eqnarray*}
\left|H(s,t,x,y)\right|&\leq& C_{1}\rho^{-1}\phi_{C,\rho}(y-x)
\\
\left| \widetilde{p}\otimes H^{(r)}(s,t,x,y)\right| &\leq&C_{1}^{r+1}\frac{%
\rho ^{r}}{\Gamma (1+\frac{r}{2})}\phi _{C,\rho }(y-x),
\end{eqnarray*}
\emph{where $\rho^{2}=t-s,$
$\phi_{C,\rho}(u)=\rho^{-d}\phi_{C}(u/\rho)$ and
\[
\phi_{C}(u)=\frac{\exp(-C\left\Vert u^{2}\right\Vert )}{\int\exp(-C\left%
\Vert v^{2}\right\Vert dv)}.
\]
}

\section{\label{sec:Parchain} Parametrix method for Markov chains.}

For any $0\leq jh\leq T,$ $x,y\in \Bbb{R}^{d}$ we consider an
additional family of ''frozen'' Markov chains defined for $jh\leq
ih\leq T$ \ as
\begin{equation}
\widetilde{X}_{i+1,h}=\widetilde{X}_{i,h}+m\left( ih,y\right) h+\sqrt{h}%
\widetilde{\xi }_{i+1,h},\;\widetilde{X}_{j,h}=x\in \Bbb{R}^{d},\,\,j\leq
i\leq n,  \label{Froz001}
\end{equation}
where $\widetilde{\xi }_{j+1,h},...,\widetilde{\xi }_{n,h}$ is an innovation
sequence such that the conditional density of $\widetilde{\xi }_{i+1,h}$
given $\widetilde{X}_{i,h}=x_{i},...,\widetilde{X}_{0,h}=x_{0}$ equals to $%
q\left( ih,y,\cdot \right) .$ Let us introduce the infinitesimal operators
corresponding to Markov chains (\ref{eq:001}) and (\ref{Froz001})
respectively,
\begin{eqnarray*}
L_{h}f\left( jh,kh,x,y\right) &=&h^{-1}\left( \int p_{h}\left(
jh,\left( j+1\right) h,x,z\right) f\left( \left( j+1\right)
h,kh,z,y\right) dz\right . \\ && \left .-f\left( \left( j+1\right)
h,kh,z,y\right) \right),\\
\widetilde{L}_{h}f\left( jh,kh,x,y\right)&=&h^{-1}\left( \int \widetilde{p}%
_{h}^{y}\left( jh,\left( j+1\right) h,x,z\right) f\left( \left(
j+1\right) h,kh,z,y\right) dz\right . \\ && \left .-f\left( \left(
j+1\right) h,kh,z,y\right) \right) ,
\end{eqnarray*}
where $\widetilde{p}_{h}^{y}\left( jh,j^{\prime }h,x,\cdot \right) $ denotes
the conditional density of $\widetilde{X}_{j^{\prime },h}$ given $\widetilde{%
X}_{j,h}=x.$ As before for any $z$ denote $\widetilde{p}_{h}\left(
jh,j^{\prime }h,x,z\right) =\widetilde{p}_{h}^{z}\left( jh,j^{\prime
}h,x,z\right) ,$ where the variable $z$ acts here twise: as the arument of
the density and as defining quantity of the process $\widetilde{X}_{i,h}.$
For technical convenience the terms $f\left( \left( j+1\right)
h,kh,z,y\right) $ on the right hand side of $L_{h}f$ and $\widetilde{L}_{h}f$
appear instead of $f\left( jh,kh,z,y\right) .$

In analogy with the definition of $H$ we put, for $k>j,$%
\[
H_{h}\left( jh,kh,x,y\right) =\left( L_{h}-\widetilde{L}_{h}\right)
\widetilde{p}_{h}\left( jh,kh,x,y\right) .
\]
We also shall use the convolution type binary operation $\otimes _{h}:$%
\[
g\otimes _{h}f\left( jh,kh,x,y\right) =\sum_{i=j}^{k-1}h\int_{\Bbb{R}%
^{d}}g\left( jh,ih,x,z\right) f\left( ih,kh,z,y\right) dz,
\]
where $0\leq j<k\leq n.$ Write $g\otimes _{h}H_{h}^{\left( 0\right) }=g$ and
$g\otimes _{h}H_{h}^{\left( r\right) }=\left( g\otimes _{h}H_{h}^{\left(
r-1\right) }\right) \otimes _{h}H_{h}$ for $r=1,...,n.$ For the higher order
convolutions we use the convention $\sum_{i=j}^{l}=0$ for $l<j.$ One can
show the following analog of the ''parametrix '' expansion for $p_{h}$ {[}%
see Konakov and Mammen (2000){]}.

\textbf{Lemma 3.} \emph{\label{km-3} Let $0\leq jh<kh\leq T.$ It holds
\[
p_{h}(jh,kh,x,y)=\sum_{r=0}^{k-j}\widetilde{p}_{h}%
\otimes_{h}H_{h}^{(r)}(jh,kh,x,y),
\]
}

\emph{where
\[
p_{h}(jh,jh,x,y)=p_{h}(kh,kh,x,y)=\delta(y-x)
\]
and $\delta$ is the Dirac delta symbol.}

\section{\label{sec:bonedg} Bounds on $\widetilde{p}_{h}-\widetilde{p%
}$ based on Edgeworth expansions.}

In this subsection we will develop some tools that are helpful for
the comparison of the expansion of $p$ (see Lemma 1) and the
expansion of $p_{h}$ ( see Lemma 3). These expansions are simple
expressions in $\tilde{p}$ or $\tilde{p}_{h}$, respectively. Recall that $%
\tilde{p}$ is a Gaussian density, see (\ref{eq:002}), and that $\tilde{p}%
_{h} $ is the density of a sum of independent variables. The densities $%
\tilde{p}$ and $\tilde{p}_{h}$ can be compared by application of the
classical Edgeworth expansions. This is done in Lemma 5. This is the
essential step for the comparison of the expansions of $p$ and $p_{h}$. For
the proof of Lemma 5 we need one additional lemma which is a technical tool
for the further considerations. To formulate this lemma we need some
additional notations. Suppose $X\in \Bbb{R}^{d}$ be a random vector \ having
a density $q(\mathbf{x}),\mathbf{x\in }\Bbb{R}^{d},$ $EX=0,$$Cov(X,X)=\Sigma
,$ where $\Sigma $ be a positively definite $d\times d$ matrix . Denote $%
A=\left\| a_{ij}\right\| =\Sigma ^{-1/2}$ and \ let $\chi _{\nu }(Z)$ be a
cumulant of the order $\nu =(\nu _{1},...,\nu _{d})$ of a random vector $%
Z\in \Bbb{R}^{d} $, $\phi (x)$ denotes a function in $\Bbb{R}^{d}$ such that
$D_{x}^{\nu }\phi (x)$ \ exist and continuous for $\left| \nu \right| =4$,
and $A^{-1}=\left\| a^{ij}\right\| =\Sigma ^{1/2}.$

\textbf{Lemma 4.} \label{km-4}The following relation holds for $s=3$ and for
$s=4$%
\[
\sum_{\left|\nu\right|=s}\frac{\chi_{\nu}(AX)D_{z}^{\nu}\phi(z)}{\nu!}%
=\sum_{\left|\nu\right|=s}\frac{\chi_{\nu}(X)D_{x}^{\nu}\phi(Ax)}{\nu!}
\]
where $z=Ax.$

\emph{\textsc{Proof of Lemma 4.}} For $\left| \nu \right| =3,\nu =(\nu
_{1},...,\nu _{d}),$ each cumulant $\chi _{\nu }(AX)$ is a linear
combination of $\chi _{\mu }(X)$ with $\left| \mu \right| =3$ and with
coefficients depending only on $a_{ij}$. It follows from the following
relation
\[
\chi _{\nu }(AX)=\mu _{\nu }(AX)=\int (a_{11}x_{1}+...+a_{1d}x_{d})^{\nu
_{1}}\times ...\times (a_{11}x_{1}+...+a_{1d}x_{d})^{\nu _{d}}q(\mathbf{x})d%
\mathbf{x}.
\]
Analogously, from the usual differentiation rule of a composite function and
from the relation $\phi (z)=\phi (Ax)$ , $x=A^{-1}z$, it follows that $%
D_{z}^{\nu }\phi (z)=D_{z}^{\nu }\phi (Ax)$ is a linear combination of $%
D_{x}^{\nu }\phi (Ax)$ with coefficients depending only on $a^{ij}$. As a
result of such substitutions we obtain that
\begin{eqnarray*}
&&\sum_{\left| \nu \right| =3}\frac{\chi _{\nu }(AX)D_{z}^{\nu }\phi (z)}{\nu !%
}=\frac{1}{3!}\sum_{j=1}^{d}\left[ \sum_{\left| \mu \right|
=3}\frac{3!}{\mu _{1}!...\mu _{d}!}a_{1j}^{\mu _{1}}...a_{dj}^{\mu
_{d}}\chi _{\mu }(X)\right] \\&&\qquad \times \left[ \sum_{\left|
\mu ^{\prime }\right| =3}\frac{3!}{\mu _{1}^{\prime }!...\mu
_{d}^{\prime }!}(a^{j1})^{\mu _{1}^{\prime }}...(a^{jd})^{\mu
_{d}^{\prime }}D_{x}^{\mu ^{\prime }}\phi (AX)\right]
\\&&\qquad
+\frac{1}{2!1!}\sum_{\{i\neq ,j\}}\left[
\sum_{l=1}^{d}\sum_{\left| \mu \right| =2}\frac{2!}{\mu
_{1}!...\mu _{d}!}a_{1j}^{\mu _{1}}...a_{dj}^{\mu _{d}}a_{il}\chi
_{\mu +e_{l}}(X)\right]
\\&&\qquad
\times \left[ \sum_{l^{\prime }=1}^{d}\sum_{\left| \mu ^{\prime }\right| =2}%
\frac{2!}{\mu _{1}^{\prime }!...\mu _{d}^{\prime }!}(a^{j1})^{\mu
_{1}^{\prime }}...(a^{jd})^{\mu _{d}^{\prime }}a^{il^{\prime }}D_{x}^{\mu
^{\prime }+e_{l^{\prime }}}\phi (AX)\right]
\\&&\qquad
+\frac{1}{3!}\sum_{\{i\neq j\neq k\}}\left[
\sum_{l,q=1}^{d}\sum_{\left| \mu \right| =1}\frac{1}{\mu
_{1}!...\mu _{d}!}a_{1j}^{\mu _{1}}...a_{dj}^{\mu
_{d}}a_{il}a_{kq}\chi _{\mu +e_{l}+e_{q}}(X)\right]
\end{eqnarray*} \begin{eqnarray*}
\times \left[ \sum_{l^{\prime },q^{\prime }=1}^{d}\sum_{\left| \mu
^{\prime }\right| =1}\frac{1}{\mu _{1}^{\prime }!...\mu
_{d}^{\prime }!}(a^{j1})^{\mu _{1}^{\prime }}...(a^{jd})^{\mu
_{d}^{\prime }}a^{il^{\prime }}a^{kq^{\prime }}D_{x}^{\mu ^{\prime
}+e_{l^{\prime }}+e_{q^{\prime }}}\phi (AX)\right]
\end{eqnarray*}
where $\sum_{\{i\neq ,j\}}$ ( $\sum_{\{i\neq j\neq k\}}$) denotes
the sum over all different pairs (triples) of $i,j\in
\{1,2,...,d\}$ ( of $i,j,k\in \{1,2,...,d\}$) and $e_{i}\in
\Bbb{R}^{d}$ denotes the vector whose $i-$th coordinate is equal
to 1 and other coordinates are zero. Collecting the similar terms
in the last equation we obtain that for $\nu =3e_{k}$, $\nu
^{\prime }=3e_{l}$ the coefficient before $\chi _{\nu
}(X)D_{x}^{\nu
^{\prime }}\phi (AX)$ is equal to $\frac{1}{3!}%
(a_{1k}a^{l1}+...+a_{dk}a^{ld})^{3}=\frac{1}{3!}\delta _{kl}$, for $\nu
=e_{q}+2e_{r},\text{$\nu ^{\prime }=e_{l}+2e_{n}$ , }q\neq r,$ the
coefficient before $\chi _{\nu }(X)D_{x}^{\nu ^{\prime }}\phi (AX)$ \ is
equal to $\frac{1}{2!}%
(a_{1q}a^{l1}+...+a_{dq}a^{ld})(a_{1r}a^{n1}+...+a_{dr}a^{nd})^{2}=
\frac{1}{2!}\delta _{ql}\delta _{rn}$, \ in particular , for $l=n$
\ the last expression is equal to zero. For \ \ $\ \nu
=e_{q}+e_{r}+e_{n},\nu ^{\prime }=e_{q^{\prime }}+e_{r^{\prime
}}+e_{n^{\prime }}$ \ $q\neq r,q\neq n,r\neq n,$ the coefficient before $%
\chi _{\nu }(X)D_{x}^{\nu ^{\prime }}\phi (AX)$ \ is equal to $%
(a_{1q}a^{q^{\prime }1}+...+a_{dq}a^{q^{\prime }d})\times
(a_{1r}a^{r^{\prime }1}+...+a_{dr}a^{r^{\prime }d})\times $ $%
(a_{1n}a^{n^{\prime }1}+...+a_{dn}a^{n^{\prime }d})=\delta _{qq^{\prime
}}\delta _{rr^{\prime }}\delta _{nn^{\prime }}$ $.$ This proves lemma for $%
\left| \nu \right| =3.$ The proof for $\left| \nu \right| =4$ is quite
similar. For this case we use the relation which enabes to express a
cumulant $\chi _{\nu }(AX)$ as $\mu _{\nu }(AX)$ plus a second order
polinomial of the moments $\mu _{\nu ^{\prime }}(AX)$, $\left| \nu ^{\prime
}\right| =2.$ A necessary correction term for $\mu _{\nu }(X)$ to get a $%
\chi _{\nu }(X)$ comes from the derivation of $D_{z}^{\nu }\phi (z)$. This
completes the proof of the lemma.

In Lemma 6 bounds will be given for derivatives of $\ \widetilde{p}_{h}$.
The proof of this lemma also makes essential use of \ Edgeworth expansions.
The following lemma is a higher order extension of the results in Section
3.3 in Konakov and Mammen (2000). Denote
\begin{equation}
\mu _{j,k}(y)=h\sum_{i=j}^{k-1}m(ih,y),V_{j,k}(y)=h\sum_{i=j}^{k-1}\sigma
(ih,y).  \label{eq:003a}
\end{equation}

\textbf{Lemma 5.} \emph{\label{km-5} The following bound holds with a
constant $C$ for $\nu =(\nu _{1},...\nu _{p})^{T}$ with $0\leq \left| \nu
\right| \leq 6$}
\[
\left| D_{z}^{\nu }\widetilde{p}_{h}(jh,kh,x,y)-D_{z}^{\nu }\widetilde{p}%
(jh,kh,x,y)\right.
\]
\[
\left. -\sqrt{h}D_{z}^{\nu }\widetilde{\pi }_{1}(jh,kh,x,y)-hD_{z}^{\nu }%
\widetilde{\pi }_{2}(jh,kh,x,y)\right|
\]
\[
\leq Ch^{3/2}\rho ^{-3}\zeta _{\rho }^{S-\left| \nu \right| }(y-x)
\]
\emph{for all $j<k,x$ and $y$. Here $D_{z}^{\nu }$ denotes the partial
differential operator of order $\nu $ with respect to $%
z=V_{j,k}^{-1/2}(y)(y-x-\mu _{j,k}(y))$. The quantity $\rho $ denotes again
the term $\rho =[h(k-j)]^{1/2}$ }$and$\emph{\ $the\ functions$ \ }$%
\widetilde{\pi }_{1}$\emph{\ and \ }$\widetilde{\pi }_{2}$ $\ are$ $\
defined $ $\ in$ $\ (\ref{eq:pi1})$ \ $and$ \ $(\ref{eq:pi2}).$ \emph{We
write $\zeta _{\rho }^{k}(\cdot )=\rho ^{-p}\zeta ^{k}(\cdot /\rho )$ where
\[
\zeta ^{k}(z)=\frac{[1+\left\| z\right\| ^{k}]^{-1}}{\int [1+\left\|
z^{\prime }\right\| ^{k}]^{-1}dz^{\prime }}.
\]
}

\textsc{Proof of Lemma 5.} We note first that $\tilde{p}_{h}(jh,kh,x,\bullet
)$ is the density of the vector
\[
x+\mu _{j,k}(y)+h^{1/2}\sum_{i=j}^{k-1}\widetilde{\xi }_{i+1,h},
\]
where, as above in the definition of the {}``frozen'' Markov chain $\tilde{Y}%
_{n}$, $\widetilde{\xi }_{i+1,h}$ is a sequence of independent variables
with densities $q(ih,y,\cdot ),$ $\mu _{j,k}(y)=\sum_{i=j}^{k-1}hm(ih,y).$
Let $f_{h}(\cdot )$ be the density of the normalized sum
\[
h^{1/2}\left[ V_{j,k}(y)\right] ^{-1/2}\sum_{i=j}^{k-1}\widetilde{\xi }%
_{i+1,h}.
\]
Clearly, we have
\[
\tilde{p}_{h}(jh,kh,x,\cdot )=\det \left[ V_{j,k}(y)\right] ^{-1/2}f_{n}\{%
\left[ V_{j,k}(y)\right] ^{-1/2}[\cdot -x-\mu _{j,k}(y)]\}.
\]
We now argue that an Edgeworth expansion holds for $f_{h}$. This implies the
following expansion for $\tilde{p}_{h}(jh,kh,x,\cdot )$%
\begin{equation}
\tilde{p}_{h}(jh,kh,x,\cdot )  \label{eq:10}
\end{equation}
\[
=\det \left[ V_{j,k}(y)\right] ^{-1/2}[\sum_{r=0}^{S-3}(k-j)^{-r/2}P_{r}(-%
\phi :\{\bar{\chi}_{\beta ,r}\})\{\left[ V_{j,k}(y)\right] ^{-1/2}[\cdot
-x-\mu _{j,k}(y)]\}
\]
\[
+[k-j]^{-(S-2)/2}O([1+\left\| \{[V_{j,k}(y)]^{-1/2}[\cdot -x-\mu
_{j,k}(y)]\}\right\| ^{S}]^{-1})]
\]
with standard notations, see Bhattacharya and Rao (1976), p. 53. In
particular, $P_{r}$ denotes a product of a standard normal density with a
polynomial that has coefficients depending only on cumulants of order $\leq
r+2$. Expansion (\ref{eq:10}) follows from Theorem 19.3 in Bhattacharya and
Rao (1976). This can be seen as in the proof of Lemma 3.7 in Konakov and
Mammen (2000a).

It follows from (\ref{eq:10}) and Condition (A3) that
\[
\left| \tilde{p}_{h}(jh,kh,x,y)-\tilde{p}(jh,kh,x,y)\right.
\left. -h^{1/2}\widehat{\pi }_{1}(jh,kh,x,y)-h\widehat{\pi }%
_{2}(jh,kh,x,y)\right|
\]
\begin{equation}
\leq Ch^{3/2}\rho ^{-3}\zeta _{\rho }^{S-|\nu |}(y-x),  \label{eq:11}
\end{equation}
where
\begin{eqnarray*}
&&\tilde{p}(jh,kh,x,y)=\det \left[ V_{j,k}(y)\right] ^{-1/2}(2\pi
)^{-p/2}
\\&& \qquad
\exp \{-\frac{1}{2}(y-x-\mu _{j,k}(y))^{T}\left[ V_{j,k}(y)\right]
^{-1}(y-x-\mu _{j,k}(y))\},
\\&&
\widehat{\pi }_{1}(jh,kh,x,y)=-\rho ^{-1}\det \left[
V_{j,k}(y)\right] ^{-1/2}\sum_{\left| \nu \right|
=3}\frac{\overline{\chi }_{\nu ,j,k}(y)}{\nu !}D_{z}^{\nu }\phi
\left\{ \left[ V_{j,k}(y)\right] ^{-1/2}(y-x-\mu
_{j,k}(y))\right\} ,
\\&&
\widehat{\pi }_{2}(jh,kh,x,y)=\rho ^{-2}\det \left[ V_{j,k}(y)\right] ^{-1/2}%
\left[ \sum_{\left| \nu \right| =4}\frac{\overline{\chi }_{\nu ,j,k}(y)}{\nu
!}D_{z}^{\nu }\phi \left\{ \left[ V_{j,k}(y)\right] ^{-1/2}(y-x-\mu
_{j,k}(y))\right\} \right.
\\&& \qquad
\left. +\frac{1}{2}\left\{ \sum_{|\nu |=3}\frac{\overline{\chi
}_{\nu
,j,k}(y)}{\nu !}D_{z}^{\nu }\right\} ^{2}\phi \left\{ \left[ V_{j,k}(y)%
\right] ^{-1/2}(y-x-\mu _{j,k}(y))\right\} \right] ,
\end{eqnarray*}
where $\ \overline{\chi }_{\nu ,j,k}(y)=\frac{1}{k-j}\sum_{i=j}^{k-1}\chi
_{\nu ,j,k,i}(y),$ $\chi _{\nu ,j,k,i}(y)=\nu -$th cumulant of $\ \rho \left[
V_{j,k}(y)\right] ^{-1/2}\widetilde{\xi }_{i+1,h}=\rho ^{\left| \nu \right|
}\times \{\nu -$th cumulant of $\ \left[ V_{j,k}(y)\right] ^{-1/2}\widetilde{%
\xi }_{i+1,h}\},$ and $D_{z}^{\nu }\phi (z)$ denotes the $\nu -$th \
derivative of $\phi $ with respect to $z=\left[ V_{j,k}(y)\right] ^{-1/2}$ $%
(y-x-\mu _{j,k}(y))$\ . It follows from the (conditional) independence of \ $%
\widetilde{\xi }_{i+1,h},i=j,...,k-1,$ that $\ \overline{\chi }_{\nu
,j,k}(y)=\frac{\rho ^{\left| \nu \right| }}{k-j}h^{-\left| \nu \right|
/2}\times \chi _{\nu }(AX),$ where $A=h^{1/2}\left[ V_{j,k}(y)\right]
^{-1/2}=\Sigma ^{-1/2},\Sigma =Cov(X,X),X=\sum_{i=j}^{k-1}\widetilde{\xi }%
_{i+1,h}.$ By Lemma 4 for $s=3,4$
\begin{eqnarray} \nonumber
\sum_{\left| \nu \right| =s}\frac{\overline{\chi }_{\nu ,j,k}(y)}{\nu !}%
D_{z}^{\nu }\phi (z)&=&\rho ^{s}\frac{1}{k-j}\sum_{\left| \nu \right| =s}\frac{%
\chi _{\nu }(AX)}{\nu !}D_{h^{1/2}z}^{\nu }\phi _{h}(h^{1/2}z)
\\ \nonumber &=&(-1)^{s}\rho ^{s}\sum_{\left| \nu \right| =s}\frac{\overline{\chi
}_{\nu }(X)}{\nu !}D_{x}^{\nu }\phi _{h}(A(y-x-\mu _{j,k}(y)))
\\&=&(-1)^{s}\rho ^{s}\sum_{\left| \nu \right| =s}\frac{\overline{\chi
}_{\nu }(X)}{\nu !}D_{x}^{\nu }\phi (\left[ V_{j,k}(y)\right]
^{-1/2}(y-x-\mu _{j,k}(y))),  \label{eq:12}
\end{eqnarray}
where we put $\phi _{h}(z)=\phi (h^{-1/2}z),\overline{\chi }_{\nu }(X)=\frac{%
1}{k-j}\sum_{i=j}^{k-1}\chi _{\nu }(ih,y).$ It follows from (\ref{eq:12})
and the condition \textbf{B1 }that up to the error term in the right hand
side of \ (\ref{eq:11}) the \ functions $\widehat{\pi }_{1}$ and $\widehat{%
\pi }_{2}$ coincide with the functions $\widetilde{\pi }_{1}$ and \ $%
\widetilde{\pi }_{2}$ given at the beginning of Section 4. For $\nu =0$ the
statement of the lemma immediately follows from (\ref{eq:11}). For $\nu >0$
one proceeds similarly. See the remark at the end of the proof of Lemma 3.7
in Konakov and Mammen (2000).

The following lemma was proved in Konakov and Molchanov (1984) (Lemma 4 on
page 68).

\textbf{Lemma 6. }\textit{Let }$L(d)$ \textit{be the set of symmetric
matrices, and let }$D_{\lambda ^{+},\lambda ^{-}}\subset L(d),0<\lambda
^{-}<\lambda ^{+}<\infty ,\mathit{\ }$\textit{be the open subset
distinguished by the inequalities \ }$\Lambda \in D_{\lambda ^{+},\lambda
^{-}}$ $\Leftrightarrow $ $\lambda ^{-}I\leq \Lambda \leq \lambda ^{+}I;$ $%
A=A(\Lambda )$ \ \textit{is a} \textit{solution of the} \textit{equation }$%
A^{2}=\Lambda $ \ with \ $A=A^{T\text{ \ }}$and \ $a_{ij}(\Lambda )$ \textit{%
\ the corresponding matrix elements. \ Then for any \ }$k,l\leq d,i,j\leq d$
\ \textit{and \ }$\Lambda \subset D_{\lambda ^{+},\lambda ^{-}}$ \ \textit{%
we have that}
\[
\left| \frac{\partial a_{ij}(\Lambda )}{\partial \lambda _{kl}}\right| \leq
\frac{1}{2\sqrt{\lambda ^{-}}}.
\]

From Lemma 5 we get the following corollary. The statement of \ the next
lemma is an extension of \ Lemma 3.7 in Mammen and Konakov (2000) where the
result has been shown for $0\leq |b|\leq 2,a=0$.

\textbf{Lemma 7. }\textit{The following bounds hold:}
\[
\left| D_{y}^{a}D_{x}^{b}\widetilde{p}_{h}(jh,kh,x,y)\right| \leq C\rho
^{-\left| a\right| -\left| b\right| }\varsigma _{\rho }^{S-\left| a\right|
}(y-x)
\]
\textit{for all }$j<k,$ \textit{for all }$x$ \ and \ $y$ \textit{and for all
}$a,b$ \ \textit{with \ }$0\leq \left| a\right| +\left| b\right| \leq 6.$
\textit{Here, }$\rho =[(k-j)h]^{1/2}$ . \textit{The exponent }$S$ \ \textit{%
has been defined in Assumption }\textbf{A3. }

\textsc{Proof of Lemma 7.} \ If $\ A=\left\| a_{ij}\right\| $ and $B=\left\|
b_{kl}\right\| $ \ and elements $a_{ij}(B)$ \ are smooth functions of $%
b_{kl} $ \ then an inequality $\left| \frac{\partial A}{\partial B}\right|
\leq C$ will mean that \ $\left| \frac{\partial a_{ij}}{\partial b_{kl}}%
\right| \leq C$ for all $1\leq i,j\leq d,1\leq k,l\leq d.$ To obtain the
assertion of the lemma we have to estimate the derivatives $%
D_{y}^{a}D_{x}^{b}z,$where \ \emph{$z=V_{j,k}^{-1/2}(y)(y-x-\mu _{j,k}(y)).$
} Note that \emph{$z=z(V_{j,k}^{-1/2},\mu _{j,k},x,y),$ }where $%
V_{j,k}^{-1/2}=V_{j,k}^{-1/2}(y) $ and \emph{\ }$\mu _{j,k}=\mu _{j,k}(y).$
It follows from the conditions (B1) and (\ref{eq:003a}) that
\begin{equation}
\left| \frac{\partial \mu _{j,k}(y)}{\partial y}\right| \leq C\rho
^{2},\left| \frac{\partial V_{j,k}(y)}{\partial y}\right| \leq C\rho ^{2}.
\label{eq:12a}
\end{equation}
It follows from Lemma 6 that
\begin{equation}
\left| \frac{\partial V_{j,k}^{1/2}}{\partial V_{j,k}}\right| \leq C.
\label{eq:12b}
\end{equation}
From inequalities (3.16) \ in Konakov and Mammen (2000) and from the
representation of an inverse matrix in terms of cofactors divided by
determinant we obtain that
\begin{equation}
\left| \frac{\partial V_{j,k}^{-1/2}}{\partial V_{j,k}^{1/2}}\right| \leq
C\rho ^{-2}.  \label{eq:12c}
\end{equation}

\bigskip From (\ref{eq:12a})-(\ref{eq:12c}) and from the differentiation
rule of a composite function we get
\begin{equation}
\left| \frac{\partial V_{j,k}^{-1/2}(y)}{\partial y}\right| \leq C.
\label{eq:12d}
\end{equation}

Inequalities (3.16) in Konakov and Mammen (2000), (\ref{eq:12d})
and the differentiation rule of a composite function imply
\begin{equation}
\left| \frac{\partial z}{\partial \mu _{j,k}}\right| \leq \frac{C}{\rho }%
,\left| \frac{\partial z}{\partial V_{j,k}^{-1/2}}\right| =(y-x-\mu
_{j,k}(y)),\left| \frac{\partial z}{\partial y}\right| \leq \frac{C}{\rho }.
\label{eq:12e}
\end{equation}
Inequalities (3.16) in Konakov and Mammen (2000) also imply
\begin{equation}
\left| \frac{\partial z}{\partial x}\right| \leq \frac{C}{\rho }.
\label{eq:12f}
\end{equation}
The assertion of \ Lemma 7 for \ $a=e_{i},b=e_{j},1\leq i,j\leq
d,$ follows from Lemma 5 and from (\ref{eq:12e}), (\ref{eq:12f}).
For other values of $a $ and $b$ one has to repeat this arguments.

\section{\label{sec:bonker} Bounds on operator kernels used in the
parametrix expansions.}

In this section we will present bounds for operator kernels appearing in the
expansions based on the parametrix method. In Lemma 8 we compare the
infinitesimal operators $L_{h}$ and $\tilde{L}_{h}$ with the differential
operators $L$ and $\tilde{L}$. We give an approximation for the error if, in
the definition of $H_{h}=(L_{h}-\tilde{L}_{h})\tilde{p}_{h}$, the terms $%
L_{h}$ and $\tilde{L}_{h}$ are replaced by $L$ or $\tilde{L}$, respectively.
We show that this term can be approximated by $K_{h}+M_{h}$, where $K_{h}=(L-%
\tilde{L})\tilde{p}_{h}$ and where $M_{h}$ is defined in Remark 1 after
Lemma 8 . Bounds on $H_{h}$, $K_{h}$, and $M_{h}$ are given in Lemma 9.
These bounds will be used in the proof of our theorem to show that in the
expansion of $p_{h}$ the terms $\tilde{p}_{h}\otimes _{h}H_{h}^{(r)}$ can be
replaced by $\tilde{p}_{h}\otimes _{h}(K_{h}+M_{h}^{\prime })^{(r)}$, $%
M_{h}^{\prime }$ is defined in Lemma \ 9.

\textbf{Lemma 8.} \emph{The following bound holds with a constant
$C$}
\begin{eqnarray*}
&&\left| H_{h}(jh,kh,x,y)-K_{h}^{\prime }(jh,kh,x,y)-M_{h}^{\prime
}(jh,kh,x,y)-R_{h}(jh,kh,x,y)\right| \\&& \qquad \qquad \leq
Ch^{3/2}\rho ^{-1}\zeta _{\rho }^{S}(y-x)
\end{eqnarray*}
\emph{with $\zeta _{\rho }^{S}$ as in Lemma 5 \ for all $j<k$, $x$
and $y$ . For $j < k-1$ we define }
\begin{eqnarray*}
&&K_{h}^{\prime }(jh,kh,x,y)=(L-\widetilde{L})\lambda
(x),M_{h}^{\prime }(jh,kh,x,y)\\&& \qquad =M_{h,1}(jh,kh,x,y)
+M_{h,2}(jh,kh,x,y)+M_{h,3}^{\prime }(jh,kh,x,y)
\\&&
M_{h,1}(jh,kh,x,y)=h^{1/2}\sum_{\left| \nu \right|
=3}\frac{D_{x}^{\nu }\lambda (x)}{\nu !}(\chi _{\nu }(jh,x)-\chi
_{\nu }(jh,y)),
\\&&
M_{h,2}(jh,kh,x,y)=h\sum_{\left| \nu \right| =4}\frac{D_{x}^{\nu }\lambda (x)%
}{\nu !}(\chi _{\nu }(jh,x)-\chi _{\nu }(jh,y)),
\\&&
M_{h,3}^{\prime }(jh,kh,x,y)=\frac{h}{2}(L_{\ast }^{2}-\widetilde{L}%
^{2})\lambda (x),
\\&&
R_{h}(jh,kh,x,y)=h^{3/2}\sum_{\left| \nu \right|
=4}\frac{D_{x}^{\nu }\lambda (x)}{\nu !}\sum_{r=1}^{d}\nu
_{r}[m_{r}(jh,x)\mu _{\nu -e_{r}}(jh,x)-m_{r}(jh,y)\mu _{\nu
-e_{r}}(jh,y)]
\\&& \qquad
+5\sum_{\left| \nu \right| =5}\frac{1}{\nu !}%
\sum_{k=1}^{d}(m_{k}(jh,x)-m_{k}(jh,y))\left\{ \nu _{k}\int q(jh,x,\theta )%
\widetilde{h}^{\nu -e_{k}}(\theta )\right.
\\&& \qquad \qquad
\left. \times \left[ \int_{0}^{1}(1-u)^{4}D^{\nu }\lambda (x+u\widetilde{h}%
(\theta ))du\right] d\theta \right.\\ && \left.\qquad  +\int
q(jh,x,\theta )\widetilde{h}^{\nu }(\theta )\left[
\int_{0}^{1}(1-u)^{4}uD^{\nu +e_{k}}\lambda
(x+u\widetilde{h}(\theta ))du\right] d\theta \right\}
\\&& \qquad
+h^{2}\sum_{\left| \nu \right| =4}\frac{D_{x}^{\nu }\lambda (x)}{\nu !}%
\sum_{\left| \nu ^{\prime }\right| =2}\nu !N(\nu ,\nu ^{\prime })[m^{\nu
^{\prime }}(jh,x)\mu _{\nu -\nu ^{\prime }}(jh,x)-m^{\nu ^{\prime
}}(jh,y)\mu _{\nu -\nu ^{\prime }}(jh,y)].
\end{eqnarray*}
\emph{Here }$L_{\ast }$\emph{\ is defined analogously to }$%
\widetilde{L}$\emph{\ \ but with the coefficients ''frozen'' at the point x,
e}$_{r}$\emph{\ \ denotes a $p$ - dimensional vector with $r-th$ element
equal to $1$ and with all other elements equal to $0.$ Furthermore, for $%
\left| \nu \right| =4,\left| \nu ^{\prime }\right| =2$ we define}
\[
N(\nu ,\nu ^{\prime })=2^{\chi \lbrack \nu ^{\prime }!=1]+\chi \lbrack (\nu
-\nu ^{\prime })!=1]-2}
\]
\emph{where $\chi (\cdot )$ means an indicator function. We put
$m(x)^{\nu }=m_{1}(x)^{\nu _{1}}\cdot ...\cdot m_{p}(x)^{\nu
_{p}}$ and $m(x)^{\nu }=0,$ }$\nu !=0$ \emph{\ and $\mu _{\nu
}(x)=0$ if at least one of the coordinates of $\nu $ is negative.
We define also the following functions}
\begin{eqnarray*}
\lambda (x)&=&\widetilde{p}_{h}((j+1)h,kh,x,y),
\\
\widetilde{h}(\theta )&=&m(jh,y)h+\theta h^{1/2}.
\end{eqnarray*}
\emph{Here again $\rho $ denotes the term $\rho =\left[
h(k-j)\right] ^{1/2}. $ For $j=k-1$ and $l=1,...,3$ we define
\[
K_{h}^{\prime
}(jh,kh,x,y)=M_{h,1}(jh,kh,x,y)=M_{h,2}(jh,kh,x,y)=M_{h,3}^{\prime
}(jh,kh,x,y)=0.
\]
}\textsc{Proof of Lemma 8. }As in the proof of Lemma 3.9 in
Konakov and Mammen (2000) we have
\[
H_{h}(jh,kh,x,y)=H_{h}^{1}(jh,kh,x,y)-H_{h}^{2}(jh,kh,x,y),
\]
\bigskip where
\begin{eqnarray}
H_{h}^{1}(jh,kh,x,y)&=&h^{-1}\int q(jh,x,\theta )[\lambda
(x+h(\theta ))-\lambda (x)]d\theta  \label{eq:12g}
\\
H_{h}^{2}(jh,kh,x,y)&=&h^{-1}\int q(jh,x,\theta )[\lambda (x+\widetilde{h}%
(\theta ))-\lambda (x)]d\theta ,  \label{eq:12h}
\\ \nonumber
h(\theta )&=&m(jh,x)h+\theta h^{1/2},
\\ \nonumber
\widetilde{h}(\theta )&=&m(jh,y)h+\theta h^{1/2}.
\end{eqnarray}
For \ $[\lambda (x+h(\theta ))-\lambda (x)]$ \ and \ $[\lambda (x+\widetilde{%
h}(\theta ))-\lambda (x)]$ \ in (\ref{eq:12g}), (\ref{eq:12h}) we use now
the Tailor expansion up to the order 5\ with the remaining term in integral
form. To pass from \ moments to cumulants we use \ the well known relations
(see e.g. relation (6.11) on page 46 in Bhattacharya and Rao (1986)). After
long but simple calculations we come to the conclusion of the lemma.

\textbf{Remark 1. }We show now that the function $K_{h}^{\prime
}(jh,kh,x,y)+M_{h,3}^{\prime }(jh,kh,x,y)$ \ in Lemma 8 \ is equal
\ to $\
K_{h}(jh,kh,x,y)+\frac{h}{2}(L_{\ast }^{2}-2L\widetilde{L}+\widetilde{L}%
^{2})\lambda (x)+M_{h,3}^{\prime \prime }(jh,kh,x,y)$ \ where
\begin{eqnarray} \nonumber
&&M_{h,3}^{\prime \prime }(jh,kh,x,y)=-h^{2}\sum_{\left| \mu \right| =2}\frac{%
m^{\mu }(jh,y)}{\mu !}(L-\widetilde{L})D^{\mu }\lambda (x)
\\ && \qquad
-3\sum_{\left| \mu \right| =3}\int_{0}^{1}(1-\delta )^{2}d\delta \int
q(jh,y,\theta )\frac{\widetilde{h}(\theta )^{\mu }}{\mu !}(L-\widetilde{L}%
)D^{\mu }\lambda (x+\delta \widetilde{h}(\theta ))d\theta .  \label{eq:012h}
\end{eqnarray}

Thus in Lemma 8 we can replace $K_{h}^{\prime
}(jh,kh,x,y)+M_{h}^{\prime
}(jh,kh,x,y)$ \ by \ $K_{h}(jh,kh,x,y)+M_{h}(jh,kh,x,y)$ where $%
K_{h}(jh,kh,x,y)=(L-\widetilde{L})\widetilde{p}_{h}(jh,kh,x,y),$ $%
M_{h}(jh,kh,x,y)=\frac{h}{2}(L_{\ast }^{2}-2L\widetilde{L}+\widetilde{L}%
^{2})\lambda (x)+M_{h}^{\prime \prime },M_{h}^{\prime \prime
}=M_{h,1}(jh,kh,x,y)+M_{h,2}(jh,kh,x,y)+M_{h,3}^{\prime \prime }(jh,kh,x,y)$
 and
\[
\max \{\left| M_{h}^{\prime }(jh,kh,x,y)\right| ,\left|
M_{h}(jh,kh,x,y)\right| \}\leq C\rho ^{-1}\zeta _{\rho }(y-x),
\]
$\rho ^{2}=kh-jh.$ To show this we note that
\[
\widetilde{p}_{h}(jh,kh,x,y)=\int q(jh,y,\theta )\lambda (x+\widetilde{h}%
(\theta ))d\theta
\]
where $\widetilde{h}(\theta )=m(jh,y)h+h^{1/2}\theta .$ From the Tailor
expansion we get
\begin{eqnarray*}
&&\widetilde{p}_{h}(jh,kh,x,y)=\lambda (x)+h\widetilde{L}\lambda
(x)+h^{2}\sum_{\left| \mu \right| =2}\frac{m^{\mu }(jh,y)}{\mu
!}D^{\mu }\lambda (x)\\ && \qquad +3\sum_{\left| \mu \right|
=3}\int_{0}^{1}(1-\delta )^{2}d\delta \int q(jh,y,\theta
)\frac{\widetilde{h}(\theta )^{\mu }}{\mu !}D^{\mu }\lambda
(x+\delta \widetilde{h}(\theta ))d\theta
\end{eqnarray*}
and, hence,
\[
K_{h}^{\prime }(jh,kh,x,y)=K_{h}(jh,kh,x,y)+(L-\widetilde{L})[\lambda (x)-%
\widetilde{p}_{h}(jh,kh,x,y)]
\]
\begin{equation}
=K_{h}(jh,kh,x,y)+h(\widetilde{L}^{2}-L\widetilde{L})\lambda
(x)+M_{h,3}^{\prime \prime }(jh,kh,x,y).  \label{eq:12i}
\end{equation}
Note that
\begin{eqnarray*}
h(\widetilde{L}^{2}-L\widetilde{L})\lambda (x)+M_{h,3}^{\prime
}(jh,kh,x,y)&=&h(\widetilde{L}^{2}-L\widetilde{L})\lambda (x)+\frac{h}{2}%
(L_{\ast }^{2}-\widetilde{L}^{2})\lambda (x)
\\&=&\frac{h}{2}(L_{\ast }^{2}-2L\widetilde{L}+\widetilde{L}^{2})\lambda (x)
\end{eqnarray*}
and from the definitions of the operators $L,\widetilde{L}$ and $%
L_{\ast }$ and Lipschitz conditions on the coefficients $m(t,x)$ and $%
\sigma (t,x)$ we obtain that
\begin{equation}
\left| \frac{h}{2}(L_{\ast }^{2}-2L\widetilde{L}+\widetilde{L}^{2})\lambda
(x)\right| \leq Ch\rho ^{-3}\zeta _{\rho }(y-x).  \label{eq:12j}
\end{equation}
Analogously, we have
\begin{eqnarray}
&&\left| h^{2}\sum_{\left| \mu \right| =2}\frac{m^{\mu }(jh,y)}{\mu !}(L-%
\widetilde{L})D^{\mu }\lambda (x)\right| \leq Ch^{2}\rho ^{-3}\zeta _{\rho
}(y-x),  \label{eq:12k}
\\ \label{eq:12l}
&&\left| 3\sum_{\left| \mu \right| =3}\int_{0}^{1}(1-\delta
)^{2}d\delta \int
q(jh,y,\theta )\frac{\widetilde{h}(\theta )^{\mu }}{\mu !}(L-\widetilde{L}%
)D^{\mu }\lambda (x+\delta \widetilde{h}(\theta ))d\theta \right|
\\ \nonumber && \qquad \leq Ch^{3/2}\rho ^{-4}\zeta _{\rho }(y-x).
\end{eqnarray}
Now (\ref{eq:12i})-(\ref{eq:12l}) imply the assertion of this remark.

\textbf{Lemma 9. }\textit{The following bounds holds:}
\begin{eqnarray} \nonumber
&&\left| \sum_{r=0}^{n}\widetilde{p}_{h}\otimes
_{h}(K_{h}+M_{h}+R_{h})^{(r)}(0,T,x,y)-\sum_{r=0}^{n}\widetilde{p}%
_{h}\otimes _{h}(K_{h}+M_{h})^{(r)}(0,T,x,y)\right|
\\   \label{eq:13} && \qquad
\leq C(\varepsilon )hn^{-1/2+\varepsilon }\zeta
_{\sqrt{T}}^{S}(y-x), \end{eqnarray} \textit{where}
$\lim_{\varepsilon \longrightarrow 0}C(\varepsilon )=+\infty . $

\textsc{Proof of Lemma 9. \ }\ For $r=1$ \ we will show that for any $%
\varepsilon >0$ with $\rho ^{2}=kh$
\begin{eqnarray}  \label{eq:14}
&&
\left| \widetilde{p}_{h}\otimes _{h}(K_{h}+M_{h}+R_{h})(0,kh,x,y)-\widetilde{%
p}_{h}\otimes _{h}(K_{h}+M_{h})(0,kh,x,y)\right|
\\ \nonumber
&& =\left| \widetilde{p}_{h}\otimes _{h}R_{h}(0,kh,x,y)\right|
\leq Ch^{3/2-\varepsilon }(kh)^{-1/2+\varepsilon
}B(\frac{1}{2},\varepsilon )\zeta _{\rho }^{S}(y-x). \label{eq:14}
\end{eqnarray}
Clearly, for an estimate of $\widetilde{p}_{h}\otimes
_{h}R_{h}(0,T,x,y)$
 it suffices to estimate
\[
I_{1}=h^{3/2}\sum_{j=0}^{k-2}h\int \widetilde{p}%
_{h}(0,kh,x,z)(f(jh,z)-f(jh,y))D_{z}^{\nu }\widetilde{p}%
_{h}((j+1)h,kh,z,y)dz
\]
for $\nu $ with $\left| \nu \right| =4,$ and
\begin{eqnarray*}
I_{2}&=&h^{2}\sum_{j=0}^{k-2}h\int \widetilde{p}%
_{h}(0,jh,x,z)(f(jh,z)-f(jh,y))\int q(jh,z,\theta )\widetilde{h}^{\nu
-e_{k}}(\theta )
\\ && \qquad
\times \int_{0}^{1}(1-u)^{4}D_{z}^{\nu }\lambda (z+u\widetilde{h}(\theta
))dud\theta dz
\end{eqnarray*}
for $\nu $ with $\left| \nu \right| =5,$ $1\leq k\leq d.$ Here
$f(t,x)$ is a function whose first and second derivatives with
respect to $x$ are continuous and bounded, uniformly in $t$ and
$x$. After integration by parts we obtain for $1\leq l,s\leq d$
\begin{eqnarray*}
I_{1}&=&-h^{3/2}\sum_{j=0}^{k-2}h\int D_{z}^{e_{k}}\widetilde{p}%
_{h}(0,jh,x,z)(f(jh,z)-f(jh,y))D_{z}^{\nu -e_{k}}\widetilde{p}%
_{h}((j+1)h,kh,z,y)dz
\\&&
+h^{3/2}\sum_{j=0}^{k-2}h\int D_{z}^{e_{s}}\widetilde{p}%
_{h}(0,jh,x,z)D_{z}^{e_{k}}f(jh,z)D_{z}^{\nu -e_{k}-e_{s}}\widetilde{p}%
_{h}((j+1)h,kh,z,y)dz
\\ &&
+h^{3/2}\sum_{j=0}^{k-2}h\int \widetilde{p}%
_{h}(0,jh,x,z)D_{z}^{e_{k}+e_{s}}f(jh,z)D_{z}^{\nu -e_{k}-e_{s}}\widetilde{p}%
_{h}((j+1)h,kh,z,y)dz.
\end {eqnarray*}
Hence,
\begin{equation}
\left| I_{1}\right| \leq Ch^{3/2}\sum_{j=0}^{k-2}h\frac{1}{\sqrt{jh}(kh-jh)}%
\zeta _{\rho }(y-x)\leq Ch^{3/2-\varepsilon }(kh)^{-1/2+\varepsilon }B(\frac{%
1}{2},\varepsilon )\zeta _{\sqrt{T}}^{S}(y-x).  \label{eq:15}
\end{equation}
In the same way, we obtain after integration by parts with $1\leq
l,s\leq d$
\begin{eqnarray*}
I_{2}&=&-h^{2}\sum_{j=0}^{k-2}h\int_{0}^{1}(1-u)^{4}du\int d\theta
\left(m(jh,y)h^{1/2}+\theta\right )^{\nu -e_{l}}\\ && \qquad
\int D_{z}^{e_{l}}\widetilde{p}%
_{h}(0,jh,x,z)(f(jh,z)-f(jh,y))
\\ && \qquad
\times q(jh,z,\theta )D_{z}^{\nu -e_{l}}\widetilde{p}_{h}((j+1)h,kh,z+u%
\widetilde{h}(\theta ),y)dz \end{eqnarray*}
\begin{eqnarray*} && +h^{2}\sum_{j=0}^{k-2}h\int_{0}^{1}(1-u)^{4}du\int d\theta
(m(jh,y)h^{1/2}+\theta )^{\nu -e_{l}}
\\ && \qquad
\times \int D_{z}^{e_{s}}[\widetilde{p}_{h}(0,jh,x,z)D^{e_{l}}f(jh,z)q(jh,z,%
\theta )]\\ && \qquad
\times D_{z}^{\nu -e_{l}-e_{s}}\widetilde{p}_{h}((j+1)h,kh,z+u\widetilde{h}%
(\theta ),y)dz
\\ &&
-h^{2}\sum_{j=0}^{k-2}h\int_{0}^{1}(1-u)^{4}du\int d\theta
(m(jh,y)h^{1/2}+\theta )^{\nu -e_{l}}\\ && \qquad
\times \widetilde{p}%
_{h}(0,jh,x,z)(f(jh,z)-f(jh,y))
\\ && \qquad
\times D_{z}^{e_{l}}q(jh,z,\theta )D_{z}^{\nu -e_{l}}\widetilde{p}%
_{h}((j+1)h,kh,z+u\widetilde{h}(\theta ),y)dz.
\end{eqnarray*}
It follows from this equation that
\begin{equation}
\left| I_{2}\right| \leq Ch^{3/2-\varepsilon }(kh)^{-1/2+\varepsilon }B(%
\frac{1}{2},\varepsilon )\zeta _{\rho }^{S}(y-x).  \label{eq:17}
\end{equation}
Claim (\ref{eq:14}) follows now from (\ref{eq:15}) and (\ref{eq:17}). For $%
r\geq 2$ \ we use the identity
\begin{eqnarray} \label{eq:18}
&&\widetilde{p}_{h}\otimes _{h}(K_{h}+M_{h}+R_{h})^{(r)}(0,T,x,y)-\widetilde{p}%
_{h}\otimes _{h}(K_{h}+M_{h})^{(r)}(0,T,x,y)
\\ \nonumber
&&=\left[ \widetilde{p}_{h}\otimes _{h}(K_{h}+M_{h}+R_{h})^{(r-1)}-\widetilde{p%
}_{h}\otimes _{h}(K_{h}+M_{h})^{(r-1)}\right] \otimes
_{h}(K_{h}+M_{h})(0,T,x,y)
\\ \nonumber
&& \qquad +\widetilde{p}_{h}\otimes
_{h}(K_{h}+M_{h}+R_{h})^{(r-1)}\otimes _{h}R_{h}(0,T,x,y)\\
\nonumber &&=I+II.
\end{eqnarray}
For $r=2$ we obtain from (\ref{eq:14}) and from the simple
estimate $\left| (K_{h}+M_{h})(jh,kh,z,y)\right| \leq C\rho
_{2}^{-1}\zeta _{\rho _{2}}^{S}(y-z)$ with $\rho _{2}^{2}=kh-jh,$
that
\begin{eqnarray*}
\left| I\right| &=&\left| \left[ \widetilde{p}_{h}\otimes
_{h}(K_{h}+M_{h}+R_{h})-\widetilde{p}_{h}\otimes
_{h}(K_{h}+M_{h})\right] \otimes
_{h}(K_{h}+M_{h})(0,kh,x,y)\right|
\\ &
\leq & C^{2}h^{3/2-\varepsilon }B(\frac{1}{2},\varepsilon
)\sum_{j=0}^{k-2}h(jh)^{-1/2+\varepsilon }(kh-jh)^{-1/2}\int \zeta
_{\rho _{1}}^{S}(z-x)\zeta _{\rho _{2}}^{S}(y-z)dz
\\ &
\leq & C^{2}h^{3/2-\varepsilon }B(\frac{1}{2},\varepsilon )B(\frac{1}{2}%
,\varepsilon +\frac{1}{2})(kh)^{\varepsilon }\zeta _{\rho
}^{S}(y-x),\rho ^{2}=kh. \end{eqnarray*}

 For $r\geq 3$ we obtain by
induction
\begin{eqnarray} \nonumber
\left| I\right| &=&\left| \left[ \widetilde{p}_{h}\otimes
_{h}(K_{h}+M_{h}+R_{h})^{(r-1)}-\widetilde{p}_{h}\otimes
_{h}(K_{h}+M_{h})^{(r-1)}\right] \right .\\ \nonumber &&
\left.\qquad \otimes _{h}(K_{h}+M_{h})(0,kh,x,y)\right| \\
\nonumber &
\leq& C^{r}h^{3/2-\varepsilon }B(\frac{1}{2},\varepsilon )B(\frac{1}{2}%
,\varepsilon +\frac{1}{2})...B(\frac{1}{2},\varepsilon +\frac{r-1}{2}%
)(kh)^{\varepsilon +(r-2)/2}\zeta _{\rho }^{S}(y-x)
\\ &
\leq &\Gamma (\varepsilon )h^{3/2-\varepsilon }\frac{[C\Gamma (1/2)]^{r}}{%
\Gamma (\varepsilon +\frac{r}{2})}(kh)^{\varepsilon +(r-2)/2}\zeta _{\rho
}^{S}(y-x),\rho ^{2}=kh.  \label{eq:19}
\end{eqnarray}
To estimate $II$ \ we use the following estimates
\begin{eqnarray}
&&\left| D_{y}^{a}D_{x}^{b}\widetilde{p}_{h}(jh,kh,x,y)\right|
\leq C\rho ^{-\left| a\right| -\left| b\right| }\varsigma _{\rho
}^{S-\left| a\right|
}(y-x),D_{x}^{b}\widetilde{p}_{h}(jh,kh,x,x+v)\leq C\zeta _{\rho
}^{S}(v), \label{eq:20}
\\ &&
\left| D_{x}^{b}(K_{h}+M_{h}+R_{h})(jh,kh,x+v,x)\right| \leq C\rho
^{-1}\zeta _{\rho }^{S}(v).  \label{eq:21}
\end{eqnarray}
The inequalities (\ref{eq:20}) and (\ref{eq:21}) are obained by
using the same arguments as is the proof of Lemma 7. Using these
inequalities and mimicking the proof of Theorem 2.3 in Konakov and
Mammen (2002) we obtain the following bounds for $r \geq 0$
\begin{eqnarray} \nonumber
&&\left| D_{x}^{b}D_{y}^{a}\widetilde{p}_{h}\otimes
_{h}(K_{h}+M_{h}+R_{h})^{(r)}(0,kh,x,y)\right|
\\ && \qquad \nonumber
\leq C^{r}(kh)^{-\left| a\right| -\left| b\right| +r}B(\frac{1}{2},\frac{1}{2%
})B(1,\frac{1}{2})...B(\frac{r}{2},\frac{1}{2})\varsigma _{\rho }^{S-\left|
a\right| }(y-x)
\\ && \qquad
\leq \frac{\lbrack C\Gamma (1/2)]^{r}}{\Gamma
(\frac{r+1}{2})}(kh)^{-\left| a\right| -\left| b\right|
+r}\varsigma _{\rho }^{S-\left| a\right| }(y-x).  \label{eq:22}
\end{eqnarray}
Inequality (\ref{eq:22}) allows us to estimate $II=[\widetilde{p}%
_{h}\otimes _{h}(K_{h}+M_{h}^{\prime \prime
}+R_{h})^{(r-1)}]\otimes _{h}R_{h}(0,kh,x,y)$. For this aim it
suffices to estimate
\begin{eqnarray} \nonumber
&&h^{3/2}\sum h\int [\widetilde{p}_{h}\otimes
_{h}(K_{h}+M_{h}+R_{h})^{(r-1)}](0,jh,x,z)
\\ && \qquad
\times D^{\nu }\widetilde{p}_{h}((j+1)h,kh,z,y)(f(jh,z)-f(jh,y))dz
\label{eq:23}
\end{eqnarray}
for $r\geq 2$ ,\ $\left| \nu \right| =4,$ \ and \
\begin{eqnarray} \nonumber
&&\sum_{j=0}^{n-2}h\int [\widetilde{p}_{h}\otimes
_{h}(K_{h}+M_{h}+R_{h})^{(r-1)}](0,jh,x,z)(f(jh,z)-f(jh,y))
\\ && \qquad
\times \int q(jh,z,\theta )\widetilde{h}^{\nu -e_{l}}(\theta
)\int_{0}^{1}(1-u)^{4}D^{\nu }\widetilde{p}_{h}((j+1)h,kh,z+u\widetilde{h}%
(\theta ),y)dud\theta dz  \label{eq:23a}
\end{eqnarray}
for $r\geq 2,$ $\left| \nu \right| =5,$ $1\leq l\leq d$. Here $
f(t,x)$ is a function whose first and second derivatives with
respect to $x$ are continuous and bounded, uniformly in $t$ and
$x$. The upper bound for (\ref{eq:23}) follows from (\ref{eq:22})
by integration by parts, exactly in the same way as it was done to
obtain the upper bound for $I_{1}$ (see (\ref{eq:15})). This gives
an estimate
\begin{eqnarray} \nonumber
&& h^{3/2}\left| \sum h\int [\widetilde{p}_{h}\otimes
_{h}(K_{h}+M_{h}+R_{h})^{(r-1)}](0,jh,x,z)\right.
\\ && \qquad  \nonumber
\left. \times D^{\nu }\widetilde{p}_{h}((j+1)h,kh,z,y)(f(jh,z)-f(jh,y))dz%
\right|
\\ &&
\leq \Gamma (\varepsilon )h^{3/2-\varepsilon }\frac{[C\Gamma (1/2)]^{r}}{%
\Gamma (\frac{r+1}{2})}(kh)^{\varepsilon +(r-2)/2}\zeta _{\rho
}^{S}(y-x). \label{eq:24}
\end{eqnarray}
The upper bound for (\ref{eq:23a}) also follows from (\ref{eq:22})
 by integration by parts in the same way as it was done to obtain
an upper bound for $I_{2}$ (see (\ref{eq:17})). This gives for
(\ref{eq:23a}) the same estimate as in (\ref{eq:24}) and, hence,
\begin{equation}
\left| II\right| \leq C\Gamma (\varepsilon )h^{3/2-\varepsilon }\frac{%
[C\Gamma (1/2)]^{r}}{\Gamma (\frac{r+1}{2})}(kh)^{\varepsilon
+(r-2)/2}\zeta _{\rho }^{S}(y-x) . \label{eq:25}
\end{equation}
The assertion of the lemma follows now from (\ref{eq:14}),
(\ref{eq:18}), (\ref{eq:19}) \ and (\ref{eq:25}).

\textbf{Lemma \ 10. }\textit{Let
}$A(s,t,x,y),B(s,t,x,y),C(s,t,x,y)$ some
functions with the absolute value less than $C(t-s)^{-1/2}\zeta _{\sqrt{t-s%
}}(y-x)$ \textit{for a constant }$C.$ \textit{Then}
\begin{eqnarray*}
&& \sum_{r=0}^{\infty }A\otimes
_{h}(B+C)^{(r)}(ih,jh,x,y)-\sum_{r=0}^{\infty }A\otimes
_{h}B^{(r)}(ih,jh,x,y)
\\ &&
=\sum_{r=1}^{\infty }\left[ A\otimes _{h}\Phi \right] \otimes _{h}\left[
C\otimes _{h}\Phi \right] ^{(r)}(ih,jh,x,y),
\end{eqnarray*}
\textit{where }$\ \Phi =\sum_{r=0}^{\infty }B^{(r)}.$

\textsc{Proof of Lemma 10. }Under the conditions of the lemma all
series are absolutely convergent. The assertion of the lemma
follows from linearity of the operation $\otimes _{h}$ and after
permutation of the summands in the absolutely convergent series.

\section{\label{sec:mainres} Proof of the main result.}

We now come to the proof of \ Theorem 1. The main tools for the
proof have been given in Sections 4 and 5. From Lemmas 1 and 2 we
get that
\[
p(0,T,x,y)=\sum_{r=0}^{n}\tilde{p}\otimes H^{(r)}(0,T,x,y)+o(h^{2}T)\phi _{C,%
\sqrt{T}}(y-x).
\]
With Lemma 3 this gives
\begin{equation}
p(0,T,x,y)-p_{h}(0,T,x,y)=T_{1}+...+T_{7}+o(h^{2}T)\phi _{C,\sqrt{T}}(y-x),
\label{eq:26}
\end{equation}
where
\begin{eqnarray*}
&&T_{1}=\sum_{r=0}^{n}\tilde{p}\otimes H^{(r)}(0,T,x,y)-\sum_{r=0}^{n}\tilde{p}%
\otimes _{h}H^{(r)}(0,T,x,y),
\\&&
T_{2}=\sum_{r=0}^{n}\tilde{p}\otimes _{h}H^{(r)}(0,T,x,y)-\sum_{r=0}^{n}%
\tilde{p}\otimes _{h}(H+M_{h}^{\prime \prime }+\sqrt{h}%
N_{1})^{(r)}(0,T,x,y),
\\&&
T_{3}=\sum_{r=0}^{n}\tilde{p}\otimes _{h}(H+M_{h}^{\prime \prime }+\sqrt{h}%
N_{1})^{(r)}(0,T,x,y)\\&& \qquad -\sum_{r=0}^{n}\tilde{p}\otimes _{h}(H+M_{h}+\sqrt{h}%
N_{1})^{(r)}(0,T,x,y),
\\&&
T_{4}=\sum_{r=0}^{n}\tilde{p}\otimes _{h}(H+M_{h}+\sqrt{h}%
N_{1})^{(r)}(0,T,x,y)-\sum_{r=0}^{n}\tilde{p}\otimes
_{h}(K_{h}+M_{h})^{(r)}(0,T,x,y),
\\&&
T_{5}=\sum_{r=0}^{n}\tilde{p}\otimes
_{h}(K_{h}+M_{h})^{(r)}(0,T,x,y)-\sum_{r=0}^{n}\tilde{p}_{h}\otimes
_{h}(K_{h}+M_{h})^{(r)}(0,T,x,y),
\\&&
T_{6}=\sum_{r=0}^{n}\tilde{p}_{h}\otimes
_{h}(K_{h}+M_{h})^{(r)}(0,T,x,y)-\sum_{r=0}^{n}\tilde{p}_{h}\otimes
_{h}(K_{h}+M_{h}+R_{h})^{(r)}(0,T,x,y),
\\&&
T_{7}=\sum_{r=0}^{n}\tilde{p}_{h}\otimes
_{h}(K_{h}+M_{h}+R_{h})^{(r)}(0,T,x,y)-\sum_{r=0}^{n}\tilde{p}_{h}\otimes
_{h}H_{h}^{(r)}(0,T,x,y).
\end{eqnarray*}
Here we put $N_{1}(s,t,x,y)=(L-\widetilde{L})\widetilde{\pi }_{1}(s,t,x,y).$

We now discuss the asymptotic behavior of the terms
$T_{1},...,T_{7}$.

\noindent \textit{Asymptotic treatment of the term $T_{1}$.} Using
Remark 1 in Konakov (2006) we get that
\begin{eqnarray} \nonumber
T_{1}&=&\frac{h}{2}[p\otimes _{h}(L^{2}-2L\widetilde{L}+\widetilde{L}^{2})%
\widetilde{p}\otimes _{h}\Phi )](0,T,x,y)
\\ &&
+\frac{h}{2}[p\otimes _{h}(L^{\prime }-\widetilde{L}^{\prime })\widetilde{p}%
\otimes _{h}\Phi ](0,T,x,y) +R_{T}(0,T,x,y),  \label{eq:28}
\end{eqnarray}
where for any $0<\varepsilon <1/2$%
\begin{equation}
\left| R_{T}(0,T,x,y)\right| \leq C(\varepsilon )(hn^{-1/2+\varepsilon }+h%
\sqrt{T})\phi _{C,\sqrt{T}}(y-x),  \label{eq:29}
\end{equation}
$\Phi (s,t,x,y)=\sum_{r=0}^{\infty }H^{(r)}(s,t,x,y)$. Here, the
summand $H^{(0)}(s,t,x,y)$ was introduced to simplify notations.
We define $g\otimes _{h}H^{(0)}(s,t,x,y)=g(s,t,x,y)$ for a function $%
g$. Note, that in the homogenous case it holds that $\sigma
_{ij}(s,x)=\sigma _{ij}(x),m_{i}(s,x)=m_{i}(x)$ and that the
second summand in (\ref{eq:28}) is equal to $0.$

\textit{Asymptotic treatment of the term $T_{2}$. }We will show
that with a positive constant $\delta >0$
\begin{eqnarray} \label{eq:T201}
&&\left| T_{2}-3\sum_{r=0}^{\infty }\tilde{p}\otimes
_{h}H^{(r)}(0,T,x,y)+\sum_{r=0}^{\infty }\tilde{p}\otimes _{h}(H+M_{h,1}+%
\sqrt{h}N_{1})^{(r)}(0,T,x,y)\right.
\\ \nonumber &&
\left. +\sum_{r=0}^{\infty }\tilde{p}\otimes
_{h}(H+M_{h,2})^{(r)}(0,T,x,y)+\sum_{r=0}^{\infty
}\tilde{p}\otimes _{h}(H+M_{h,3}^{\prime \prime
})^{(r)}(0,T,x,y)\right|\\  && \leq Chn^{-\delta }\zeta
_{\sqrt{T}}(y-x).\nonumber
\end{eqnarray}
It suffices to discuss the
case $r\geq 2$ because for $r=1,2$ the left hand side of (\ref{eq:T201}%
) is equal to zero. For $r\geq 2$, (\ref{eq:T201}) immediately
follows from the following bounds
\begin{eqnarray} \label{eq:T202}
&& \left| \tilde{p}\otimes _{h}(H+M_{h}^{\prime \prime }+\sqrt{h}%
N_{1})^{(r)}(0,T,x,y)\right.
\\ && \nonumber
\qquad -\tilde{p}\otimes
_{h}(H+M_{h,1}+M_{h,2}+\sqrt{h}N_{1})^{(r)}(0,T,x,y)
\\ && \nonumber \qquad
\left. -[\tilde{p}\otimes _{h}(H+M_{h,3}^{\prime \prime })^{(r)}-\tilde{p}%
\otimes _{h}H^{(r)}](0,T,x,y)\right|
\\ && \nonumber
\leq C(\varepsilon )h^{3/2-2\varepsilon }\frac{C^{r}}{\Gamma (\frac{r-1}{2})}%
T^{3\varepsilon +\frac{r-3}{2}}\zeta _{\sqrt{kh}}(v-x),
\\ \label{eq:T203} &&
\left| \tilde{p}\otimes _{h}(H+M_{h,1}+M_{h,2}+\sqrt{h}N_{1})^{(r)}(0,T,x,y)%
\right.
\\ && \nonumber \qquad
-\tilde{p}\otimes _{h}(H+M_{h,1}+\sqrt{h}N_{1})^{(r)}(0,T,x,y)
\\ && \nonumber \qquad
\left. -[\tilde{p}\otimes _{h}(H+M_{h,2})^{(r)}-\tilde{p}\otimes
_{h}H^{(r)}](0,T,x,y)\right|
\\ && \nonumber \leq
C(\varepsilon )h^{3/2-2\varepsilon }\frac{C^{r}}{\Gamma (\frac{r-1}{2})}%
T^{3\varepsilon +\frac{r-4}{2}}\zeta _{\sqrt{T}}(v-x),
\end{eqnarray}
for all sufficiently small $\varepsilon >0$ with a function $C$
that fulfils $\lim_{\varepsilon \longrightarrow 0}C(\varepsilon
)=+\infty$. We will prove first the bound (\ref{eq:T202}). Denote
the expression under the sign of the absolute value in
(\ref{eq:T202}) by $\Gamma _{r}$. Note that $\Gamma _{0}=\Gamma
_{1}=0$. For $r\geq 2$ we make use of the following recurrence
formula
\begin{eqnarray} \label{eq:T203aa}
&&\Gamma _{r}=\Gamma _{r-1}\otimes _{h}H+\left[ \tilde{p}\otimes
_{h}(H+M_{h}^{\prime \prime }+\sqrt{h}N_{1})^{(r-1)}\right.
\\ \nonumber && \qquad
\left. -\tilde{p}\otimes _{h}(H+M_{h,1}+M_{h,2}+\sqrt{h}N_{1})^{(r-1)}\right]
\otimes _{h}(M_{h}^{\prime \prime }+\sqrt{h}N_{1})
\\ \nonumber && \qquad
+\left[ \tilde{p}\otimes _{h}(H+M_{h,1}+M_{h,2}+\sqrt{h}N_{1})^{(r-1)}-%
\tilde{p}\otimes _{h}(H+M_{h,3}^{\prime \prime })^{(r-1)}\right]
\otimes _{h}M_{h,3}^{\prime \prime }
\\ \nonumber && =I+II+III.
\end{eqnarray}
We start with the estimation of $II$. First we will estimate
\begin{equation}
\left| \tilde{p}\otimes _{h}(H+M_{h}^{\prime \prime }+\sqrt{h}N_{1})^{(r-1)}-%
\tilde{p}\otimes _{h}(H+M_{h,1}+M_{h,2}+\sqrt{h}N_{1})^{(r-1)}\right| .
\label{eq:T203a}
\end{equation}
For $r=2$ we have $\tilde{p}\otimes _{h}M_{h,3}^{\prime \prime
}(0,kh,x,y)$. It follows from (\ref{eq:012h}) that it is enough to
estimate
\begin{equation}
J_{1}=h^{2}\sum_{i=0}^{k-2}h\int \widetilde{p}%
(0,ih,x,v)(f(ih,v)-f(ih,y))D_{v}^{\nu }\widetilde{p}_{h}((i+1)h,kh,v,y)dv
\label{eq:T204}
\end{equation}
for $\left| \nu \right| =4$ and
\begin{eqnarray} \label{eq:T205} &&
J_{2}=h^{3/2}\sum_{i=0}^{k-2}h\int \widetilde{p}(0,ih,x,v)(f(ih,v)-f(ih,y))%
\int q(ih,v,\theta )\theta ^{\nu }\int_{0}^{1}(1-\delta )^{2}
\\ && \nonumber \qquad
\times D_{v}^{\nu
+e_{l}+e_{q}}\widetilde{p}_{h}((i+1)h,kh,v+\delta
\widetilde{h}(\theta ),y)d\delta d\theta dv
\end{eqnarray}
for $\left| \nu \right| =3.$ Here $f(t,x)$ is a function for which
$D_{x}^{\nu }f(t,x)$ is bounded uniformly in $(t,x)$ for $\left|
\nu \right| =0,1,2,3$. An estimate for$J_{1}$ follows from
(\ref{eq:15}). This gives
\begin{equation}
\left| J_{1}\right| \leq Ch^{2-\varepsilon }(kh)^{\varepsilon -1/2}B(\frac{1%
}{2},\varepsilon )\zeta _{\sqrt{kh}}^{S}(y-x).  \label{eq:T206}
\end{equation}
The estimate for $J_{2}$ \ can be obtained analogously to the estimate of $\
I_{2}$ \ (see (\ref{eq:17})). Integrating by parts we get
\begin{eqnarray*}
J_{2}&=&h^{3/2}\sum_{j=0}^{k-2}h\int_{0}^{1}(1-\delta )^{2}d\delta
\int d\theta \cdot \theta ^{\nu }\int
D_{v}^{e_{l}+e_{q}}[\widetilde{p}(0,ih,x,v)
\\ &&
\times (f(ih,v)-f(ih,y))q(ih,v,\theta )]D_{v}^{\nu }\widetilde{p}%
_{h}((i+1)h,kh,v+\delta \widetilde{h}(\theta ),y)dv.
\end{eqnarray*}
The derivative
\[
D_{v}^{e_{l}+e_{q}}[\widetilde{p}(0,ih,x,v)(f(ih,v)-f(ih,y))q(ih,v,\theta )]
\]
is a sum of 9 summands. Integrating by parts once more for
summands which contain $D_{\nu }^{\mu }\widetilde{p}(0,ih,x,v)$
with $\left| \mu \right| <2$, we obtain
\begin{eqnarray} \label{eq:T207}
\left| J_{2}\right| &\leq& Ch^{3/2-2\varepsilon }\zeta _{\sqrt{kh}%
}^{S-3}(y-x)\int \psi (\theta )\left\| \theta \right\|
^{3}(h^{(S-3)/2}\left\| \theta \right\| ^{S-3}+1)d\theta
\\ && \qquad \nonumber
\sum_{j=0}^{k-2}h\frac{1}{(ih)^{1-\varepsilon }}\times \frac{1}{%
(kh-ih)^{1-\varepsilon }}\\  \nonumber&\leq &Ch^{3/2-2\varepsilon
}B(\varepsilon ,\varepsilon )(kh)^{2\varepsilon -1}\zeta
_{\sqrt{kh}}^{S-3}(y-x)
\end{eqnarray}
for any $\varepsilon \in (0,1/4).$ It follows from (\ref{eq:T206})
and (\ref
{eq:T207}) that for $r=2$ \ (\ref{eq:T203a}) does not exceed $%
Ch^{3/2-2\varepsilon }B(\varepsilon ,\varepsilon
)(kh)^{\varepsilon -1}\zeta _{\sqrt{kh}}^{S-3}(y-x).$ For $r\geq
3$ we use the recurrence relation \begin{eqnarray} \label{eq:T208}
&& \tilde{p}\otimes _{h}(H+M_{h}^{\prime \prime }+\sqrt{h}N_{1})^{(r-1)}-\tilde{%
p}\otimes _{h}(H+M_{h,1}+M_{h,2}+\sqrt{h}N_{1})^{(r-1)}
\\ \nonumber &&
=\left[ \tilde{p}\otimes _{h}(H+M_{h}^{\prime \prime }+\sqrt{h}%
N_{1})^{(r-2)}-\tilde{p}\otimes _{h}(H+M_{h,1}+M_{h,2}+\sqrt{h}N_{1})^{(r-2)}%
\right]
\\ \nonumber && \qquad
\otimes _{h}(H+M_{h}^{\prime \prime
}+\sqrt{h}N_{1})+[\tilde{p}\otimes
_{h}(H+M_{h,1}+M_{h,2}+\sqrt{h}N_{1})^{(r-2)}]\otimes
_{h}M_{h,3}^{\prime \prime }
\\ \nonumber &&
=I^{\prime }+II^{\prime }.
\end{eqnarray}
From (\ref{eq:T208}) we obtain for $r=3$
\begin{eqnarray*}
\left| I^{\prime }\right| &\leq& Ch^{3/2-2\varepsilon
}B(\varepsilon
,\varepsilon )\zeta _{\sqrt{kh}}^{S-3}(y-x)\sum_{i=0}^{k-2}h(ih)^{%
\varepsilon -1}(kh-ih)^{-1/2}
\\ &\leq& Ch^{3/2-2\varepsilon }B(\varepsilon ,\varepsilon )B(\frac{1}{2}%
,\varepsilon )(kh)^{\varepsilon -1/2}\zeta
_{\sqrt{kh}}^{S-3}(y-x).
\end{eqnarray*}
For the estimate of $II^{\prime }$ we use the following estimates
\begin{eqnarray}
&& \left|
D_{v}^{a}D_{x}^{b}(H+M_{h,1}+M_{h,2}+\sqrt{h}N_{1})(jh,kh,x,v)\right|
\leq C\rho ^{-1-\left| a\right| -\left| b\right| }\zeta _{\rho
}(v-x), \label{eq:T209}
\\ &&
\left|
D_{x}^{b}(H+M_{h,1}+M_{h,2}+\sqrt{h}N_{1})(jh,kh,x,x+v)\right|
\leq C\rho ^{-1}\zeta _{\rho }(v-x).  \label{eq:T210}
\end{eqnarray}
To prove (\ref{eq:T209}) it is enough to get the corresponding
estimates for summands in $M_{h,1},M_{h,2}$ and $ \sqrt{h}N_{1}$
\begin{eqnarray}  \label{eq:T211}
&& h^{1/2}D_{v}^{a}D_{x}^{b}\left[ D_{x}^{\nu }\widetilde{p}%
_{h}((j+1)h,kh,x,v)(f(jh,x)-f(jh,v))\right]\\ \nonumber && \leq
C\rho ^{-1-\left| a\right| -\left| b\right| }\zeta _{\rho
}(v-x)\mbox{ for }\left| \nu \right| =3, \\ &&  \label{eq:T212}
hD_{v}^{a}D_{x}^{b}\left[ D_{x}^{\nu }\widetilde{p}%
_{h}((j+1)h,kh,x,v)(f(jh,x)-f(jh,v))\right] \\ \nonumber && \leq
C\rho ^{-1-\left| a\right| -\left| b\right| }\zeta _{\rho
}(v-x)\mbox{ for }\left| \nu \right| =4,
\\ && \label{eq:T213}
h^{1/2}D_{v}^{a}D_{x}^{b}\left[ D_{x}^{\nu +e_{p}+e_{q}}\widetilde{p}%
_{h}((j+1)h,kh,x,v)\rho ^{2}(f(jh,x)-f(jh,v))\right] \\ \nonumber
&&\leq C\rho ^{-1-\left| a\right| -\left| b\right| }\zeta _{\rho
}(v-x)\mbox{ for }\left| \nu \right| =3,
\end{eqnarray}
for a function $f(t,x)$ with $\left| a\right| +\left| b\right| $ \
derivatives w.r.t.\ $x$ that are uniformly bounded w.r.t.\ $t$.
These estimates are direct consequences of Lemma 7. To prove
(\ref{eq:T210}) it is enough to show the corresponding estimates
for summands in $M_{h,1},M_{h,2}$
and $\ \sqrt{h}N_{1}$%
\begin{eqnarray} \label{eq:T214}
&&h^{1/2}D_{x}^{b}\left[ D_{x}^{\nu
}\widetilde{p}_{h}((j+1)h,kh,x,y)\mid
_{y=x+v}(f(jh,x)-f(jh,x+v))\right] \\ \nonumber &&\leq C\rho
^{-1}\zeta _{\rho }(v-x)\mbox { for }\left| \nu \right| =3,
\\  \label{eq:T215} &&
hD_{x}^{b}\left[ D_{x}^{\nu }\widetilde{p}_{h}((j+1)h,kh,x,y)\mid
_{y=x+v}(f(jh,x)-f(jh,x+v))\right]
\\ \nonumber &&
\leq C\rho ^{-1}\zeta _{\rho }(v-x) \mbox { for } \left| \nu
\right| =4,
\\ \label{eq:T216} &&
h^{1/2}D_{x}^{b}\left[ D_{x}^{\nu +e_{p}+e_{q}}\widetilde{p}%
_{h}((j+1)h,kh,x,v)\mid _{y=x+v}\rho ^{2}(f(jh,x)-f(jh,x+v))\right]
\\ \nonumber &&
\leq C\rho ^{-1}\zeta _{\rho }(v-x) \mbox { for } \left| \nu
\right| =3.
\end{eqnarray}
These estimates also follow from the estimates obtained in the proof of
Lemma 7. If \ $z(V_{j,k}^{-1/2}(y),\mu
_{j,k}(y),x,y)=V_{j,k}^{-1/2}(y)(y-x-\mu _{j,k}(y))$ \ then
\begin{equation}
\left| \frac{\partial z}{\partial y}\right| \leq \frac{C}{\rho },\left|
\frac{\partial z}{\partial x}\right| \leq \frac{C}{\rho }  \label{eq:T216a}
\end{equation}
and for $z(V_{j,k}^{-1/2}(x+v),\mu
_{j,k}(x+v),x,x+v)=V_{j,k}^{-1/2}(x+v)(v-\mu _{j,k}(x+v))$ \ we have
\begin{equation}
\left| \frac{\partial z}{\partial x}\right| \leq \left| \frac{\partial z}{%
\partial V_{j,k}^{-1/2}}\right| \left| \frac{\partial V_{j,k}^{-1/2}}{%
\partial x}\right| +\left| \frac{\partial z}{\partial \mu _{j,k}}\right|
\left| \frac{\partial \mu _{j,k}}{\partial x}\right| \leq C(\left\|
v\right\| +1).  \label{eq:T216b}
\end{equation}
Now with the inequalities (\ref{eq:T209}),(\ref{eq:T210}) we can
proceed like in the proof of Theorem 2.3 in Konakov and Mammen
(2002). This gives the following estimate for $r\geq 3$
\begin{eqnarray} \label{eq:T217}
&&\left| D_{v}^{a}D_{x}^{b}[\widetilde{p}\otimes _{h}(H+M_{h,1}+M_{h,2}+\sqrt{h%
}N_{1})^{(r-2)}](jh,kh,x,v)\right|
\\ \nonumber &&
\leq C^{r}B(1,\frac{1}{2})\times ...\times
B(\frac{r-1}{2},\frac{1}{2})\rho ^{r-2-\left| a\right| -\left|
b\right| }\zeta _{\rho }(v-x).
\end{eqnarray}
Now we denote $\widetilde{p}_{1,r}=\widetilde{p}\otimes
_{h}(H+M_{h,1}+M_{h,2}+\sqrt{h}N_{1})^{(r)}$ and
$\widetilde{p}_{0}=\widetilde{p}$. For an estimate of
$\widetilde{p}_{1,r-2}\otimes _{h}M_{h,3}^{\prime \prime }$ it
suffices to make the same calculations using integration by parts
as it was done above for $J_{1}$ and $J_{2}$. This gives
\begin{eqnarray} \label{eq:T217a}
\left| II^{\prime }\right| &\leq& \left|
\widetilde{p}_{1,r-2}\otimes _{h}M_{h,3}^{\prime \prime
}(0,kh,x,y)\right|
\\ \nonumber &
\leq& C^{r}h^{3/2-2\varepsilon }B(1,\frac{1}{2})B(1+\frac{1}{2},\frac{1}{2}%
)\\ \nonumber && \qquad \times ...\times
B(1+\frac{r-3}{2},\frac{1}{2})B(\frac{r-2}{2},\varepsilon
)(kh)^{\varepsilon +\frac{r-4}{2}}\zeta _{\sqrt{kh}}(v-x)
\end{eqnarray}
and by induction for $r\geq 3$
\begin{eqnarray} \label{eq:T217b}
&&
\left| I^{\prime }\right| \leq C^{r}h^{3/2-2\varepsilon }B(\varepsilon ,%
\frac{1}{2})B(\varepsilon +\frac{1}{2},\frac{1}{2})\\ \nonumber &&
\qquad \times ...\times B(\varepsilon
+\frac{r-3}{2},\frac{1}{2})B(\varepsilon ,\varepsilon
)(kh)^{\varepsilon +\frac{r-4}{2}}\zeta _{\sqrt{kh}}(v-x).
\end{eqnarray} Comparing (\ref{eq:T217a}) and (\ref{eq:T217b}) we obtain
that for $r\geq 3$%
\begin{eqnarray} \label{eq:T217c}
&&
\left| \tilde{p}\otimes _{h}(H+M_{h}^{\prime \prime }+\sqrt{h}N_{1})^{(r-1)}-%
\tilde{p}\otimes _{h}(H+M_{h,1}+M_{h,2}+\sqrt{h}N_{1})^{(r-1)}\right|
\\ \nonumber &&
\leq C^{r}h^{3/2-2\varepsilon }B(\varepsilon ,\frac{1}{2})B(\varepsilon +\frac{1}{%
2},\frac{1}{2})\\ \nonumber && \qquad \times ...\times B(\varepsilon +\frac{r-3}{2},\frac{1}{2}%
)B(\varepsilon ,\varepsilon )(kh)^{\varepsilon +\frac{r-4}{2}}\zeta _{\sqrt{%
kh}}(v-x).  \end{eqnarray}
From (\ref{eq:T217c}) we get the following estimate for $II$%
\begin{equation}
\left| II\right| \leq C^{r}h^{3/2-2\varepsilon }B(\varepsilon ,\varepsilon
)B(\varepsilon ,\frac{1}{2})B(\varepsilon +\frac{1}{2},\frac{1}{2})\times
...\times B(\varepsilon +\frac{r-2}{2},\frac{1}{2})T^{\varepsilon +\frac{r-3%
}{2}}\zeta _{\sqrt{T}}(v-x).  \label{eq:T218}
\end{equation}
To estimate $III$ note that the same inequalities (\ref{eq:T209}), (%
\ref{eq:T210}),(\ref{eq:T217}) hold for $H+M_{h,3}^{\prime \prime
},$ namely,
\begin{eqnarray} &&
\left| D_{v}^{a}D_{x}^{b}(H+M_{h,3}^{\prime \prime })(jh,kh,x,v)\right| \leq
C\rho ^{-1-\left| a\right| -\left| b\right| }\zeta _{\rho }(v-x)
\label{eq:T219}
\\ &&
\left| D_{x}^{b}(H+M_{h,3}^{\prime \prime })(jh,kh,x,x+v)\right| \leq C\rho
^{-1}\zeta _{\rho }(v-x)  \label{eq:T220}
\\ &&
\left| D_{v}^{a}D_{x}^{b}[\widetilde{p}\otimes
_{h}(H+M_{h,3}^{\prime \prime })^{(r)}](jh,kh,x,v)\right|
\label{eq:T221}
\\  \nonumber && \qquad
\leq C^{r}B(1,\frac{1}{2})\times ...\times
B(\frac{r+1}{2},\frac{1}{2})\rho ^{r-\left| a\right| -\left|
b\right| }\zeta _{\rho }(v-x).
\end{eqnarray}
To prove (\ref{eq:T219})- (\ref{eq:T221}) it is enough to get the
corresponding estimates for summands in $M_{h,3}^{\prime \prime }$
(see (\ref{eq:012h})). These estimates can be proved by the same
arguments as used in the proof of (\ref{eq:T209}), (\ref{eq:T210})
and (\ref{eq:T217}). To estimate $III$ we have now to estimate
$\widetilde{p}_{1,r}\otimes _{h}M_{h,3}^{\prime \prime }$ and
$\widetilde{p}_{2,r}\otimes _{h}M_{h,3}^{\prime \prime }$ where
\[
\widetilde{p}_{2,r}=\widetilde{p}\otimes _{h}(H+M_{h,3}^{\prime \prime
})^{(r)}.
\]
Using integration by parts and inequality (\ref{eq:T221}), we
obtain for $\widetilde{p}_{2,r-1}\otimes _{h}M_{h,3}^{\prime
\prime }$ the same estimate as for $\widetilde{p}_{1,r}\otimes
_{h}M_{h,3}^{\prime \prime }$. It holds for $i=1,2$
\begin{eqnarray*}
&&\left| \widetilde{p}_{2,r}\otimes _{h}M_{h,3}^{\prime \prime
}(0,kh,x,y)\right|
\\ &&
\qquad \leq C^{r}h^{3/2-2\varepsilon }B(1,\frac{1}{2})\times ...\times B(\frac{r+1}{%
2},\frac{1}{2})B(\frac{r}{2},\varepsilon )(kh)^{\varepsilon +\frac{r-2}{2}%
}\zeta _{\sqrt{kh}}(y-x).\end{eqnarray*} Hence it holds for $r\geq
2$
\begin{eqnarray}  \label{eq:T223} &&
\left| III\right| \leq \left| \widetilde{p}_{1,r-1}\otimes
_{h}M_{h,3}^{\prime \prime }(0,kh,x,y)\right| +\left| \widetilde{p}%
_{2,r-1}\otimes _{h}M_{h,3}^{\prime \prime }(0,kh,x,y)\right|
\\ \nonumber &&
\leq C^{r}h^{3/2-2\varepsilon }B(1,\frac{1}{2})\times ...\times B(\frac{r}{2}%
,\frac{1}{2})B(\frac{r-1}{2},\varepsilon )(kh)^{\varepsilon +\frac{r-3}{2}%
}\zeta _{\sqrt{kh}}(y-x).
\end{eqnarray}
From (\ref{eq:T203aa}), (\ref{eq:T218}) and (\ref{eq:T223}) we get
for $r\geq 2$
\[
\left| \Gamma _{r}(0,kh,x,y)\right| \leq C^{r}h^{3/2-2\varepsilon
}B(\varepsilon ,\varepsilon )B(\varepsilon ,\frac{1}{2})\times ...\times
B(\varepsilon +\frac{r-2}{2},\frac{1}{2})(kh)^{\varepsilon +\frac{r-3}{2}%
}\zeta _{\sqrt{kh}}(v-x).
\]
In particular
\begin{equation}
\left| \Gamma _{r}(0,T,x,y)\right| \leq h^{3/2-2\varepsilon }\frac{\Gamma
^{3}(\varepsilon )}{\Gamma (2\varepsilon )}\frac{C^{r}}{\Gamma (\varepsilon +%
\frac{r-1}{2})}T^{\varepsilon +\frac{r-3}{2}}\zeta
_{\sqrt{kh}}(v-x) \label{eq:T224}
\end{equation}
for any $\varepsilon \in (0,1/4)$ and $r\geq 2$. Now we shall
estimate the left hand side of (\ref{eq:T203}). Denote the
expression under the sign of the absolute value in (\ref{eq:T203})
by $\digamma _{r}$. Note that $\digamma _{0}=\digamma _{1}=0$. For
$r\geq 2$ we make use of the following recurrence formula
\begin{eqnarray*}
\digamma _{r}&=&\digamma _{r-1}\otimes _{h}H+\left[
\widetilde{p}\otimes
_{h}(H+M_{h,1}+M_{h,2}+\sqrt{h}N_{1})^{(r-1)}\right.
\\ && \qquad
\left. -\widetilde{p}\otimes _{h}(H+M_{h,1}+\sqrt{h}N_{1})^{(r-1)}\right]
\otimes _{h}(M_{h,1}+M_{h,2}+\sqrt{h}N_{1})
\\ &&
+\left[ \widetilde{p}\otimes _{h}(H+M_{h,1}+\sqrt{h}N_{1})^{(r-1)}-%
\widetilde{p}\otimes _{h}(H+M_{h,1})^{(r-1)}\right] \otimes
_{h}M_{h,2}\\ &=&I+II+III.
\end{eqnarray*}
We start again from the estimation of
\[
A_{r-1}=\widetilde{p}\otimes _{h}(H+M_{h,1}+M_{h,2}+\sqrt{h}N_{1})^{(r-1)}-%
\widetilde{p}\otimes _{h}(H+M_{h,1}+\sqrt{h}N_{1})^{(r-1)}.
\]
For $r=2$ we have $A_{1}=(\widetilde{p}\otimes
_{h}M_{h,2})(0,kh,x,y)$. It is enough to estimate
\[
J_{3}=h\sum_{i=0}^{k-2}h\int \widetilde{p}(0,ih,x,v)(f(ih,v)-f(ih,y))D_{v}^{%
\nu }\widetilde{p}_{h}((i+1)h,kh,v,y)dv
\]
for $\left| \nu \right| =4$. Analogously to (\ref{eq:15}) we
obtain that
\[
\left| J_{3}\right| \leq Ch^{1-\varepsilon }(kh)^{-1/2+\varepsilon }B(\frac{1%
}{2},\varepsilon )\zeta _{\sqrt{kh}}^{S}(y-x)
\]
and, hence,
\begin{equation}
\left| A_{1}\right| \leq Ch^{1-\varepsilon }(kh)^{-1/2+\varepsilon }B(\frac{1%
}{2},\varepsilon )\zeta _{\sqrt{kh}}^{S}(y-x).  \label{eq:T224a}
\end{equation}
For $r\geq 3$ we use the recurrence relation
\begin{eqnarray} \label{eq:T225}
 A_{r-1}&=&A_{r-2}\otimes _{h}(H+M_{h,1}+M_{h,2}+\sqrt{h}N_{1})
\\ \nonumber && \qquad
+\left[ \widetilde{p}\otimes
_{h}(H+M_{h,1}+\sqrt{h}N_{1})^{(r-2)}\right] \otimes
_{h}M_{h,2}
\\ \nonumber &=&I^{\prime }+II^{\prime }.
\end{eqnarray}
From (\ref{eq:T224a}) and (\ref{eq:T225}) we obtain for $r=3$%
\begin{eqnarray} \label{eq:T225a}
\left| I^{\prime }\right| &\leq& Ch^{1-\varepsilon
}B(\frac{1}{2},\varepsilon )\zeta
_{\sqrt{kh}}^{S}(y-x)\sum_{i=0}^{k-2}h(ih)^{\varepsilon
-1/2}(kh-ih)^{-1/2}
\\ \nonumber
&\leq& Ch^{1-\varepsilon }B(\frac{1}{2},\varepsilon
)B(\frac{1}{2},\varepsilon +\frac{1}{2})(kh)^{\varepsilon }\zeta
_{\sqrt{kh}}^{S}(y-x)  .
\end{eqnarray}
To estimate $II^{\prime }$ \ we use the inequality for $r\geq 3$
\begin{eqnarray} \label{eq:T226}
&& \left| D_{v}^{a}D_{x}^{b}[\widetilde{p}\otimes _{h}(H+M_{h,1}+\sqrt{h}%
N_{1})^{(r-2)}](jh,kh,x,v)\right|
\\ \nonumber && \qquad
\leq C^{r}B(1,\frac{1}{2})\times ...\times
B(\frac{r-1}{2},\frac{1}{2})\rho ^{r-2-\left| a\right| -\left|
b\right| }\zeta _{\rho }(v-x). \end{eqnarray} This inequality is a
direct consequence of (\ref{eq:T217}). We have
\begin{equation}
\left| II^{\prime }\right| \leq Ch^{1-\varepsilon }B(1,\varepsilon
)(kh)^{\varepsilon }\zeta _{\sqrt{kh}}^{S}(y-x).  \label{eq:T226a}
\end{equation}
Comparing (\ref{eq:T225a}) and (\ref{eq:T226a}) we obtain that
$\left|
A_{2}\right| \leq $\ $C^{2}h^{1-\varepsilon }B(\frac{1}{2},\varepsilon )B(%
\frac{1}{2},\varepsilon +\frac{1}{2})(kh)^{\varepsilon }\zeta _{\sqrt{kh}%
}^{S}(y-x)$. By induction we easily get that for $r\geq 2$
\begin{equation}
\left| A_{r-1}(0,kh,x,y)\right| \leq C^{r}h^{1-\varepsilon }B(\frac{1}{2}%
,\varepsilon )\times ...\times B(\frac{1}{2},\varepsilon +\frac{r-2}{2}%
)(kh)^{\varepsilon +\frac{r-3}{2}}\zeta _{\sqrt{kh}}^{S}(y-x).
\label{eq:T227}
\end{equation}
For a bound of
$
A_{r-1}\otimes _{h}(M_{h,1}+M_{h,2}+\sqrt{h}N_{1})
$
it is enough to estimate
\begin{eqnarray}  \label{eq:T228}
&& J_{4}=h^{1/2}\sum_{i=0}^{k-2}h\int
A_{r-1}(0,ih,x,v)(f(ih,v)-f(ih,y))\\ \nonumber && \qquad
 D_{v}^{\nu }\widetilde{p}%
_{h}((i+1)h,kh,v,y)dv \mbox{ for }\left| \nu \right| =3, \\
&& \label{eq:T229}
J_{5}=h\sum_{i=0}^{k-2}h\int A_{r-1}(0,ih,x,v)(f(ih,v)-f(ih,y))
\\ \nonumber && \qquad D_{v}^{\nu }%
\widetilde{p}_{h}((i+1)h,kh,v,y)dv \mbox{ for } \left| \nu \right|
=4,
\\ \label{eq:T230} &&
J_{6}=h^{1/2}\sum_{i=0}^{k-2}h\int
A_{r-1}(0,ih,x,v)(f(ih,v)-f(ih,y))(kh-ih)
\\ \nonumber && \qquad
D_{v}^{\nu +e_{p}+e_{q}}\widetilde{p%
}(ih,kh,v,y)dv \mbox{ for }\left| \nu \right| =3.\end{eqnarray}
 It
follows from (\ref{eq:T227}) that
\begin{eqnarray*}
&&\left| J_{4}\right| \leq C^{r}h^{3/2-2\varepsilon
}B(\frac{1}{2},\varepsilon )\times ...\times
B(\frac{1}{2},\varepsilon +\frac{r-2}{2})\\ && \qquad \times
B(\varepsilon
,\varepsilon +\frac{r-1}{2})(kh)^{2\varepsilon +\frac{r-3}{2}}\zeta _{\sqrt{%
kh}}^{S}(y-x).
\end{eqnarray*}
Clearly, the same estimate holds for $J_{5}$ and $J_{6}$. Thus we
obtain
\begin{equation}
\left| II\right| \leq C^{r}h^{3/2-2\varepsilon }B(\frac{1}{2},\varepsilon
)\times ...\times B(\frac{1}{2},\varepsilon +\frac{r-2}{2})B(\varepsilon
,\varepsilon +\frac{r-1}{2})(kh)^{2\varepsilon +\frac{r-3}{2}}\zeta _{\sqrt{%
kh}}^{S}(y-x).  \label{eq:T231}
\end{equation}
Now we estimate $III$. We denote
\[
B_{r-1}=\widetilde{p}\otimes _{h}(H+M_{h,1}+\sqrt{h}N_{1})^{(r-1)}-%
\widetilde{p}\otimes _{h}(H+M_{h,1})^{(r-1)}.
\]
Using the recurrence equation $B_{0}=0$ and
\[
B_{r-1}=B_{r-2}\otimes
_{h}(H+M_{h,1}+\sqrt{h}N_{1})+\widetilde{p}\otimes
_{h}(H+M_{h,1})^{(r-2)}\otimes _{h}\sqrt{h}N_{1}
\]
we obtain that
\begin{equation}
III=\sum_{l=0}^{r-2}\widetilde{p}_{3,l}\otimes _{h}\sqrt{h}N_{1}\otimes
_{h}(H+M_{h,1}+\sqrt{h}N_{1})^{(r-l-2)}\otimes _{h}M_{h,2}(0,T,x,y),
\label{eq:T232}
\end{equation}
where $\widetilde{p}_{3,l}=\widetilde{p}\otimes
_{h}(H+M_{h,1})^{(l)}.$ To estimate $III$ \ it's enough to
estimate a typical \ term \ in the last sum. Thus we need a bound
for
\begin{eqnarray}  \label{eq:T233}
&&h^{3/2}\sum_{k=0}^{n-2}h\int \left\{ \sum_{j=0}^{k-1}h\left[
\int \sum_{i=0}^{j-1}h\int
\widetilde{p}_{3,l}(0,ih,x,w)(jh-ih)\right. \right.
\\ \nonumber && \qquad
\left.  \times D_{w}^{\mu
+e_{n}+e_{m}}\widetilde{p}(ih,jh,w,z)(g(ih,w)-g(ih,z))dw\right]
\\ \nonumber && \qquad \left.
\times  (H+M_{h,1}+\sqrt{h}%
N_{1})^{(r-l-2)}(ih,kh,z,v)dz\right\} (f(kh,v)-f(kh,y))\\
\nonumber && \qquad \times D_{v}^{\nu
}\widetilde{p}_{h}((k+1)h,T,v,y)dv.
\end{eqnarray}
For an estimate of (\ref{eq:T233}) we apply two times integration
by parts in the internal integral$\int ...dw$ and then we make two
times an integration by parts in $\int ...dv$. We use also the
following estimates for $0\leq l\leq r-3$
\begin{eqnarray} \label{eq:T234}
&& \left| D_{w}^{a}D_{x}^{b}\widetilde{p}_{3,l}(0,ih,x,w)\right|
\\ \nonumber && \qquad
\leq C^{l}B(1,\frac{1}{2})\times ...\times B(\frac{l+1}{2},\frac{1}{2})(ih)^{%
\frac{l-\left| a\right| -\left| b\right| }{2}}\zeta _{\sqrt{ih}}(w-x),
\\ \label{eq:T235} &&
\left| D_{v}^{a}D_{z}^{b}(H+M_{h,1}+\sqrt{h}N_{1})^{(r-l-2)}(ih,kh,z,v)%
\right|
\\ \nonumber && \qquad
\leq C^{r-l-2}B(\frac{1}{2},\frac{1}{2})\times ...\times B(\frac{1}{2},\frac{%
r-l-3}{2})(kh-ih)^{\frac{r-l-4-\left| a\right| -\left| b\right| }{2}}\zeta _{%
\sqrt{kh-ih}}(v-z).
\end{eqnarray}
We put here $B(\frac{1}{2},0)=1$. This gives the following
estimate for
$0\leq l\leq r-3,r\geq 2$%
\begin{eqnarray} \label{eq:T236}
&&
\left| \widetilde{p}_{3,l}\otimes _{h}\sqrt{h}N_{1}\otimes _{h}(H+M_{h,1}+%
\sqrt{h}N_{1})^{(r-l-2)}\otimes _{h}M_{h,2}(0,T,x,y)\right|
\\ \nonumber && \qquad
\leq C^{r}h^{3/2-3\varepsilon }\frac{\Gamma (\varepsilon )}{\Gamma
(3\varepsilon +\frac{r-1}{2})}T^{3\varepsilon +\frac{r-3}{2}}\zeta _{\sqrt{T}%
}(y-x).  \end{eqnarray} For $l=r-2$ we have to estimate
\[
\widetilde{p}_{3,r-2}\otimes _{h}\sqrt{h}N_{1}{}_{h}\otimes
_{h}M_{h,2}(0,T,x,y).
\]
This is a finite sum of terms corresponding to the different summands in $%
N_{1,h}$ and $M_{h,2}$. A typical term can be bounded by
\begin{eqnarray} \label{eq:T237}
&&
h^{3/2}\sum_{k=0}^{n-2}h\int \left\{ \sum_{j=0}^{k-1}h\int \widetilde{p}%
_{3,r-2}(0,jh,x,w)(kh-jh)D_{w}^{\mu +e_{n}+e_{m}}\widetilde{p}%
(jh,kh,w,v)\right.
\\ \nonumber && \qquad
\left. \times (g(jh,w)-g(jh,v))dw\right\} (f(kh,v)-f(kh,y))D_{v}^{\nu }%
\widetilde{p}_{h}((k+1)h,T,v,y)dv . \end{eqnarray}  We apply again
integration by parts and after direct calculations we obtain the
following estimate for $r \geq 2$
\begin{eqnarray} \label{eq:T238}
&& \left| \widetilde{p}_{3,r-2}\otimes
_{h}\sqrt{h}N_{1}{}_{h}\otimes _{h}M_{h,2}(0,T,x,y)\right|
\\ \nonumber && \qquad
\leq C^{r}h^{3/2-3\varepsilon }\frac{\Gamma (\varepsilon +\frac{r-2}{2})}{%
\Gamma (3\varepsilon +\frac{r-2}{2})}\Gamma ^{2}(\varepsilon )\frac{1}{%
\Gamma (\frac{r}{2})}T^{3\varepsilon +\frac{r-4}{2}}\zeta _{\sqrt{T}%
}(y-x).  \end{eqnarray} The inequalities (\ref{eq:T202}),
(\ref{eq:T203}) follow from (\ref {eq:T224}), (\ref{eq:T236}) and
(\ref{eq:T238}).

\textit{Asymptotic treatment of the term $T_{3}$. }We will show that
\begin{eqnarray} \label{eq:T301}
&& \left| T_{3}-\left[ \sum_{r=0}^{\infty }\tilde{p}\otimes
_{h}(H+A)^{(r)}(0,T,x,y)-\sum_{r=0}^{\infty }\tilde{p}\otimes
_{h}H^{(r)}(0,T,x,y)\right] \right|
\\ \nonumber && \qquad
\leq Chn^{-\delta }\zeta _{\sqrt{T}}(y-x),  \end{eqnarray}
where $A=M_{h}^{\prime \prime }-M_{h}=-\frac{h}{2}(L_{\ast }^{2}-2L%
\widetilde{L}+\widetilde{L}^{2})\lambda (x)$. Denote
\begin{eqnarray*}
C_{r}&=&\tilde{p}\otimes _{h}(H+M_{h}^{\prime \prime }+\sqrt{h}%
N_{1})^{(r)}(0,T,x,y)
\\ && \nonumber \qquad
-\tilde{p}\otimes _{h}(H+M_{h}+\sqrt{h}N_{1})^{(r)}(0,T,x,y)
\\ && \nonumber \qquad
-[\tilde{p}\otimes _{h}(H+A)^{(r)}-\tilde{p}\otimes
_{h}H^{(r)}](0,T,x,y). \end{eqnarray*} Analogously to
(\ref{eq:T203aa}) we have the recurrence relation
\begin{eqnarray} \label{T302}
C_{r}&=&C_{r-1}\otimes _{h}H+\left[ \tilde{p}\otimes
_{h}(H+M_{h}^{\prime \prime }+\sqrt{h}N_{1})^{(r-1)}\right.
\\ && \nonumber \qquad
\left. -\tilde{p}\otimes
_{h}(H+M_{h}+\sqrt{h}N_{1})^{(r-1)}\right] \otimes
_{h}(M_{h}^{\prime \prime }+\sqrt{h}N_{1})
\\ && \nonumber \qquad
+\left[ \tilde{p}\otimes _{h}(H+M_{h}+\sqrt{h}N_{1})^{(r-1)}-\tilde{p}%
\otimes _{h}(H+A)^{(r-1)}\right] \otimes _{h}A
\\  \nonumber \qquad&=&I+II+III.
\end{eqnarray}
Denote
\[
D_{r-1}=\tilde{p}\otimes _{h}(H+M_{h}+\sqrt{h}N_{1})^{(r-1)}-\tilde{p}%
\otimes _{h}(H+A)^{(r-1)}.
\]
Clearly
\[
D_{r-1}=D_{r-2}\otimes _{h}(H+M_{h}+\sqrt{h}N_{1})+\tilde{p}_{h}\otimes
_{h}(H+A)^{(r-2)}\otimes _{h}(M_{h}-A+\sqrt{h}N_{1}).
\]
Iterating we obtain
\begin{eqnarray} \label{T303}
&&III=D_{r-1}\otimes _{h}A\\ \nonumber &&=
\sum_{l=0}^{r-2}\widetilde{p}_{4,l}\otimes _{h}(M_{h}-A+\sqrt{h}%
N_{1})\otimes _{h}(H+M_{h}+\sqrt{h}N_{1})^{(r-l-2)}\otimes _{h}A(0,T,x,y),
\end{eqnarray}
where $\widetilde{p}_{4,l}=\widetilde{p}\otimes _{h}(H+A)^{(l)}.$ \ The sum (%
\ref{T303}) \ can be estimated exactly in the same way as the sum \ (\ref
{eq:T232}). This gives
\begin{equation}
\left| III\right| \leq C(\varepsilon )h^{3/2-2\varepsilon }\frac{C^{r}}{%
\Gamma (\frac{r-1}{2})}T^{3\varepsilon +\frac{r-4}{2}}\zeta _{\sqrt{T}%
}(v-x),r=2,3...  \label{T304}
\end{equation}
To estimate $II$ \ denote
\[
E_{r-1}=\tilde{p}\otimes _{h}(H+M_{h}^{\prime \prime }+\sqrt{h}%
N_{1})^{(r-1)}-\tilde{p}\otimes _{h}(H+M_{h}+\sqrt{h}N_{1})^{(r-1)}.
\]
For $r=2$ \ we have $E_{1}=\tilde{p}\otimes _{h}A$ \ and analogously to (%
\ref{eq:T206})
\begin{equation}
\left| E_{1}\right| \leq Ch^{1-\varepsilon }(kh)^{\varepsilon -1/2}B(\frac{1%
}{2},\varepsilon )\zeta _{\sqrt{kh}}^{S}(y-x).  \label{T305}
\end{equation}
Analogously to (\ref{eq:T208}) \ for $r\geq 3$ \ we use the recurrence
relation \
\begin{eqnarray} \label{T306}
E_{r-1}&=&E_{r-2}\otimes _{h}(H+M_{h}^{\prime \prime }+\sqrt{h}N_{1})+[\tilde{p%
}\otimes _{h}(H+M_{h}+\sqrt{h}N_{1})^{(r-2)}]\otimes _{h}A
\\ \nonumber
&=&I^{\prime }+II^{\prime }.
\end{eqnarray}
Both terms in (\ref{T306}) \ are analogous to the corresponding terms in (%
\ref{eq:T225}) and may be estimated analogously. This gives the
following estimates for $r\geq 3$
\begin{eqnarray} \nonumber
\left| E_{r-1}\right| &\leq& C^{r}h^{1-\varepsilon
}B(\frac{1}{2},\varepsilon )\times ...\times
B(\frac{1}{2},\varepsilon +\frac{r-2}{2})(kh)^{\varepsilon
+\frac{r-3}{2}}\zeta _{\sqrt{kh}}^{S}(y-x),\\ \label{T307}
\left| II\right| &=&\left| E_{r-1}\otimes _{h}(M_{h}^{\prime \prime }+\sqrt{h}%
N_{1})(0,T,x,y)\right|\\ \nonumber
&\leq& C(\varepsilon )h^{3/2-2\varepsilon }\frac{C^{r}}{\Gamma (\frac{r-1}{2})}%
T^{3\varepsilon +\frac{r-4}{2}}\zeta _{\sqrt{T}}(v-x).
\end{eqnarray}
The desired estimate (\ref{eq:T301}) follows from (\ref{T302}),
(\ref{T304}) and (\ref{T307}).

\textit{Asymptotic treatment of the term $T_{4}$.} We will show
that
\begin{eqnarray} \label{eq:3a}
T_{4}&=&\sum_{r=1}^{\infty }\tilde{p}\otimes _{h}H^{(r)}(0,T,x,y)
\\ \nonumber && \qquad
-\sum_{r=1}^{\infty }\tilde{p}\otimes
_{h}[H+hN_{2}]^{(r)}(0,T,x,y)+R_{h}^{\ast }(x,y),
\end{eqnarray}
with $N_{2}(s,t,x,y)=(L-\widetilde{L})\widetilde{\pi
}_{2}(s,t,x,y),\left| R_{h}^{\ast }(x,y)\right| \leq Chn^{-\delta
}\zeta _{\sqrt{T}}^{S}(y-x)$  for $\delta >0$ small enough  and
with a constant $C$  depending on $\delta $. For the proof of
(\ref{eq:3a}) it suffices to show that for $\delta $ small enough
\begin{eqnarray} \label{eq:3b} &&
\left| \sum_{r=1}^{n}\tilde{p}\otimes _{h}(H+M_{h}+\sqrt{h}%
N_{1}+hN_{2})^{(r)}(0,T,x,y)\right.
\\ \nonumber && \qquad
\left. -\sum_{r=1}^{n}\tilde{p}\otimes
_{h}(K_{h}+M_{h})^{(r)}(0,T,x,y)\right|
\\ \nonumber && \leq \left[ \sum_{k=1}^{n}\frac{C^{k}%
}{\Gamma (\frac{k}{2})}\right] hn^{-\delta }\zeta
_{\sqrt{T}}^{S}(y-x),
\\ \label{eq:3c} &&
\left| \sum_{r=1}^{n}\tilde{p}\otimes _{h}(H+M_{h}+\sqrt{h}%
N_{1})^{(r)}(0,T,x,y)\right.
\\ \nonumber && \qquad
-\sum_{r=1}^{n}\tilde{p}\otimes _{h}(H+M_{h}+\sqrt{h}%
N_{1}+hN_{2})^{(r)}(0,T,x,y)
\\ \nonumber && \qquad
\left. -\left[ \sum_{r=1}^{n}\tilde{p}\otimes
_{h}H^{(r)}(0,T,x,y)-\sum_{r=1}^{n}\tilde{p}\otimes
_{h}[H+hN_{2}]^{(r)}(0,T,x,y)\right] \right|
\\ \nonumber &&
\leq \left[ \sum_{k=1}^{n}\frac{C^{k}}{\Gamma
(\frac{k}{2})}\right] Chn^{-\delta }\zeta _{\sqrt{T}}^{S}(y-x).
\end{eqnarray}
Denote $D_{3,0}\equiv 0$ and
\begin{eqnarray*}
D_{3,m}(0,jh,x,y)&=&\sum_{r=1}^{m}\tilde{p}\otimes
_{h}(K_{h}+M_{h})^{(r)}(0,jh,x,y)
\\ \nonumber &&
-\sum_{r=1}^{m}\tilde{p}\otimes _{h}(H+M_{h}+\sqrt{h}%
N_{1}+hN_{2})^{(r)}(0,jh,x,y).
\end{eqnarray*}
Then (\ref{eq:3b}) can be rewritten as
\[
\left| D_{3,n}(0,T,x,y)\right| \leq Chn^{-\delta }\zeta _{\sqrt{T}%
}^{S}(y-x).
\]
We now make iterative use of
\begin{equation}
D_{3,m}=D_{3,m-1}\otimes _{h}(H+M_{h}+\sqrt{h}N_{1}+hN_{2})+g_{m-1},
\label{eq:3d}
\end{equation}
for $m=1,2,...,$ \ where
\begin{eqnarray*}
g_{m}(0,jh,x,y)&=&-\left[ \sum_{r=0}^{m}\tilde{p}\otimes
_{h}(K_{h}+M_{h})^{(r)}\right] \otimes _{h}(H-K_{h}+\sqrt{h}%
N_{1}+hN_{2})(0,jh,x,y)
\\
&=&S_{h,m}\otimes _{h}(L-\widetilde{L})d_{h}(0,jh,x,y)
\end{eqnarray*}
with
\begin{eqnarray*}
g_{0}(0,jh,x,y)&=&-\tilde{p}\otimes _{h}(H-K_{h}+\sqrt{h}%
N_{1}+hN_{2})(0,jh,x,y),
\\
d_{h}&=&\widetilde{p}_{h}-\widetilde{p}-\sqrt{h}\widetilde{\pi }_{1}-h%
\widetilde{\pi }_{2},
\\
S_{h,m}(0,ih,x,y)&=&\sum_{r=0}^{m}\tilde{p}\otimes
_{h}(K_{h}+M_{h})^{(r)}(0,ih,x,y).
\end{eqnarray*}
Iterative application of (\ref{eq:3d}) gives
\[
D_{3,n}(0,T,x,y)=\sum_{r=0}^{n-1}g_{r}\otimes _{h}(H+M_{h}+\sqrt{h}%
N_{1}+hN_{2})^{(n-r-1)}(0,T,x,y).
\]
To prove (\ref{eq:3b}) we will show that
\begin{eqnarray} \label{eq:3e}
&&\left| g_{r}\otimes _{h}(H+M_{h}+\sqrt{h}N_{1}+hN_{2})^{(n-r-1)}(0,T,x,y)%
\right|
\\ \nonumber && \qquad
\leq \frac{C^{n-r}}{\Gamma (\frac{n-r}{2})}hn^{-\delta }\zeta _{\sqrt{T}%
}^{S}(y-x).  \label{eq:3e}
\end{eqnarray}
For this purpose we decompose the left handside of (\ref{eq:3e})
into four terms
\begin{eqnarray*}
a_{r,1}&=&\sum_{0\leq i\leq n/2}h\int g_{r}(0,ih,x,u)(H+M_{h}+\sqrt{h}%
N_{1}+hN_{2})^{(n-r-1)}(ih,T,u,y)du,
\\
a_{r,2}&=&\sum_{n/2<i\leq n}h^{2}\sum_{0\leq k\leq i/2}\int \int
S_{h,r}(0,kh,x,v)(L-\widetilde{L})d_{h}(kh,ih,v,u)
\\ &&
\times (H+M_{h}+\sqrt{h}N_{1}+hN_{2})^{(n-r-1)}(ih,T,u,y)dvdu,
\\
a_{r,3}&=&\sum_{n/2<i\leq n}h^{2}\sum_{i/2<k\leq i-n^{\delta
^{\prime }}}\int \int
(L^{T}-\widetilde{L}^{T})S_{h,r}(0,kh,x,v)d_{h}(kh,ih,v,u)
\\ &&
\times (H+M_{h}+\sqrt{h}N_{1}+hN_{2})^{(n-r-1)}(ih,T,u,y)dvdu,
\\
a_{r,4}&=&\sum_{n/2<i\leq n}h^{2}\sum_{i-n^{\delta ^{\prime
}}<k\leq i-1}\int \int
(L^{T}-\widetilde{L}^{T})S_{h,r}(0,kh,x,v)d_{h}(kh,ih,v,u)
\\ &&
\times (H+M_{h}+\sqrt{h}N_{1}+hN_{2})^{(n-r-1)}(ih,T,u,y)dvdu.
\end{eqnarray*}
Here $L^{T}$ and $\widetilde{L}^{T}$ denote the adjoint operators
of $L$ and $\widetilde{L}$, and $\delta ^{\prime }$ satisfies
inequalities $2\varkappa <\delta ^{\prime
}<\frac{3}{5}(1-\varkappa ),$ where $\varkappa $ is defined in
(B2). For the proof of
(\ref{eq:3e}) it suffices to show for $l=1,2,3,4$%
\begin{eqnarray} \label{eq:3f}
\left| a_{r,l}\right| &\leq& hn^{-\delta
}C^{n-r}B(1,\frac{1}{2})\times ...\times
B(\frac{n-r-1}{2},\frac{1}{2})\zeta _{\sqrt{T}}^{S}(y-x)
\\ \nonumber
&\leq & \left[ \frac{C^{n-r}}{\Gamma (\frac{n-r}{2})}\right]
hn^{-\delta }\zeta _{\sqrt{T}}^{S}(y-x)
\end{eqnarray}
for some $\delta >0.$

\textit{Proof \ of \ (\ref{eq:3f}) for }$l=2.$ Note that \ $k\leq
i/2,i>n/2$ \ imply \ $ih-kh\geq \frac{T}{4}$. The claim follows
from the inequalities
\begin{eqnarray} \label{eq:3g}
&& \max \{\left| K_{h}(ih,jh,x,y)\right| ,\left|
M_{h}(ih,jh,x,y)\right| ,\left| \sqrt{h}N_{1}(ih,jh,x,y)\right|
\\ \nonumber && \qquad
\left| hN_{2}(ih,jh,x,y)\right| ,\left| H(ih,jh,x,y)\right| \}
\\ \nonumber &&
\leq C\rho ^{-1}\zeta _{\rho }(y-x)\mbox{ with }\rho ^{2}=jh-ih
\mbox { for }0\leq i<j\leq n,
\end{eqnarray}
%\\ \label{eq:3h} &&
\begin{eqnarray} \label{eq:3h} &&
\left| S_{h,m}(0,kh,x,v)\right| \leq C\zeta
_{\sqrt{kh}}^{S-2}(v-x),
\\ \label{eq:03a} &&
\left| (L-\widetilde{L})d_{h}(kh,ih,v,u)\right| \leq
Ch^{3/2}(ih-kh)^{-2}\zeta
_{\sqrt{ih-kh}}^{S-8}(u-v)
\\ \nonumber &&
=O(hn^{-1/2+3/2\varkappa })\zeta
_{\sqrt{ih-kh}}^{S-8}(u-v),
\\ \label{eq:03b} &&
\left| (H+M_{h}+\sqrt{h}N_{1}+hN_{2})^{(n-r-1)}(ih,T,u,y)\right|
\\ \nonumber &&
\leq C^{n-r}\rho ^{n-r-3}B(\frac{1}{2},\frac{1}{2})\times ...\times B(\frac{%
n-r-2}{2},\frac{1}{2})\zeta _{\sqrt{T-ih}}^{S-2}(y-u)
\\ \nonumber &&
\leq \left[ \frac{C^{n-r}}{\Gamma (\frac{n-r-1}{2})}\right]
(T-ih)^{-1/2}\zeta _{\sqrt{T-ih}}^{S-2}(y-u)
\end{eqnarray}
for $n-r-3=-1,0,1,...,n-3$ with $\rho ^{2}=T-ih$. We put
$B(\frac{1}{2},0)=1)$.
Inequality (\ref{eq:3g}) follows from the definitions of the functions $%
K_{h},...,H$. Inequalities (\ref{eq:3h}) and (\ref{eq:03b}) can be
proved by the same method as used in the proof of Theorem 2.3 in
Konakov and Mammen (2002) (pp. 282 - 284). Inequality
(\ref{eq:03a}) follows from the inequality $ih-kh\geq
\frac{T}{4}$, Lemma 5 and the arguments used in the proof of Lemma
7.

\textit{Proof \ of \ (\ref{eq:3f}) for }$l=3.$ Note that $n/2<i,k>i/2$ \
imply $\ kh>\frac{T}{4}.$ We use the following inequalities
\begin{eqnarray} \label{eq:03ba} &&
\left| d_{h}(kh,ih,v,u)\right| \leq Ch^{3/2}(ih-kh)^{-3/2}\zeta _{\sqrt{ih-kh%
}}^{S-6}(u-v),
\\ \nonumber &&
\left| (L^{T}-\widetilde{L}^{T})S_{h,r}(0,kh,x,v)\right| \leq CT^{-1}\zeta _{%
\sqrt{kh}}^{S-2}(v-x),
\\ \label{eq:03c} &&
\left| h\sum_{i/2<k\leq i-n^{\delta ^{\prime }}}\int (L^{T}-\widetilde{L}%
^{T})S_{h,r}(0,kh,x,v)d_{h}(kh,ih,v,u)dv\right|
\\ \nonumber &&
\leq Ch^{3/2}T^{-1}\sum_{i/2<k\leq i-n^{\delta ^{\prime }}}h\frac{1}{%
(ih-kh)^{3/2}}\zeta _{\sqrt{ih}}(u-x)
\\ \nonumber &&
\leq Ch^{3/2}T^{-1}\int_{ih/2}^{ih-n^{\delta ^{\prime }}h}\frac{du}{%
(ih-u)^{3/2}}\zeta _{\sqrt{ih}}(u-x)\leq Ch^{3/2}T^{-3/2}n^{(1-\delta
^{\prime })/2}\zeta _{\sqrt{ih}}(u-x)
\\ \nonumber &&
\leq Chn^{-\delta ^{\prime \prime }}\zeta _{\sqrt{ih}}(u-x),
\end{eqnarray}
where $\delta ^{\prime \prime }=\delta ^{\prime }/2-\varkappa >0.$
Claim (\ref {eq:3f}) for $l=3$ now follows from (\ref{eq:03c}) and
(\ref{eq:03b}).

\textit{Proof of  (\ref{eq:3f}) for }$l=4$. For $i-n^{\delta
^{\prime }}<k\leq i-1,n/2<i$ we have $ih>T/2,$ $kh>T/3,$
$(i-k)<n^{\delta ^{\prime }}$ for sufficiently large $n$. The
integral
\[
\int (L^{T}-\widetilde{L}^{T})S_{h,r}(0,kh,x,v)\widetilde{p}%
_{h}(kh,ih,v,u)dv
\]
is a finite sum of integrals. We show how to estimate a typical
term of this sum. The other terms can be estimated analogously. We
consider for
fixed $j,l$%
\begin{eqnarray} \label{eq:03d} &&
\int \frac{\partial ^{2}S_{h,r}(0,kh,x,v)}{\partial v_{j}\partial v_{l}}%
(\sigma _{jl}(kh,v)-\sigma
_{jl}(kh,u))h^{-d/2}
\\ \nonumber &&
\qquad \times
q^{(i-k)}[kh,u,h^{-1/2}(u-v-h\sum_{l=k}^{i-1}m(lh,u))]dv
\\ \nonumber &&
=\int \frac{\partial ^{2}S_{h,r}(0,kh,x,v)}{\partial v_{j}\partial v_{l}}%
\mid _{v=u^{\ast }-\sqrt{h}w}
\\ \nonumber && \qquad
\times \lbrack \sigma _{jl}(kh,u^{\ast }-\sqrt{h}w)-\sigma
_{jl}(kh,u)]q^{(i-k)}(kh,u,w)dw,
\end{eqnarray}
where $u^{\ast }=u-h\sum_{l=k}^{i-1}m(lh,u).$ Now using a Tailor
expansion we obtain that the right hand side of (\ref{eq:03d}) is
equal to
\begin{eqnarray*} &&
\int \left[ \frac{\partial ^{2}S_{h,r}(0,kh,x,u^{\ast })}{\partial
v_{j}\partial v_{l}}-\sqrt{h}\sum_{\left| \nu \right|
=1}\frac{w^{\nu }}{\nu !}\int_{0}^{1}D_{v}^{\nu }\frac{\partial
^{2}S_{h,r}(0,kh,x,u^{\ast }-\delta \sqrt{h}w)}{\partial
v_{j}\partial v_{l}}d\delta \right]
\\ &&
\times \left[ -\sqrt{h}\sum_{\left| \nu \right| =1}\frac{[w+\sqrt{h}%
\sum_{l=k}^{i-1}m(lh,u)]^{\nu }}{\nu !}D_{u}^{\nu }\sigma
_{jl}(kh,u)\right.
\\ && \qquad +2h\sum_{\left| \nu \right| =2}\frac{[w+\sqrt{h}%
\sum_{l=k}^{i-1}m(lh,u)]^{\nu }}{\nu !}
\\ &&
\left. \qquad \times \int_{0}^{1}D_{u}^{\nu }\sigma _{jl}(kh,u-\delta \sqrt{h}%
w-\delta h\sum_{l=k}^{i-1}m(lh,u))d\delta \right] q^{(i-k)}(kh,u,w)dw.
\end{eqnarray*}
Note that
\begin{eqnarray}\nonumber &&
-\sqrt{h}\int \frac{\partial ^{2}S_{h,r}(0,kh,x,u^{\ast
})}{\partial
v_{j}\partial v_{l}}(w_{p}+\sqrt{h}%
\sum_{l=k}^{i-1}m_{p}(lh,u))q^{(i-k)}(kh,u,w)dw
\\ \nonumber &&
=-h\frac{\partial ^{2}S_{h,r}(0,kh,x,u^{\ast })}{\partial v_{j}\partial v_{l}%
}\sum_{l=k}^{i-1}m_{p}(lh,u),
\\  \label{eq:03e} &&
 h\int \frac{\partial ^{2}S_{h,r}(0,kh,x,u^{\ast
})}{\partial v_{j}\partial
v_{l}}(w_{p}+\sqrt{h}\sum_{l=k}^{i-1}m_{p}(lh,u))(w_{q}+\sqrt{h}%
\sum_{l=k}^{i-1}m_{q}(lh,u))
\\ \nonumber && \qquad
\times \left\{ \int_{0}^{1}D_{u}^{\nu }\sigma _{jl}(kh,u)d\delta
+\int_{0}^{1}\left[ D_{u}^{\nu }\sigma _{jl}(kh,u-\delta \sqrt{h}w
\right . \right . \\ \nonumber && \qquad \qquad \left. \left .
-\delta h\sum_{l=k}^{i-1}m(lh,u))-D_{u}^{\nu }\sigma
_{jl}(kh,u)\right] d\delta \right\} q^{(i-k)}(kh,u,w)dw
\\ \nonumber &&
=h\frac{\partial ^{2}S_{h,r}(0,kh,x,u^{\ast })}{%
\partial v_{j}\partial v_{l}}D_{u}^{\nu }\sigma _{jl}(kh,u)\int
w_{p}w_{q}q^{(i-k)}(kh,u,w)dw
\\ \nonumber && \qquad
+h^{2}\frac{\partial ^{2}S_{h,r}(0,kh,x,u^{\ast })}{\partial
v_{j}\partial
v_{l}}\sum_{l=k}^{i-1}m_{p}(lh,u)\sum_{l=k}^{i-1}m_{q}(lh,u)+R,
\end{eqnarray}
where by (A3$^{\prime }$) we have for $%
j_{0}<(i-k)<n^{\delta ^{\prime }},w^{\prime }=(i-k)^{-1/2}w$ \ \textbf{\ }
\begin{eqnarray} \nonumber
\left| R\right| &\leq& Ch^{3/2}\left| \frac{\partial
^{2}S_{h,r}(0,kh,x,u^{\ast })}{\partial v_{j}\partial
v_{l}}\right| \int \left( n^{\delta ^{\prime }/2}\left\| w^{\prime
}\right\| +O(T^{1/2}n^{-1/2+\delta ^{\prime }})\right) ^{3}\psi
(w^{\prime })dw^{\prime }
\\ \label{eq:03f}
&\leq& Ch\zeta _{\sqrt{ih}}(u-x)(h^{S}n^{S\delta ^{\prime
}}+1)T^{-3/2}n^{-1/2+3\delta ^{\prime }/2}\int \left\| w^{\prime
}\right\| ^{3}\psi (w^{\prime })dw^{\prime } \\ \nonumber &\leq&
Chn^{-1/2+\varkappa /2+3\delta ^{\prime }/2}\zeta
_{\sqrt{ih}}(u-x)\leq Chn^{-1/2(1-3\delta ^{\prime }-\varkappa
)}\zeta _{\sqrt{ih}}(u-x),
\end{eqnarray}
We obtain analogously
\begin{eqnarray*}
&& h\int
w_{p}(w_{q}+\sqrt{h}\sum_{l=k}^{i-1}m_{q}(lh,u))D_{u}^{e_{q}}\sigma
_{jl}(kh,u)\int_{0}^{1}\frac{\partial ^{3}S_{h,r}(0,kh,x,u^{\ast
}-\delta \sqrt{h}w)}{\partial v_{p}\partial v_{j}\partial
v_{l}}d\delta
\\ && \qquad
\times q^{(i-k)}(kh,u,w)dw
\\ &&
=h\frac{\partial \sigma _{jl}(kh,u)}{\partial u_{q}%
}\frac{\partial ^{3}S_{h,r}(0,kh,x,u^{\ast })}{\partial
v_{p}\partial v_{j}\partial v_{l}}\int
w_{p}w_{q}q^{(i-k)}(kh,u,w)dw+R,
\end{eqnarray*}
where
\[
\left| R\right| \leq Chn^{-1/2(1-3\delta ^{\prime }-3\varkappa )}\zeta _{%
\sqrt{ih}}(u-x),1-3\delta ^{\prime }-3\varkappa >0
\]
and, for $1-3\delta ^{\prime }-2\varkappa >0$
\begin{eqnarray} \label{eq:03g} &&
\left| h^{3/2}\int w_{p}\int_{0}^{1}\frac{\partial
^{3}S_{h,r}(0,kh,x,u^{\ast }-\delta \sqrt{h}w)}{\partial
v_{p}\partial v_{j}\partial v_{l}}d\delta \right .\\ \nonumber
&&\left. \qquad \int_{0}^{1}\frac{\partial ^{2}\sigma
_{jl}(kh,u-\delta \sqrt{h}w-\delta
h\sum_{l=k}^{i-1}m_{q}(lh,u))}{\partial u_{r}\partial
u_{s}}d\delta \right.
\\ \nonumber && \qquad \qquad
\left. \times (w_{r}+\sqrt{h}\sum_{l=k}^{i-1}m_{r}(lh,u))(w_{s}+\sqrt{h}%
\sum_{l=k}^{i-1}m_{s}(lh,u))q^{(i-k)}(kh,u,w)dw\right|
\\ \nonumber &&
\leq Chn^{-1/2(1-3\delta ^{\prime }-2\varkappa )}\zeta _{\sqrt{ih}%
}(u-x).
\end{eqnarray}
For $1\leq i-k\leq j_{0}$ \ the same estimates remain true \
because the following bound holds \
\begin{equation}
\int \left\| w\right\| ^{S}q^{(j)}(t,x,w)dw\leq C(j_{0}).  \label{eq:03h}
\end{equation}
The same estimates hold for $\widetilde{p}(kh,ih,v,u)$ \ with \ $\phi
^{(i-k)}(kh,u,w)$ \ instead of $q^{(i-k)}(kh,u,w),$ where $\phi (kh,u,w)$ \
is a gaussian density with the mean \ $0$ \ and \ with \ the covariance
matrix equal to \ $\sigma (kh,u).$ The first two moments of $q^{(i-k)}$ \
and $\phi ^{(i-k)}$ \ coinside so after substraction we obtain uniformly for
\ $i-n^{\delta ^{\prime }}<k\leq i-1$
\begin{eqnarray} \label{eq:03i} &&
\left| \sum_{n/2<i\leq n}h^{2}\sum_{i-n^{\delta ^{\prime }}<k\leq
i-1}\int
\int (L^{T}-\widetilde{L}^{T})S_{h,r}(0,kh,x,v)\right. \\ \nonumber &&
\qquad \left. \times(\widetilde{p}_{h}(kh,ih,v,u)-%
\widetilde{p}(kh,ih,v,u))dv\right.
\\ \nonumber && \qquad
\left. \times
(H+M_{h}+\sqrt{h}N_{1}+hN_{2})^{(n-r-1)}(ih,T,u,y)du\right|
\\ \nonumber &&
\leq \left[ \frac{C^{n-r}}{\Gamma (\frac{n-r-1}{2})}\right]
hT^{3/2}n^{-3/2(1-\varkappa -5\delta ^{\prime }/3)}\zeta
_{\sqrt{ih}}(u-x).
\end{eqnarray}
To estimate the other terms in $\ d_{h}(kh,ih,v,u)$ we need bounds
for the following expressions
\begin{eqnarray*} &&
h\sum_{i-n^{\delta ^{\prime }}<k\leq i-1}\int (L^{T}-\widetilde{L}%
^{T})S_{h,r}(0,kh,x,v)\sqrt{h}(ih-kh)\\
&& \qquad \times D_{v}^{\nu }\widetilde{p}(kh,ih,v,u)dv \mbox{ for
} \left| \nu \right| =3,
\\
&&
h\sum_{i-n^{\delta ^{\prime }}<k\leq i-1}\int (L^{T}-\widetilde{L}%
^{T})S_{h,r}(0,kh,x,v)h(ih-kh)\\
&& \qquad \times D_{v}^{\nu }\widetilde{p}(kh,ih,v,u)dv \mbox{ for
} \left| \nu \right| =4,
\\
&&
h\sum_{i-n^{\delta ^{\prime }}<k\leq i-1}\int (L^{T}-\widetilde{L}%
^{T})S_{h,r}(0,kh,x,v)h(ih-kh)^{2}\\
&& \qquad \times D_{v}^{\nu }\widetilde{p}(kh,ih,v,u)dv \mbox{ for
} \left| \nu \right| =6.
\end{eqnarray*}
We have
\begin{eqnarray} \label{eq:03j} &&
\left| h\sum_{i-n^{\delta ^{\prime }}<k\leq i-1}\int (L^{T}-\widetilde{L}%
^{T})S_{h,r}(0,kh,x,v)\sqrt{h}(ih-kh)\right.
\\ \nonumber && \qquad \left. D_{v}^{\nu }\widetilde{p}%
(kh,ih,v,u)dv\right|
\\ \nonumber &&
=\left| h\sum_{i-n^{\delta ^{\prime }}<k\leq i-1}\int D^{e_{p}+e_{q}}(L^{T}-%
\widetilde{L}^{T})S_{h,r}(0,kh,x,v)\sqrt{h}(ih-kh)\right.
\\ \nonumber && \qquad \left. D_{v}^{\nu -e_{p}-e_{q}}%
\widetilde{p}(kh,ih,v,u)dv\right| \\\nonumber  && \leq
CT^{-2}n^{\delta ^{\prime }}h^{3/2}\sum_{i-n^{\delta ^{\prime
}}<k\leq i-1}\frac{h}{\sqrt{ih-kh}}\zeta _{\sqrt{ih}}(u-x)
\\ \nonumber &&
\leq Chn^{-(1-\varkappa -3\delta ^{\prime }/2)}\zeta
_{\sqrt{ih}}(u-x). \nonumber
\end{eqnarray}
Clearly, the same estimate (\ref{eq:03j}) holds for \ $\left| \nu
\right| =4$ \ and $\left| \nu \right| =6.$ \ Now (\ref{eq:3f}) for
$l=4$ \ follows from this remark and (\ref{eq:03i}) and
(\ref{eq:03j}).

\textit{Proof \ of \ (\ref{eq:3f}) for }$l=1.$ Note that for this case $%
T-ih\geq T/2.$%
\begin{eqnarray}  \label{eq:03ja} &&
a_{r,1}=\sum_{0\leq i\leq n/2}h^{2}\sum_{0\leq k\leq i-1}\int \int (L^{T}-%
\widetilde{L}^{T})S_{h,r}(0,kh,x,v)\\ \nonumber && \qquad \times
d_{h}(kh,ih,v,u)\Psi _{h,r}(ih,T,u,y)dvdu
\\ \nonumber &&
=\sum_{0\leq k\leq n/2-1}h\int (L^{T}-\widetilde{L}^{T})S_{h,r}(0,kh,x,v)%
\\ \nonumber && \qquad \times \left\{ \sum_{k+1\leq i\leq k+n^{\delta ^{\prime }}}h\int
d_{h}(kh,ih,v,u)\Psi _{h,r}(ih,T,u,y)du\right.
\\ \nonumber && \qquad
\left. +\sum_{k+n^{\delta ^{\prime }}<i\leq n/2}h\int
d_{h}(kh,ih,v,u)\Psi _{h,r}(ih,T,u,y)du\right\} dv,
\end{eqnarray}
where we denote
\[
\Psi _{h,r}(ih,T,u,y)=(H+M_{h}+\sqrt{h}N_{1}+hN_{2})^{(n-r-1)}(ih,T,u,y).
\]
We consider
\begin{eqnarray*} &&
\sum_{k+1\leq i\leq k+n^{\delta ^{\prime }}}h\int
h^{-d/2}q^{(i-k)}(kh,u,h^{-1/2}[u-v-h\sum_{l=k}^{i-1}m(lh,u)])\Psi
_{h,r}(ih,T,u,y)du
\\&&
=\sum_{k+1\leq i\leq k+n^{\delta ^{\prime }}}h\int \left\{ q^{(i-k)}(kh,v,w)+%
\sqrt{h}\sum_{\left| \nu \right| =1}(w+\sqrt{h}\sum_{l=k}^{i-1}m(lh,u))^{\nu
}D_{v}^{\nu }q^{(i-k)}(kh,v,w)\right.
\\&& \qquad
+h\sum_{\left| \nu \right|
=2}\frac{(w+\sqrt{h}\sum_{l=k}^{i-1}m(lh,u))^{\nu
}}{\nu !}D_{v}^{\nu }q^{(i-k)}(kh,v,w)
\\&& \qquad
+3h^{3/2}\sum_{\left| \nu \right| =3}%
\frac{(w+\sqrt{h}\sum_{l=k}^{i-1}m(lh,u))^{\nu }}{\nu !}
\\&& \qquad
\left. \times \int_{0}^{1}(1-\delta )^{2}D_{v}^{\nu
}q^{(i-k)}(kh,v+\delta h^{1/2}w+\delta
h\sum_{l=k}^{i-1}m(lh,u),w)d\delta \right\}
\\&& \qquad
\times \left\{ \Psi _{h,r}(ih,T,v,y)+\sqrt{h}\sum_{\left| \nu \right| =1}(w+%
\sqrt{h}\sum_{l=k}^{i-1}m(lh,u))^{\nu }D_{v}^{\nu }\Psi
_{h,r}(ih,T,v,y)\right.
\\&& \qquad
+h\sum_{\left| \nu \right|
=2}\frac{(w+\sqrt{h}\sum_{l=k}^{i-1}m(lh,u))^{\nu }}{\nu
!}D_{v}^{\nu }\Psi _{h,r}(ih,T,v,y)
\\&& \qquad +3h^{3/2}\sum_{\left| \nu
\right| =3}\frac{(w+\sqrt{h}\sum_{l=k}^{i-1}m(lh,u))^{\nu }}{\nu
!}
\\&& \qquad
\left. \times \int_{0}^{1}(1-\delta )^{2}D_{v}^{\nu }\Psi
_{h,r}(ih,T,v+\delta h^{1/2}w+\delta
h\sum_{l=k}^{i-1}m(lh,u),y)d\delta \right\} dw
\end{eqnarray*}
This integral is a sum of $4\times 4=16$ integrals. We estimate
only two of \ them. Other integrals can be estimated by similar
methods. First, we estimate
\[
\sum_{k+1\leq i\leq k+n^{\delta ^{\prime }}}h\int q^{(i-k)}(kh,v,w)\Psi
_{h,r}(ih,T,v,y)dw=\sum_{k+1\leq i\leq k+n^{\delta ^{\prime }}}h\Psi
_{h,r}(ih,T,v,y)dw.
\]
Note that we get the same term when we replace $q^{(i-k)}(kh,v,w)$ by $%
\phi ^{(i-k)}(kh,v,w)$. After the replacement this term
disappears. Second, we estimate
\begin{eqnarray*} &&
\sum_{k+1\leq i\leq k+n^{\delta ^{\prime }}}h\int q^{(i-k)}(kh,v,w)\sqrt{h}%
\sum_{\left| \nu \right| =1}(w+\sqrt{h}\sum_{l=k}^{i-1}m(lh,u))^{\nu
}D_{v}^{\nu }\Psi _{h,r}(ih,T,v,y)dw
\\ &&
=h^{3/2}\sum_{j=1}^{d}\sum_{k+1\leq i\leq k+n^{\delta ^{\prime
}}}D_{v}^{e_{j}}\Psi _{h,r}(ih,T,v,y)\int q^{(i-k)}(kh,v,w)[w_{j}+\sqrt{h}%
\sum_{l=k}^{i-1}m_{j}(lh,v)
\\ && \qquad
+O(hn^{\delta ^{\prime }}\left\| w\right\| +h^{3/2}n^{2\delta
^{\prime }})]dw \\ &&
=h^{2}\sum_{j=1}^{d}\sum_{k+1\leq i\leq
k+n^{\delta ^{\prime }}}D_{v}^{e_{j}}\Psi
_{h,r}(ih,T,v,y)\sum_{l=k}^{i-1}m_{j}(lh,v)
\\ && \qquad
+O\left( h^{2}n^{2\delta ^{\prime }}\sum_{j=1}^{d}\sum_{k+1\leq
i\leq k+n^{\delta ^{\prime }}}h\left| D_{v}^{e_{j}}\Psi
_{h,r}(ih,T,v,y)\right| \right) \\ && \qquad +O\left(
h^{3/2}n^{\delta ^{\prime }}\sum_{j=1}^{d}\sum_{k+1\leq i\leq
k+n^{\delta ^{\prime }}}h \left| D_{v}^{e_{j}}\Psi
_{h,r}(ih,T,v,y)\right| \int q^{(i-k)}(kh,v,w)\left\| w\right\|
dw\right)
\\ &&
=h^{2}\sum_{j=1}^{d}\sum_{k+1\leq i\leq k+n^{\delta ^{\prime
}}}D_{v}^{e_{j}}\Psi
_{h,r}(ih,T,v,y)\sum_{l=k}^{i-1}m_{j}(lh,v)+R,
\end{eqnarray*}
where
\[
\left| R\right| \leq \frac{C^{n-r}}{\Gamma (\frac{n-r-1}{2})}%
T^{1/2}hn^{-3/2+2\delta ^{\prime }}\zeta _{\sqrt{T-kh}}(y-v).
\]
The first term in the right hand side of this equation will be the
same if we replace $q^{(i-k)}(kh,v,w)$ \ \ by $\phi
^{(i-k)}(kh,v,w)$. After the replacement this term disappears. For
a proof of this equation we consider the function $u(w)$ that is
defined as an implicit function and we used the following change
of variables
\[
h^{1/2}w=u-v-h\sum_{l=k}^{i-1}m(lh,u)
\]
to obtain
\[
\sqrt{h}\sum_{l=k}^{i-1}m(lh,u(w))=\sqrt{h}\sum_{l=k}^{i-1}m(lh,v)+O\left(
h(i-k)\left\| w\right\| +h^{3/2}(i-k)^{2}\right)
\]
because of $\ (i-k)\leq n^{\delta ^{\prime }}$. By similar methods
we get
\begin{eqnarray} \label{eq:03m} &&
\left| \sum_{k+1\leq i\leq k+n^{\delta ^{\prime }}}h\int [\sqrt{h}\widetilde{%
\pi }_{1}(kh,ih,v,u)+h\widetilde{\pi }_{2}(kh,ih,v,u)]\Psi
_{h,r}(ih,T,u,y)du\right|
\\ \nonumber && \qquad
\leq \frac{C^{n-r}}{\Gamma (\frac{n-r-1}{2})}hn^{-3/2+2\delta
^{\prime }+\varkappa /2}\zeta _{\sqrt{T-kh}}(y-v).
\end{eqnarray}
It remains to estimate
\[
\sum_{k+n^{\delta ^{\prime }}<i\leq n/2}h\int d_{h}(kh,ih,v,u)\Psi
_{h,r}(ih,T,u,y)du.
\]
From (\ref{eq:03b}) and (\ref{eq:03ba}) we obtain
\begin{eqnarray} \label{eq:03p}
&& \left| \sum_{k+n^{\delta ^{\prime }}<i\leq n/2}h\int
d_{h}(kh,ih,v,u)\Psi
_{h,r}(ih,T,u,y)du\right| \\ \nonumber && \qquad \leq \left[ \frac{C^{n-r}}{\Gamma (\frac{n-r-1}{2})%
}\right] T^{-1/2}h^{3/2}
\int_{kh+n^{\delta ^{\prime }}h}^{T/2}\frac{du}{(u-kh)^{3/2}}\zeta _{%
\sqrt{T-kh}}(y-v)
\\ \nonumber && \qquad
\leq \left[ \frac{C^{n-r}}{\Gamma (\frac{n-r-1}{2})}\right]
T^{-1/2}hn^{-\delta ^{\prime }/2}\zeta _{\sqrt{T-kh}}(y-v)
\\ \nonumber && \qquad
\leq \left[ \frac{C^{n-r}}{\Gamma (\frac{n-r-1}{2})}\right]
hn^{-1/2(\delta ^{\prime }-\varkappa )}\zeta _{\sqrt{T-kh}}(y-v).
\end{eqnarray}
Now we substitute the estimate (\ref{eq:03p}) into
(\ref{eq:03ja}). This
gives the following estimate for any $0<\varepsilon <\varkappa $%
\begin{eqnarray} \label{eq:03q}&&
\left| \sum_{k=1}^{n/2-1}h\int
(L^{T}-\widetilde{L}^{T})S_{h,r}(0,kh,x,v) \right .
\\ \nonumber && \qquad \qquad\left .
\sum_{k+n^{\delta ^{\prime }}<i\leq n/2}h\int d_{h}(kh,ih,v,u)\Psi
_{h,r}(ih,T,u,y)du\right|
\\ \nonumber && \qquad
\leq \left[ \frac{C^{n-r}}{\Gamma (\frac{n-r-1}{2})}\right]
hn^{-1/2(\delta ^{\prime }-\varkappa )}h^{-\varepsilon
}\sum_{k=1}^{n/2}h(kh)^{\varepsilon -1}\zeta _{\sqrt{T}}(y-x)
\\ \nonumber && \qquad
\leq C(\varepsilon )\left[ \frac{C^{n-r}}{\Gamma
(\frac{n-r-1}{2})}\right] hn^{-1/2(\delta ^{\prime }-\varkappa
-\varepsilon )}\zeta _{\sqrt{T}}(y-x).
\end{eqnarray}
For $k=0$ we get with $S_{h,r}(0,0,x,v)=\delta (x-v)$ where
$\delta (\bullet )$ is the Dirac function that
\begin{eqnarray*} &&
\left| \sum_{1\leq i\leq i/2}h^{2}\int \int S_{h,r}(0,0,x,v)(L-\widetilde{L}%
)d_{h}(0,ih,v,u)\Psi _{h,r}(ih,T,u,y)du\right|
\\ && \qquad
\leq C(\varepsilon )\left[ \frac{C^{n-r}}{\Gamma
(\frac{n-r-1}{2})}\right] hn^{-(1/2-\varepsilon )}\zeta
_{\sqrt{T}}(y-x).
\end{eqnarray*}
This completes the proof (\ref{eq:3f}) for $l=1$. The estimate
(\ref {eq:3c}) may be proved by the same arguments as were used to
prove (\ref {eq:T301}).

\textit{Asymptotic treatment of the term $T_{5}.$ \ }We will show that,
\begin{eqnarray} \label{eq:T500} &&
T_{5}=-\sqrt{h}\sum_{r=0}^{\infty }\widetilde{\pi }_{1}\otimes
_{h}(H+M_{h,1}+\sqrt{h}N_{1})^{(r)}(0,T,x,y)
\\ && \qquad \nonumber
-h\sum_{r=0}^{\infty }\widetilde{\pi }_{2}\otimes
_{h}H^{(r)}(0,T,x,y)+R_{h}(x,y),
\end{eqnarray}
where $\left| R_{h}(x,y)\right| \leq Chn^{-\gamma }\zeta _{\sqrt{T}%
}^{S-2}(y-x)$ for some $\gamma >0$. Note that with $
S_{h}(s,t,x,y)=\sum_{r=1}^{n}(K_{h}+M_{h})^{(r)}(s,t,x,y)$ the
term $T_{5}$ can be rewritten as
\[
T_{5}=(\widetilde{p}-\widetilde{p}_{h})(0,T,x,y)+(\widetilde{p}-\widetilde{p}%
_{h})\otimes _{h}S_{h}(0,T,x,y).
\]
We start by showing that for $\ \varkappa <\delta <\frac{1-\varkappa }{4}$ \
\ uniformly for $x,y\in R$%
\begin{eqnarray}  \label{eq:29a} &&
\left| h\sum_{1\leq j\leq n^{\delta }}\int (\widetilde{p}_{h}-\widetilde{p}%
)(0,jh,x,u)S_{h}(jh,T,u,y)du\right|\\ && \qquad \nonumber \leq
O(hn^{-1/2(1-\varkappa -4\delta )})\zeta _{\sqrt{T}}^{S-2}(y-x)
\end{eqnarray}
for $\delta $ small enough. For the proof of (\ref{eq:29a}) we
will show
that uniformly for $1\leq j\leq n^{\delta }$ \ and for $x,y\in R^{d}$%
\begin{eqnarray} \label{eq:29b}
&& \int
\widetilde{p}_{h}(0,jh,x,u)S_{h}(jh,T,u,y)du=S_{h}(jh,T,x,y)
\\ && \qquad \nonumber
+O[h^{1/2}T^{-1/2}n^{-1/2+\delta }+h^{1/2}T^{-1}+n^{\delta
/2}h^{1/2}]\zeta _{\sqrt{T}}^{S-2}(y-x),
\\ &&  \label{eq:29c}
\int \widetilde{p}(0,jh,x,u)S_{h}(jh,T,u,y)du=S_{h}(jh,T,x,y)
\\ && \qquad \nonumber
+O[h^{1/2}T^{-1/2}n^{-1/2+\delta }+h^{1/2}T^{-1}+n^{\delta
/2}h^{1/2}]\zeta _{\sqrt{T}}^{S-2}(y-x).
\end{eqnarray}
Claim (\ref{eq:29a}) immediately follows from (\ref{eq:29b})-(\ref{eq:29c}%
). For the proof we will make use of the fact that for all $1\leq
j\leq
n^{\delta }$ and for all $x,y\in R^{d}$ and $\left| \nu \right| =1$%
\begin{equation}
\left| D_{x}^{\nu }S_{h}(jh,T,x,y)\right| \leq C(T-jh)^{-1}\zeta _{\sqrt{T-jh%
}}^{S-2}(y-x).  \label{eq:29d}
\end{equation}
Claim (\ref{eq:29d}) can be shown with the same arguments as the
proof of (5.7) in Konakov and Mammen (2002). Note that the
function $\Phi $ in that paper has a similar structure as $S_{h}.$
For $1\leq j\leq n^{\delta }$
the bound (\ref{eq:29d}) immediately implies for a constant $C^{\prime }$%
\begin{equation}
\left| D_{x}^{\nu }S_{h}(jh,T,x,y)\right| \leq C^{\prime }T^{-1}\zeta _{%
\sqrt{T}}^{S-2}(y-x).  \label{eq:29e}
\end{equation}
We have $\widetilde{p}_{h}(0,jh,x,u)=h^{-d/2}q^{(j)}[0,u,h^{-1/2}(u-x-h%
\sum_{i=0}^{j-1}m(ih,u))].$ Denote the determinant of the Jacobian
matrix of \ $u-h\sum_{i=0}^{j-1}m(ih,u)$ \ by $\Delta _{h}$.
Because of condition (A3) and (\ref{eq:29e}) it holds that for
$1\leq j\leq n^{\delta }$
\begin{eqnarray*}
&& \int \widetilde{p}_{h}(0,jh,x,u)S_{h}(jh,T,u,y)du
\\&&
=\int
h^{-d/2}q^{(j)}[0,u,h^{-1/2}(u-x-h\sum_{i=0}^{j-1}m(ih,u))]S_{h}(jh,T,u,y)du
\\&&
=\int q^{(j)}(0,x+h^{1/2}w+h\sum_{i=0}^{j-1}m(ih,u(w)),w)\left|
\Delta _{h}^{-1}\right|
\\ && \qquad
S_{h}(jh,T,x+h^{1/2}w+h\sum_{i=0}^{j-1}m(ih,u(w)),y)dw
\\ &&
=\int [q^{(j)}(0,x,w)+O(j^{-d/2}h^{1/2})(\left\| w\right\| +1)\psi
(j^{-1/2}w)][1+O(jh)][S_{h}(jh,T,x,y)
\\ && \qquad
+O(h^{1/2}T^{-1})\zeta _{\sqrt{T}}^{S-2}(y-x)(1+h^{(S-2)/2}\left\| w\right\|
^{S-2})(\left\| w\right\| +1)]dw
\\ &&
=S_{h}(jh,T,x,y)+O[h^{1/2}T^{-1/2}n^{-1/2+\delta
}+h^{1/2}T^{-1}+h^{1/2}n^{\delta /2}]\zeta _{\sqrt{T}}^{S-2}(y-x)
\end{eqnarray*}
with $u=u(w)$ in $\sum_{i=0}^{j-1}m(ih,u)$ defined by the Inverse
Function Theorem with the substitution $w=$ $h^{-1/2}(u-x-h%
\sum_{i=0}^{j-1}m(ih,u))$. This proves (\ref{eq:29b}). The proof
of (\ref {eq:29c}) is the same with obvious modifications. From
(\ref{eq:29a}) we get that for $\delta <\frac{1-\varkappa }{4}$
(where $\varkappa $ is defined in (B2))
\[
T_{5}=(\widetilde{p}-\widetilde{p}_{h})(0,T,x,y)+h\sum_{n^{\delta }<j<n}\int
(\widetilde{p}-\widetilde{p}_{h})(0,jh,x,u)S_{h}(jT,u,y)du+R_{h}(x,y)
\]
with $\left| R_{h}(x,y)\right| \leq O(hn^{-1/2(1-\varkappa -4\delta )})\zeta
_{\sqrt{T}}^{S-2}(y-x)$. We now make use of the expansion of
$%
\widetilde{p}_{h}-\widetilde{p}$ given in Lemma 5. We have with
$\rho =(jh)^{1/2}\geq h^{1/2}n^{\delta /2}$
\begin{eqnarray} \label{eq:29g} &&
\left| h\sum_{j=n^{\delta }}^{n}h^{3/2}\rho ^{-3}\int \zeta _{\rho
}^{S}(u-x)S_{h}(jh,T,u,y)du\right|\\ \nonumber && \qquad \leq
Ch^{2}T^{-\delta ^{\prime }}n^{-\delta ^{\prime \prime
}}\sum_{j=n^{\delta }}^{n}\rho ^{-2+2\delta ^{\prime }}\int \left|
\zeta _{\rho }^{S}(u-x)S_{h}(jh,T,u,y)\right| du, \end{eqnarray}
where $\delta ^{\prime }<\frac{1}{2}\delta (1-\delta )^{-1},2$
$\delta ^{\prime \prime }=\delta +2\delta \delta ^{\prime
}-2\delta ^{\prime }.$ Now we get that
\begin{equation}
h\sum_{j=n^{\delta }}^{n}\rho ^{-2+2\delta ^{\prime }}\int \left| \zeta
_{\rho }^{S}(u-x)S_{h}(jh,T,u,y)\right| du\leq CB(\delta ^{\prime
},1/2)T^{\delta ^{\prime }-1/2}\zeta _{\sqrt{T}}^{S-2}(y-x)  \label{eq:30}
\end{equation}
for a constant $C.$ This shows that for $\delta ^{\prime }>0$
small enough
\[
T_{5}=-[\sqrt{h}\widetilde{\pi }_{1}+h\widetilde{\pi }_{2}](0,T,x,y)
\]
\[
-h\sum_{n^{\delta }<j<n}\int [\sqrt{h}\widetilde{\pi }_{1}+h\widetilde{\pi }%
_{2}](0,jh,x,u)S_{h}(jh,T,u,y)du+R_{h}^{\prime }(x,y),
\]
with $\left| R_{h}^{\prime }(x,y)\right| \leq O(hn^{-(\delta ^{\prime \prime
}-\varkappa /2)})\zeta _{\sqrt{T}}^{S-2}(y-x)$ \ with a constant in $O(\cdot
)$ depending on $\delta ^{\prime }.$ It follows from (\ref{eq:29a}), (\ref
{eq:29g}) and (\ref{eq:30}) that
\begin{equation}
T_{5}=-\sum_{r=0}^{\infty }[\sqrt{h}\widetilde{\pi }_{1}+h\widetilde{\pi }%
_{2}]\otimes _{h}(K_{h}+M_{h})^{(r)}(0,T,x,y)+R_{h}^{\prime \prime }(x,y),
\label{eq:31}
\end{equation}
where $\left| R_{h}^{\prime \prime }(x,y)\right| \leq O(hn^{-(\delta
^{\prime \prime }-\varkappa /2)})\varsigma _{\sqrt{T}}^{S-2}(y-x).$ Now we
apply Lemma 10 \ with $\ A=\sqrt{h}\widetilde{\pi }_{1},$ $B=H+M_{h,1}+\sqrt{%
h}N_{1},C=(K_{h}-H-\sqrt{h}N_{1})+(M_{h}-M_{h,1})$ \ to
\begin{eqnarray} \label{eq:T501}&&
-\sum_{r=0}^{\infty }\sqrt{h}\widetilde{\pi }_{1}\otimes
_{h}(K_{h}+M_{h})^{(r)}(0,T,x,y)\\ && \nonumber \qquad
+\sum_{r=0}^{\infty }\sqrt{h}\widetilde{\pi }%
_{1}\otimes _{h}(H+M_{h,1}+\sqrt{h}N_{1})^{(r)}(0,T,x,y)
\end{eqnarray}
and with $A=h\widetilde{\pi }_{2},B=H,C=(K_{h}-H)+M_{h}$ to
\begin{equation}
-\sum_{r=0}^{\infty }h\widetilde{\pi }_{2}\otimes
_{h}(K_{h}+M_{h})^{(r)}(0,T,x,y)+\sum_{r=0}^{\infty }h\widetilde{\pi }%
_{2}\otimes _{h}H^{(r)}(0,T,x,y).  \label{eq:T502}
\end{equation}
The estimate (\ref{eq:T500}) follows from (\ref{eq:30}), (\ref{eq:T501}), (%
\ref{eq:T502}), Lemma 10, Lemma 5 , (\ref{eq:12e}) and (\ref{eq:12f}).

\textit{Asymptotic treatment of the term $T_{6}$. }By Lemma 9
\[
\left| T_{6}\right| \leq C(\varepsilon )hn^{-1/2+\varepsilon }\zeta _{\sqrt{T%
}}^{S}(y-x).
\]

\textit{Asymptotic treatment of the term $T_{7}$. } We use the
recurrence
relation for $r=2,3,...$%
\[
\widetilde{p}_{h}\otimes _{h}(K_{h}+M_{h}+R_{h})^{(r)}(0,T,x,y)-\widetilde{p}%
_{h}\otimes _{h}H_{h}{}^{(r)}(0,T,x,y)
\]
\[
=\left[ \widetilde{p}_{h}\otimes _{h}(K_{h}+M_{h}+R_{h})^{(r-1)}-\widetilde{p%
}_{h}\otimes _{h}H_{h}{}^{(r-1)}\right] \otimes _{h}H_{h}(0,T,x,y)
\]
\[
+[\widetilde{p}_{h}\otimes _{h}(K_{h}+M_{h}+R_{h})^{(r-1)}\otimes
_{h}(K_{h}+M_{h}+R_{h}-H_{h})](0,T,x,y)
\]
and we apply Lemma 8 for $r=1$. We get in the same way as in the
proof of Lemma 9
\[
\left| T_{7}\right| \leq Ch^{3/2}T^{-1/2}\zeta _{\sqrt{T}%
}^{S}(y-x)=Chn^{-1/2}\zeta _{\sqrt{T}}^{S}(y-x).
\]
\textit{Plugging in the asymptotic expansions of \
}$T_{1},...,T_{7}$ . We now plug the asymptotic expansions of \
$T_{1},...,T_{7}$ \ $\ $into (\ref {eq:26}). Using Lemma 10,
Theorem 2.1 in Konakov and Mammen (2002) we get
\begin{eqnarray} \label{eq:P01} &&
p_{h}(0,T,x,y)-p(0,T,x,y)
\\ && \nonumber
=\sqrt{h}\left[ \widetilde{\pi }_{1}+p^{d}\otimes _{h}\Re _{1}\right]
\otimes _{h}\Phi (0,T,x,y)
\\ && \nonumber \qquad
+h\left\{ \left[ \widetilde{\pi }_{2}+\widetilde{\pi }_{1}\otimes
_{h}\Phi
\otimes _{h}\Re _{1}+p^{d}\otimes _{h}\Re _{2}+p^{d}\otimes _{h}\Re _{3}%
\right] \otimes _{h}\Phi (0,T,x,y)\right.
\\ && \nonumber \qquad
+p^{d}\otimes _{h}\left( \Re _{1}\otimes _{h}\Phi \right)
^{(2)}(0,T,x,y)
\\ && \nonumber \qquad
\left. +\frac{1}{2}p\otimes _{h}(L_{\ast }^{2}-L^{2})p^{d}(0,T,x,y)-\frac{1}{%
2}p\otimes _{h}(L^{\prime }-\widetilde{L}^{\prime })p^{d}(0,T,x,y)\right\}
\\ && \nonumber \qquad +O(h^{1+\delta }\zeta _{\sqrt{T}}(y-x)),
\end{eqnarray}
where
\[
\Re _{1}(s,t,x,y)=N_{1}(s,t,x,y)+M_{1}(s,t,x,y)-\widetilde{M}_{1}(s,t,x,y),
\]
\[
\Re _{2}(s,t,x,y)=N_{2}(s,t,x,y)+\Pi _{1}(s,t,x,y)-\widetilde{\Pi }%
_{1}(s,t,x,y),
\]
\[
\Re _{3}(s,t,x,y)=\sum_{\left| \nu \right| =4}\frac{(\chi _{\nu }(s,x)-\chi
_{\nu }(s,y)}{\nu !}D_{x}^{\nu }\widetilde{p}(s,t,x,y),
\]
\[
M_{1}(s,t,x,y)=\sum_{\left| \nu \right| =3}\frac{\chi _{\nu }(s,x)}{\nu !}%
D_{x}^{\nu }\widetilde{p}(s,t,x,y),\widetilde{M}_{1}(s,t,x,y)=\sum_{\left|
\nu \right| =3}\frac{\chi _{\nu }(s,y)}{\nu !}D_{x}^{\nu }\widetilde{p}%
(s,t,x,y),
\]
\[
\Pi _{1}(s,t,x,y)=\sum_{\left| \nu \right| =3}\frac{\chi _{\nu }(s,x)}{\nu !}%
D_{x}^{\nu }\widetilde{\pi }_{1}(s,t,x,y),\widetilde{\Pi }%
_{1}(s,t,x,y)=\sum_{\left| \nu \right| =3}\frac{\chi _{\nu }(s,y)}{\nu !}%
D_{x}^{\nu }\widetilde{\pi }_{1}(s,t,x,y).
\]
For the homogenous case and $T=[0,1]$ \ (\ref{eq:P01}) coinsides \
with formula (53) \ on page 623 in Konakov and Mammen (2005).

\textit{Asymptotic replacement of }$p^{d}$ by $p$. It follows from
(11), (26) and (27) in Konakov (2006) that
\begin{equation}
\left| (p^{d}-p)(ih,jh,x,z)\right| \leq C(\varepsilon )h^{1-\varepsilon
}(jh-ih)^{\varepsilon -1/2}\phi _{\sqrt{(j-i)h}}(z-x)  \label{eq:P02}
\end{equation}
for any $0<\varepsilon <1/2.$\ Using (\ref{eq:P02}) and making an
integration by parts we can replace $p^{d}$ by $p$ in (\ref{eq:P01}%
). For example the operator $L_{\ast }^{2}-L^{2}$ is the operator
of the third order. Making an integration by parts we have for
$\left| \nu
\right| =3$%
\begin{eqnarray} \label{eq:P03} &&
\left| \sum_{i=1}^{n-1}h\int D_{z}^{\nu
}p(0,ih,x,z)(p^{d}-p)(ih,T,z,y)dz\right|\\ \nonumber && \leq
C(\varepsilon
)h^{1-\varepsilon }\sum_{i=1}^{n-1}h\frac{1}{(ih)^{3/2}}\frac{1}{%
(T-ih)^{1/2-\varepsilon }}\phi _{\sqrt{T}}(y-x)
\\ \nonumber &&
\leq C(\varepsilon )h^{1/2-2\varepsilon }T^{2\varepsilon
-1/2}B(\varepsilon ,\varepsilon +\frac{1}{2})\phi
_{\sqrt{T}}(y-x).
\end{eqnarray}
By (B2) we have $\ 0<\varkappa <1-4\varepsilon .$ This implies
\[
\left| \frac{h}{2}p\otimes _{h}(L_{\ast
}^{2}-L^{2})(p^{d}-p)(0,T,x,y)\right| \leq C(\varepsilon
)hT^{1/2}n^{-(1/2-2\varepsilon -\varkappa /2)}\phi _{\sqrt{T}}(y-x)
\]
\[
\leq C(\varepsilon )h^{1+\delta }\phi _{\sqrt{T}}(y-x)
\]
for some $0<\delta <1/2.$ \ The other terms in (\ref{eq:P01}) containing  $%
p^{d}$ can be estimated analogously. Thus we get the following
representation
\begin{eqnarray}  \label{eq:P04}&&
p_{h}(0,T,x,y)-p(0,T,x,y)
\\ \nonumber &&
=\sqrt{h}\left[ \widetilde{\pi }_{1}+p\otimes _{h}\Re _{1}\right] \otimes
_{h}\Phi (0,T,x,y)
\\ \nonumber && \qquad
+h\left\{ \left[ \widetilde{\pi }_{2}+\widetilde{\pi }_{1}\otimes
_{h}\Phi \otimes _{h}\Re _{1}+p\otimes _{h}\Re _{2}+p\otimes
_{h}\Re _{3}\right] \otimes _{h}\Phi (0,T,x,y)\right.
\\ \nonumber && \qquad
+p\otimes _{h}\left( \Re _{1}\otimes _{h}\Phi \right)
^{(2)}(0,T,x,y)
\\ \nonumber && \qquad
\left. +\frac{1}{2}p\otimes _{h}(L_{\ast }^{2}-L^{2})p(0,T,x,y)-\frac{1}{2}%
p\otimes _{h}(L^{\prime }-\widetilde{L}^{\prime
})p(0,T,x,y)\right\} \\ \nonumber && \qquad +O(h^{1+\delta }\zeta
_{\sqrt{T}}(y-x)).
\end{eqnarray}
For the further analysis we need a generalization of a binary operation $%
\otimes $ introduced in Konakov and Mammen (2005). Recall a
corresponding definition. Suppose that $s\in \lbrack 0,t-h]$ and
$t\in \{h,2h,...,T\}.$ Then the binary type operation $\otimes
_{h}^{\prime }$ is defined as follows
\[
f\otimes _{h}^{\prime }g(s,t,x,y)=\sum_{s\leq jh\leq t-h}h\int
f(s,jh,x,z)g(jh,t,z,y)dz.
\]
Note that for $s\in \{0,h,2h,...,T\}$ these two operations coincide, that is
$\otimes _{h}^{\prime }\equiv \otimes _{h}.$

\textit{Asymptotic replacement of }$(p\otimes _{h}\Re _{i})\otimes
_{h}\Phi (0,T,x,y)$ \ \ \textit{by} \ \ $p\otimes (\Re _{i}\otimes
_{h}^{\prime }\Phi )(0,T,x,y)=(p\otimes \Re _{i})\otimes _{h}\Phi
(0,T,x,y),i=1,2,3,$ \textit{\ }$[p\otimes _{h}\left( \Re
_{1}\otimes _{h}\Phi \right) ]\otimes _{h}\left( \Re _{1}\otimes
_{h}\Phi \right) (0,T,x,y)$ \ \ \textit{by \ }$p\otimes \lbrack
(\Re _{i}\otimes _{h}^{\prime }\Phi )\otimes _{h}^{\prime }(\Re
_{i}\otimes _{h}^{\prime }\Phi )](0,T,x,y), $ $p\otimes
_{h}(L_{\ast }^{2}-L^{2})p(0,T,x,y)$ \ \textit{by \ }\ $p\otimes
(L_{\ast }^{2}-L^{2})p(0,T,x,y)$ \ \textit{and} $\ \ p\otimes
_{h}(L^{\prime }-\widetilde{L}^{\prime })p(0,T,x,y)$ \ \textit{by \ \ } $%
p\otimes (L^{\prime }-\widetilde{L}^{\prime })p(0,T,x,y).$

These replacements follow from the definitions of $\Re
_{i},i=1,2,3,$ and can be proved by the same method as in Konakov
(2006), pp. 9-12, \ where we estimate the replacement error of \
$p\otimes _{h}H$ \ by $p\otimes H.$ Linearity of the operation
$\otimes _{h}$ implies that it is enough to consider \ the
functions $p\otimes _{h}\Im $ \ where $\Im (u,t,z,v)$ \ is a
function of one of the followin forms:
\[
\frac{\chi _{\nu }(u,z)-\chi _{\nu }(u,v)}{\nu !}D_{x}^{\nu }\widetilde{p}%
(u,t,z,v),\left| \nu \right| =3,4,\frac{\chi _{\nu }(u,z)-\chi _{\nu }(u,v)}{%
\nu !}D_{x}^{\nu }\widetilde{\pi }_{1}(u,t,z,v),\left| \nu \right| =3
\]
\[
(L-\widetilde{L})\widetilde{\pi }_{1}(u,t,z,v),(L-\widetilde{L})\widetilde{%
\pi }_{2}(u,t,z,v).
\]
We consider the case $\Im (u,t,z,v)=(L-\widetilde{L})\widetilde{\pi }%
_{1}(u,t,z,v).$ The other cases can be treated similarly. It is enough to
consider a typical term of $(L-\widetilde{L})\widetilde{\pi }_{1}(u,t,z,v),$
namely, we shall estimate
\[
\int_{0}^{jh}du\int p(0,u,x,z)\left( \int_{u}^{jh}\chi _{\nu }(w,v)dw\right)
D_{z}^{\nu }(L-\widetilde{L})\widetilde{p}(u,jh,z,v)dz
\]
\[
-\sum_{i=0}^{j-1}h\int p(0,ih,x,z)\left( \int_{ih}^{jh}\chi _{\nu
}(w,v)dw\right) D_{z}^{\nu }(L-\widetilde{L})\widetilde{p}(ih,jh,z,v)dz
\]
\[
\sum_{i=0}^{j-1}\int_{ih}^{(i+1)h}du\int [\lambda (u)-\lambda
(ih)]dz=\sum_{i=0}^{j-1}\int_{ih}^{(i+1)h}(u-ih)du\int \lambda ^{\prime
}(ih)dz
\]
\begin{equation}
+\frac{1}{2}\sum_{i=0}^{j-1}\int_{ih}^{(i+1)h}(u-ih)^{2}du\int_{0}^{1}(1-%
\delta )\int \lambda ^{\prime \prime }(s)\mid _{s=s_{i}}dzd\delta
du, \label{eq:P04aa}
\end{equation}
where $\lambda (u)=p(0,u,x,z)\left( \int_{u}^{jh}\chi _{\nu }(w,v)dw\right)
D_{z}^{\nu }H(u,jh,z,v),s_{i}=ih+\delta (u-ih).$ As in Konakov (2006), p.9,
we obtain that
\[
\sum_{i=0}^{j-1}\int_{ih}^{(i+1)h}(u-ih)du\int \lambda ^{\prime }(ih)dz
\]
\[
=\frac{h}{2}\sum_{i=0}^{j-1}h\int_{ih}^{jh}\chi _{\nu }(s,v)ds\int
p(0,ih,x,z)D_{z}^{\nu }A_{0}(ih,jh,z,v)dz
\]
\[
+\frac{h}{2}\sum_{i=0}^{j-1}h\int_{ih}^{jh}\chi _{\nu }(s,v)ds\int
p(0,ih,x,z)D_{z}^{\nu }H_{1}(ih,jh,z,v)dz
\]
\begin{equation}
-\frac{h}{2}\sum_{i=0}^{j-1}h\chi _{\nu }(ih,v)\int p(0,ih,x,z)D_{z}^{\nu
}H(ih,jh,z,v)dz,  \label{eq:P04a}
\end{equation}
where
\[
A_{0}(s,jh,z,v)=(L^{2}-2L\widetilde{L}+\widetilde{L}^{2})\widetilde{p}%
(s,jh,z,v),
\]
\[
H_{l}(s,t,z,v)=\frac{1}{2}\sum_{i,j=1}^{d}\left( \frac{\partial ^{l}\sigma
_{ij}(s,z)}{\partial s^{l}}-\frac{\partial ^{l}\sigma _{ij}(s,v)}{\partial
s^{l}}\right) \frac{\partial ^{2}\widetilde{p}(s,t,z,v)}{\partial
z_{i}\partial z_{j}}
\]
\[
\sum_{i=1}^{d}\left( \frac{\partial ^{l}m_{i}(s,z)}{\partial s^{l}}-\frac{%
\partial ^{l}m_{i}(s,v)}{\partial s^{l}}\right) \frac{\partial \widetilde{p}%
(s,t,z,v)}{\partial z_{i}},l=0,1,2,H_{0}\equiv H.
\]
The differential operator \ $A_{0}$ \ was calculated in Konakov
(2006), p.9. $A_{0}$ is a forth order differential operator. From
the structure of this operator and from (\ref{eq:P04a}) it is
clear that it'is enough to estimate
\begin{equation}
I\triangleq \frac{h}{2}\sum_{j=0}^{n-1}h\int
\sum_{i=0}^{j-1}h\int_{ih}^{jh}\chi _{\nu }(s,v)ds\int p(0,ih,x,z)D_{z}^{\nu
+\mu }\widetilde{p}(ih,jh,z,v)dz\Phi (jh,T,v,y)dv  \label{eq:P04b}
\end{equation}
for $\left| \nu \right| =3,\left| \mu \right| =3.$ To estimate (\ref{eq:P04b}%
) we consider three possible cases: a) $jh>T/2,ih\leq jh/2$ $\Longrightarrow
jh-ih>T/4$ \ b) $jh>T/2,ih>jh/2$ $\Longrightarrow $ $ih>T/4$ \ c) $jh<T/2$ $%
\Longrightarrow $ $T-jh>T/2.$ In the case a) we make an integration by parts
transferring two derivatives to $p(0,ih,x,z).$ This gives
\[
\left| \frac{h}{2}\sum_{i=0}^{j-1}h\int_{ih}^{jh}\chi _{\nu }(s,v)ds\int
D_{z}^{e_{k}+e_{l}}p(0,ih,x,z)D_{z}^{\nu +\mu -e_{k}-e_{l}}\widetilde{p}%
(ih,jh,z,v)dz\right|
\]
\[
\leq Ch^{1-2\varepsilon }\int_{0}^{jh}\frac{1}{u^{1-\varepsilon }}\frac{du}{%
(jh-u)^{1-\varepsilon }}\phi _{\sqrt{jh}}(v-x)
\]
\[
\leq C(\varepsilon )h^{1-2\varepsilon }(jh)^{2\varepsilon -1}\phi _{\sqrt{jh}%
}(v-x)
\]
and
\[
\left| I\right| \leq C(\varepsilon )h^{1-2\varepsilon }\int_{0}^{T}\frac{du}{%
u^{1-2\varepsilon }(T-u)^{1/2}}\phi _{\sqrt{T}}(y-x)
\]
\begin{equation}
\leq C(\varepsilon )h^{3/4}T^{1/2}n^{-(1/4-2\varepsilon -\varkappa /2)}\phi
_{\sqrt{T}}(y-x)\leq C(\varepsilon )T^{1/2-\delta }h^{3/4+\delta }\phi _{%
\sqrt{T}}(y-x),  \label{eq:P04c}
\end{equation}
where $\delta =(1/4-2\varepsilon -\varkappa /2)>0$ \ if \ $\varkappa
<1/2-4\varepsilon ,$ $0<\varepsilon <0,05$ \ (see the condition (B2)%
). In the case b) we make an integration by parts transferring \
four derivatives to $p(0,ih,x,z).$ This gives the same estimate
(\ref{eq:P04c}). At last, in the case c) we make an integration by
parts transferring three derivatives to $\Phi (jh,T,v,y)$ and one
derivative to $p(0,ih,x,z).$ This
gives the same estimate (\ref{eq:P04c}). To pass from $D_{z}^{\mu }%
\widetilde{p}(ih,jh,z,v)$ to $\ D_{v}^{\mu }\widetilde{p}(ih,jh,z,v)$ \ we
used the following estimate
\[
\left| D_{z}^{\mu }\widetilde{p}(ih,jh,z,v)+D_{v}^{\mu }\widetilde{p}%
(ih,jh,z,v)\right| \leq C\phi _{\sqrt{jh-ih}}(v-z).
\]
Clearly, the same estimate (\ref{eq:P04c}) hold true for the other summands
in the right hand side of (\ref{eq:P04a})

\[
\frac{h}{2}\left| \sum_{j=0}^{n-1}h\int \sum_{i=0}^{j-1}h\int_{ih}^{jh}\chi
_{\nu }(s,v)ds\int p(0,ih,x,z)D_{z}^{\nu }H_{1}(ih,jh,z,v)dz\Phi
(jh,T,v,y)dv\right|
\]
\[
\leq C(\varepsilon )T^{1/2-\delta }h^{3/4+\delta }\phi _{\sqrt{T}}(y-x),
\]
\[
\frac{h}{2}\sum_{i=0}^{j-1}h\int \chi _{\nu }(ih,v)\sum_{i=0}^{j-1}h\int
p(0,ih,x,z)D_{z}^{\nu }H(ih,jh,z,v)dz\Phi (jh,T,v,y)dv
\]
\[
\leq C(\varepsilon )T^{1/2-\delta }h^{3/4+\delta }\phi _{\sqrt{T}}(y-x).
\]
Now we shall estimate the second summand in the right hand side of
(\ref {eq:P04a}). Analogously to (19) in Konakov (2006) we obtain
\begin{eqnarray} \label{eq:P04d} &&
\frac{1}{2}\sum_{i=0}^{j-1}\int_{ih}^{(i+1)h}(u-ih)^{2}du\int_{0}^{1}(1-%
\delta )\int \lambda ^{\prime \prime }(s)\mid _{s=s_{i}}dzd\delta du
\\ && \nonumber =
\frac{1}{2}\sum_{i=0}^{j-1}\int_{ih}^{(i+1)h}(u-ih)^{2}\int_{0}^{1}(1-\delta
)\sum_{k=1}^{4}\int_{s}^{jh}\chi _{\nu }(\tau ,v)d\tau \\ &&
\nonumber \qquad  \times \int p(0,s,x,z)D_{z}^{\nu
}A_{k}(s,jh,z,v)\mid _{s=s_{i}}dzd\delta du
\\ &&
\nonumber \qquad
-\sum_{i=0}^{j-1}\int_{ih}^{(i+1)h}(u-ih)^{2}\int_{0}^{1}(1-\delta
)\chi _{\nu }(s,v)\\ && \nonumber \qquad  \times  \int
{}p(0,s,x,z)D_{z}^{\nu }A_{0}(s,jh,z,v)\mid _{s=s_{i}}dzd\delta du
\\ &&
\nonumber \qquad
-\sum_{i=0}^{j-1}\int_{ih}^{(i+1)h}(u-ih)^{2}\int_{0}^{1}(1-\delta
)\chi _{\nu }(s,v)\\ && \nonumber \qquad  \times \int
p(0,s,x,z)D_{z}^{\nu }H_{1}(s,jh,z,v)\mid _{s=s_{i}}dzd\delta du
\\ &&
\nonumber \qquad
-\sum_{i=0}^{j-1}\int_{ih}^{(i+1)h}(u-ih)^{2}\int_{0}^{1}(1-\delta )\frac{%
\partial \chi _{\nu }(s,v)}{\partial s}\\ &&
\nonumber \qquad  \times \int p(0,s,x,z)D_{z}^{\nu
}H(s,jh,z,v)\mid _{s=s_{i}}dzd\delta du,
\end{eqnarray}
where the operators $A_{i},i=1,2,3,4,$ are defined as follows:
\[
A_{1}(s,jh,z,v)=(L^{3}-3L^{2}\widetilde{L}+3L\widetilde{L}^{2}-\widetilde{L}%
^{3})\widetilde{p}(s,jh,z,v),
\]
\[
A_{2}(s,jh,z,v)=(L_{1}H+2LH_{1})(s,jh,z,v),
\]
\[
A_{3}(s,jh,z,v)=[(L-\widetilde{L})\widetilde{L}_{1}+2(L_{1}-\widetilde{L}%
_{1})\widetilde{L}]\widetilde{p}(s,jh,z,v),
\]
\begin{equation}
A_{4}(s,jh,z,v)=H_{2}(s,jh,z,v).  \label{eq:P04e}
\end{equation}
The operator $A_{1}$ is given in (24) in Konakov (2006). As in Konakov
(2006) it is enough to estimate for fixed $\ p,q,r,l$%
\begin{eqnarray*} &&
\frac{1}{2}\sum_{i=0}^{j-1}\int_{ih}^{(i+1)h}(u-ih)^{2}\int_{0}^{1}(1-\delta
)\sum_{k=1}^{4}\int_{s}^{jh}\chi _{\nu }(\tau ,v)d\tau \\ &&
\qquad \int
p(0,s,x,z)D_{z}^{\nu }\left( \frac{\partial ^{4}\widetilde{p}(s,jh,z,v)}{%
\partial z_{p}\partial z_{q}\partial z_{l}\partial z_{r}}\right) \mid
_{s=s_{i}}dzd\delta du.
\end{eqnarray*}
As in (25) in Konakov (2006) we obtain that (\ref{eq:P04f}) \ does not
exceed
\begin{equation}
C(\varepsilon )h^{3/2-\varepsilon }(jh)^{2\varepsilon -1}\phi _{\sqrt{jh}%
}(v-x).  \label{eq:P04f}
\end{equation}
It follows from the explicit form of these operators (\ref{eq:P04e}) \ that
the same estimate (\ref{eq:P04f}) \ holds for $A_{2},A_{3}$ \ and $A_{4}.$
The other three terms in the right hand side of (\ref{eq:P04d}) \ do not
contain the factor $\int_{s}^{jh}\chi _{\nu }(\tau ,v)d\tau $ \ and they
should be estimated separately. Clearly, it's enough to estimate the term
containing $A_{0}.$ The remaining two summands are less singular. From the
explicit form of $\ A_{0}$ (formula (15) in Konakov (2006)) we obtain that
it is enough to estimate for fixed $q,l,r$
\begin{eqnarray} && \label{eq:P04g}
\sum_{i=0}^{j-1}\int_{ih}^{(i+1)h}(u-ih)^{2}\int_{0}^{1}(1-\delta
)\chi
_{\nu }(s,v)\\ \nonumber && \qquad \int {}p(0,s,x,z)D_{z}^{\nu }\left( \frac{\partial ^{3}%
\widetilde{p}(s,jh,z,v)}{\partial z_{q}\partial z_{l}\partial
z_{r}}\right) (s,jh,z,v)\mid _{s=s_{i}}dzd\delta du.
\end{eqnarray}
Analogously to (25) in Konakov (2006) we get that (\ref{eq:P04g})
does not exceed
\begin{equation}
C(\varepsilon )h^{1-2\varepsilon }(jh)^{2\varepsilon -1}\phi _{\sqrt{jh}%
}(v-x).  \label{eq:P04h}
\end{equation}
Now from (\ref{eq:P04aa}), (\ref{eq:P04c}), (\ref{eq:P04d}), (\ref{eq:P04f})
and (\ref{eq:P04h}) we obtain that
\[
\left| \lbrack p\otimes _{h}(L-\widetilde{L})\widetilde{\pi }_{1}]\otimes
_{h}\Phi (0,T,x,y)-p\otimes \lbrack (L-\widetilde{L})\widetilde{\pi }%
_{1}\otimes _{h}^{\prime }\Phi ](0,T,x,y)\right|
\]
\begin{equation}
\leq Ch^{3/4+\delta }\phi _{\sqrt{T}}(y-x)  \label{eq:P04i}
\end{equation}
for some \ $\delta >0.$ The other replacements can be shown analogously.
Thus we come to the following representation
\[
p_{h}(0,T,x,y)-p(0,T,x,y)
\]
\[
=\sqrt{h}\left[ \widetilde{\pi }_{1}\otimes _{h}^{\prime }\Phi
(0,T,x,y)+p\otimes (\Re _{1}\otimes _{h}^{\prime }\Phi )(0,T,x,y)\right]
\]
\[
+h\left[ \widetilde{\pi }_{2}\otimes _{h}^{\prime }\Phi (0,T,x,y)+p\otimes
(\Re _{2}\otimes _{h}^{\prime }\Phi )(0,T,x,y)+p\otimes _{h}(\Re _{3}\otimes
_{h}^{\prime }\Phi )(0,T,x,y)\right]
\]
\[
+h\left[ \widetilde{\pi }_{1}\otimes _{h}^{\prime }\Phi +p\otimes \left( \Re
_{1}\otimes _{h}^{\prime }\Phi \right) \right] \otimes _{h}^{\prime }(\Re
_{1}\otimes _{h}^{\prime }\Phi )(0,T,x,y)
\]
\begin{equation}
+\frac{h}{2}p\otimes (L_{\ast }^{2}-L^{2})p(0,T,x,y)-\frac{h}{2}p\otimes
(L^{\prime }-\widetilde{L}^{\prime })p(0,T,x,y)+O(h^{1+\delta }\zeta _{\sqrt{%
T}}(y-x)).  \label{eq:P05}
\end{equation}
Now we furhter simplify our expansion of \ $p_{h}-p.$ We now show the
following expansion
\[
p_{h}(0,T,x,y)-p(0,T,x,y)
\]
\[
=\sqrt{h}(p\otimes \mathcal{F}_{1}[p_{\Delta }])(0,T,x,y)+h\left( p\otimes
\mathcal{F}_{2}[p_{\Delta }]\right) (0,T,x,y)
\]
\[
+h\left( p\otimes \mathcal{F}_{1}[p\otimes \mathcal{F}_{1}[p_{\Delta
}]]\right) (0,T,x,y)
\]
\begin{equation}
+\frac{h}{2}p\otimes (L_{\ast }^{2}-L^{2})p(0,T,x,y)-\frac{h}{2}p\otimes
(L^{\prime }-\widetilde{L}^{\prime })p(0,T,x,y)+O(h^{1+\delta }\zeta _{\sqrt{%
T}}(y-x)),  \label{eq:P06}
\end{equation}
where for \ $s\in \lbrack 0,t-h],t\in \{h,2h,...,T\}$%
\[
p_{\Delta }(s,t,z,y)=(\widetilde{p}\otimes _{h}^{\prime }\Phi )(s,t,z,y)
\]
\[
=\widetilde{p}(s,t,z,y)+\sum_{s\leq jh\leq t-h}h\int \widetilde{p}%
(s,jh,z,v)\Phi _{1}(jh,t,v,y)dv.
\]
Here $\Phi _{1}=H+H\otimes _{h}^{\prime }H+H\otimes _{h}^{\prime }H\otimes
_{h}^{\prime }H+....$We start from the calculation of $\ p\otimes \widetilde{%
L}\widetilde{\pi }_{1}(s,t,x,y).$%
\[
p\otimes \widetilde{L}\widetilde{\pi }_{1}(s,t,x,y)=\int_{s}^{t}d\tau \int
p(s,\tau ,x,v)(t-\tau )\sum_{\left| \nu \right| =3}\frac{\overline{\chi }%
_{\nu }(\tau ,t,y)}{\nu !}D_{v}^{\nu }(\widetilde{L}_{v}\widetilde{p}(\tau
,t,v,y))dv
\]
\[
=-\sum_{\left| \nu \right| =3}\frac{1}{\nu !}\int dv\left[
\int_{s}^{t}p(s,\tau ,x,v)\left( \int_{\tau }^{t}\chi _{\nu }(u,y)du\right)
\frac{\partial }{\partial \tau }D_{v}^{\nu }\widetilde{p}(\tau ,t,v,y)d\tau %
\right]
\]
\begin{equation}
=-\sum_{\left| \nu \right| =3}\frac{1}{\nu !}\int dv\int_{s}^{\frac{s+t}{2}%
}...-\sum_{\left| \nu \right| =3}\frac{1}{\nu !}\int dv\int_{\frac{s+t}{2}%
}^{t}...=I+II.  \label{eq:P07}
\end{equation}
Integrating by parts w.r.t. time variable we obtain
\[
I=-\sum_{\left| \nu \right| =3}\frac{1}{\nu !}\int dv\left[ p(s,\tau
,x,v)\left( \int_{\tau }^{t}\chi _{\nu }(u,y)du\right) D_{v}^{\nu }%
\widetilde{p}(\tau ,t,v,y)\mid _{\tau =s}^{\tau =(s+t)/2}\right.
\]
\[
\left. -\int_{s}^{\frac{s+t}{2}}D_{v}^{\nu }\widetilde{p}(\tau ,t,v,y)\left(
\frac{\partial p(s,\tau ,x,v)}{\partial \tau }\int_{\tau }^{t}\chi _{\nu
}(u,y)du-p(s,\tau ,x,v)\chi _{\nu }(\tau ,y)\right) d\tau \right]
\]
\[
=-\sum_{\left| \nu \right| =3}\frac{1}{\nu !}\int dv\left[ p(s,\frac{s+t}{2}%
,x,v)\left( \int_{\frac{s+t}{2}}^{t}\chi _{\nu }(u,y)du\right) D_{v}^{\nu }%
\widetilde{p}(\frac{s+t}{2},t,v,y)\right.
\]
\[
\left. +\sum_{\left| \nu \right| =3}\frac{1}{\nu !}\left( \int_{s}^{t}\chi
_{\nu }(u,y)du\right) D_{v}^{\nu }\widetilde{p}(s,t,x,y)\right]
+\sum_{\left| \nu \right| =3}\frac{1}{\nu !}\int_{s}^{\frac{s+t}{2}}d\tau
\left( \int_{\tau }^{t}\chi _{\nu }(u,y)du\right)
\]
\begin{equation}
\times \int L^{T}p(s,\tau ,x,v)D_{v}^{\nu }\widetilde{p}(\tau
,t,v,y)dv-\sum_{\left| \nu \right| =3}\frac{1}{\nu !}\int_{s}^{\frac{s+t}{2}%
}\chi _{\nu }(\tau ,y)d\tau \int p(s,\tau ,x,v)D_{v}^{\nu }\widetilde{p}%
(\tau ,t,v,y)dv.  \label{eq:P08}
\end{equation}
Analogously we get
\[
II=\sum_{\left| \nu \right| =3}\frac{1}{\nu !}\int p(s,\frac{s+t}{2}%
,x,v)\left( \int_{\frac{s+t}{2}}^{t}\chi _{\nu }(u,y)du\right) D_{v}^{\nu }%
\widetilde{p}(\frac{s+t}{2},t,v,y)dv
\]
\[
+\sum_{\left| \nu \right| =3}\frac{1}{\nu !}\int_{\frac{s+t}{2}}^{t}d\tau
\left( \int_{\tau }^{t}\chi _{\nu }(u,y)du\right) \int L^{T}p(s,\tau
,x,v)D_{v}^{\nu }\widetilde{p}(\tau ,t,v,y)dv
\]
\begin{equation}
-\sum_{\left| \nu \right| =3}\frac{1}{\nu !}\int_{\frac{s+t}{2}}^{t}\chi
_{\nu }(\tau ,y)d\tau \int p(s,\tau ,x,v)D_{v}^{\nu }\widetilde{p}(\tau
,t,v,y)dv.  \label{eq:P09}
\end{equation}
From (\ref{eq:P07})- (\ref{eq:P09}) we have
\begin{equation}
p\otimes \widetilde{L}\widetilde{\pi }_{1}(s,t,x,y)=\widetilde{\pi }%
_{1}(s,t,x,y)+p\otimes L\widetilde{\pi }_{1}(s,t,x,y)-p\otimes \widetilde{M}%
_{1}(s,t,x,y).  \label{eq:P10}
\end{equation}
and from (\ref{eq:P10}) we obtain
\[
\widetilde{\pi }_{1}(s,t,x,y)+p\otimes \Re _{1}(s,t,x,y)
\]
\[
=\widetilde{\pi }_{1}(s,t,x,y)+p\otimes L\widetilde{\pi }_{1}(s,t,x,y)-p%
\otimes \widetilde{L}\widetilde{\pi }_{1}(s,t,x,y)+p\otimes M_{1}(s,t,x,y)
\]
\begin{equation}
-p\otimes \widetilde{M}_{1}(s,t,x,y)=p\otimes M_{1}(s,t,x,y).  \label{eq:P11}
\end{equation}
It follows from (\ref{eq:P11}) and the definitions of \ the operations $%
\otimes $ and $\otimes _{h}^{\prime }$that \ ($\chi \lbrack s,jh]$ below
denotes the indicator of the interval $[s,jh]$)
\begin{eqnarray}  \label{eq:P12} &&
\sqrt{h}\left[ \widetilde{\pi }_{1}\otimes _{h}^{\prime }\Phi
(s,t,x,y)+(p\otimes \Re _{1})\otimes _{h}^{\prime }\Phi
(s,t,x,y)\right]
\\ \nonumber &&
=\sqrt{h}(\widetilde{\pi }_{1}+p\otimes \Re _{1})\otimes _{h}^{\prime }\Phi
(s,t,x,y)=\sqrt{h}(p\otimes M_{1})\otimes _{h}^{\prime }\Phi (s,t,x,y)
\\ \nonumber &&
=\sqrt{h}\sum_{0\leq jh\leq t-h}h\int (p\otimes
M_{1})(s,jh,x,z)\Phi (jh,t,z,y)dz
\\ \nonumber &&
=\sqrt{h}\sum_{0\leq jh\leq
t-h}h\int \left[ \int_{s}^{jh}du\int p(s,u,x,v)\right.
\\ \nonumber && \qquad
\left. \times M_{1}(u,jh,v,z)dv\right] \Phi
(jh,t,z,y)dz
\\ \nonumber && =\sqrt{h}\sum_{0\leq jh\leq t-h}h\int \left[
\int_{s}^{t}du\chi \lbrack s,jh]\int p(s,u,x,v)\right.
\\ \nonumber && \qquad
\left. \times M_{1}(u,jh,v,z)dv\right] \Phi (jh,t,z,y)dz
\\ \nonumber &&
=\sqrt{h}%
\int_{s}^{t}du\int p(s,u,x,v)\sum_{\left| \nu \right|
=3}\frac{\chi _{\nu }(u,v)}{\nu !}
\\ \nonumber && \qquad
\times D_{v}^{\nu }\left[ \sum_{0\leq jh\leq t-h}h\chi \lbrack
s,jh]\int
\widetilde{p}(u,jh,v,z)\Phi (jh,t,z,y)dz\right] dv\\ \nonumber && =\sqrt{h}%
\int_{s}^{t}du\int p(s,u,x,v)
\\ \nonumber && \qquad
\times \sum_{\left| \nu \right| =3}\frac{\chi _{\nu }(u,v)}{\nu
!}D_{v}^{\nu }p_{\Delta }(u,t,v,y)dv=\sqrt{h}(p\otimes
\mathcal{F}_{1})[p_{\Delta }](s,t,x,y).
\end{eqnarray}
Using similar arguments as in the proof of (\ref{eq:P12}) one can
show that
\[
h\left[ \widetilde{\pi }_{2}\otimes _{h}^{\prime }\Phi (s,t,x,y)+(p\otimes
\Re _{2})\otimes _{h}^{\prime }\Phi (s,t,x,y)\right.
\]
\begin{equation}
\left. +p\otimes _{h}(\Re _{3}\otimes _{h}^{\prime }\Phi
)(s,t,x,y)\right] =h(p\otimes \mathcal{F}_{2})[p_{\Delta
}](s,t,x,y)+hp\otimes \Pi _{1}\otimes _{h}^{\prime }\Phi (s,t,x,y)
. \label{eq:P13}
\end{equation}
For the first two terms in the right hand side of (\ref{eq:P05}) \ we obtain
from (\ref{eq:P12}) \ and (\ref{eq:P13})
\[
\sqrt{h}\left[ \widetilde{\pi }_{1}\otimes _{h}^{\prime }\Phi
(0,T,x,y)+(p\otimes \Re _{1})\otimes _{h}^{\prime }\Phi (0,T,x,y)\right]
\]
\[
+h\left[ \widetilde{\pi }_{2}\otimes _{h}^{\prime }\Phi (0,T,x,y)+p\otimes
(\Re _{2}\otimes _{h}^{\prime }\Phi )(0,T,x,y)+p\otimes _{h}(\Re _{3}\otimes
_{h}^{\prime }\Phi )(0,T,x,y)\right]
\]
\begin{equation}
=\sqrt{h}(p\otimes \mathcal{F}_{1})[p_{\Delta }](0,T,x,y)+h(p\otimes
\mathcal{F}_{2})[p_{\Delta }](s,t,x,y)+hp\otimes \Pi _{1}\otimes
_{h}^{\prime }\Phi (s,t,x,y).  \label{eq:P14}
\end{equation}
Using (\ref{eq:P12}) we also have
\[
h\left[ \widetilde{\pi }_{1}\otimes _{h}^{\prime }\Phi +p\otimes \left( \Re
_{1}\otimes _{h}^{\prime }\Phi \right) \right] \otimes _{h}^{\prime }(\Re
_{1}\otimes _{h}^{\prime }\Phi )(0,T,x,y)
\]
\[
=h(p\otimes \mathcal{F}_{1}[p_{\Delta }])\otimes _{h}^{\prime }(\Re
_{1}\otimes _{h}^{\prime }\Phi )(0,T,x,y)
\]
\[
=hp\otimes \mathcal{F}_{1}\left[ p_{\Delta }\otimes _{h}^{\prime }(\Re
_{1}\otimes _{h}^{\prime }\Phi )\right] (0,T,x,y).
\]
Note that
\[
hp\otimes \Pi _{1}\otimes _{h}^{\prime }\Phi (s,t,x,y)=h\int_{s}^{t}du\int
p(s,u,x,v)\sum_{\left| \nu \right| =3}\frac{\chi _{\nu }(u,v)}{\nu !}%
D_{v}^{\nu }[\widetilde{\pi }_{1}\otimes _{h}^{\prime }\Phi ](u,t,v,y)
\]
\[
=hp\otimes \mathcal{F}_{1}[\widetilde{\pi }_{1}\otimes _{h}^{\prime }\Phi
](s,t,x,y).
\]
For the proof of (\ref{eq:P06}) it remains to show that
\[
hp\otimes \mathcal{F}_{1}[\widetilde{\pi }_{1}\otimes _{h}^{\prime }\Phi
+p_{\Delta }\otimes _{h}^{\prime }(\Re _{1}\otimes _{h}^{\prime }\Phi
)](0,T,x,y)
\]
\begin{equation}
=h\left( p\otimes \mathcal{F}_{1}[p\otimes \mathcal{F}_{1}[p_{\Delta
}]]\right) (0,T,x,y)+O(h^{1+\delta }\zeta _{\sqrt{T}}(y-x)).  \label{eq:P15}
\end{equation}
We will show that
\begin{equation}
hp\otimes \mathcal{F}_{1}[(p-p_{\Delta })\otimes _{h}^{\prime }(\Re
_{1}\otimes _{h}^{\prime }\Phi )](0,T,x,y)=O(h^{1+\delta }\zeta _{\sqrt{T}%
}(y-x)) , \label{eq:P16}
\end{equation}
and
\[
hp\otimes \mathcal{F}_{1}[p\otimes _{h}^{\prime }(\Re _{1}\otimes
_{h}^{\prime }\Phi )](0,T,x,y)
\]
\begin{equation}
-hp\otimes \mathcal{F}_{1}[p\otimes (\Re _{1}\otimes _{h}^{\prime }\Phi
)](0,T,x,y)=O(h^{1+\delta }\zeta _{\sqrt{T}}(y-x)).  \label{eq:P17}
\end{equation}
Then (\ref{eq:P15}) will follow from (\ref{eq:P16}), (\ref{eq:P17}) and (\ref
{eq:P12}). The estimate (\ref{eq:P17}) can be shown analogously to (\ref
{eq:P04i}). An additional singularity arising from the derivatives in the
operator $\mathcal{F}_{1}[\cdot ]$ \ is neglected by the factor $h$ \ in (%
\ref{eq:P17}). To estimate (\ref{eq:P16}) note that \ from the definition of
$\ \Re _{1}$ and $\ \Phi $%
\begin{equation}
\left| (\Re _{1}\otimes _{h}^{\prime }\Phi )(jh,T,x,y)\right| \leq
C(\varepsilon )h^{-\varepsilon }(T-jh)^{\varepsilon -1}\phi _{\sqrt{T-jh}}.
\label{eq:P18}
\end{equation}
Then we use the following estimate which can be proved by the same method as
in Konakov (2006), pp.8-12, where an estimate for $(p-p^{d})(0,jh,x,y)$ \
was obtained.
\begin{equation}
\left| (p-p_{\Delta })(u,jh,v,z)\right| \leq Ch^{1/2}\phi _{\sqrt{jh-u}%
}(z-v).  \label{eq:P19}
\end{equation}
From (\ref{eq:P18}) \ and (\ref{eq:P19})
\begin{equation}
\left| (p-p_{\Delta })\otimes _{h}^{\prime }(\Re _{1}\otimes _{h}^{\prime
}\Phi )(u,T,v,y)\right| \leq C(\varepsilon )h^{1/2-\varepsilon
}(T-u)^{\varepsilon }\phi _{\sqrt{T-u}}(y-v)  \label{eq:P20}
\end{equation}
Now to estimate (\ref{eq:P16}) it's enough to estimate a typical term of a
corresponding sum, namely, for $\left| \nu \right| =3$ \ we have to estimate
\[
h\int_{0}^{T}du\int p(0,u,x,v)\frac{\chi _{\nu }(u,v)}{\nu !}D_{v}^{\nu
}[\sum_{\{j:u\leq jh\leq T-h\}}h\int (p-p_{\Delta })(u,jh,v,z)
\]
\[
\times (\Re _{1}\otimes _{h}^{\prime }\Phi
)(jh,T,z,y)dz]dv=h\int_{0}^{T/2}...+h\int_{T/2}^{T}...=I+II.
\]
To estimate $II$ \ we make an integration by parts transferring \ three
derivatives to $p(0,u,x,v)\frac{\chi _{\nu }(u,v)}{\nu !}.$ Using (\ref
{eq:P20}) we obtain the following \ estimate
\[
\left| II\right| \leq C(\varepsilon )h^{3/2-\varepsilon }\int_{T/2}^{T}\frac{%
(T-u)^{\varepsilon }}{u^{3/2}}du\phi _{\sqrt{T}}(y-x)
\]
\begin{equation}
\leq C(\varepsilon )h^{3/2-\varepsilon }T^{\varepsilon }\phi _{\sqrt{T}%
}(y-x).  \label{eq:P21}
\end{equation}
To estimate $I$ \ we consider two cases, namely, a) $jh-u\geq T/4$ \ and \
b) $jh-u\leq T/4$ $\Longrightarrow $ $T-jh\geq T/4$ . Analogously to (11) \
in Konakov (2006) the difference $h(p-p_{\Delta })$ can be represented \ as
\[
h(p-p_{\Delta })(u,jh,v,z)=h(p\otimes H-p\otimes _{h}^{\prime }H)(u,jh,v,z)
\]
\[
+h(p\otimes H-p\otimes _{h}^{\prime }H)\otimes _{h}^{\prime }\Phi
_{1}(u,jh,v,z)
\]
\[
=h\int_{u}^{j^{\ast }h}d\tau \int p(u,\tau ,v,z^{\prime })H(\tau
,jh,z^{\prime },z)dz^{\prime }
\]
\[
+h\sum_{i=j^{\ast }}^{j-1}\int_{ih}^{(i+1)h}d\tau \int (\lambda (\tau
,z^{\prime })-\lambda (ih,z^{\prime }))dz^{\prime }+h(p\otimes H-p\otimes
_{h}^{\prime }H)\otimes _{h}^{\prime }\Phi _{1}(u,jh,v,z)
\]
\begin{equation}
=I^{\prime }+II^{\prime }+III^{\prime }.  \label{eq:P21a}
\end{equation}
where $\lambda (\tau ,z^{\prime })=p(u,\tau ,v,z^{\prime })H(\tau
,jh,z^{\prime },z),\Phi _{1}(ih,jh,z^{\prime }z)=H(ih,jh,z^{\prime
}z)+$ $H\otimes _{h}^{\prime }$\\ $H(ih,jh,z^{\prime }z)+...,j^{\ast }=j^{\ast }(u)=[%
\frac{u}{h}]+1$ (with a convention $[x]=x-1$ for integer $x$%
). \ For $I^{\prime },$ case a), we have $jh-\tau >T/5$ \ for \ $n\ \ $large
enough. With a substitution $v+v^{\prime }=z^{\prime }$ \ we obtain
\[
\left| D_{v}^{\nu }h\int_{u}^{j^{\ast }h}d\tau \int p(u,\tau ,v,z^{\prime
})H(\tau ,jh,z^{\prime },z)dz^{\prime }\right|
\]
\[
=\left| D_{v}^{\nu }h\int_{u}^{j^{\ast }h}d\tau \int p(u,\tau ,v,v+v^{\prime
})H(\tau ,jh,v+v^{\prime },z)dv^{\prime }\right|
\]
\[
\leq Ch\int_{u}^{j^{\ast }h}\frac{d\tau }{(jh-\tau )^{2}}\phi _{\sqrt{jh-u}%
}(z-v)\leq Ch^{2}T^{-2}\phi _{\sqrt{jh-u}}(z-v)
\]
\begin{equation}
\leq Cn^{-2}\phi _{\sqrt{jh-u}}(z-v)=CT^{2}h^{2}\phi _{\sqrt{jh-u}}(z-v)
\label{eq:P22}
\end{equation}
To get (\ref{eq:P22}) we used an estimate from Friedman (1964) (Theorem 7,
page 260)
\[
\left| D_{v}^{\nu }p(u,\tau ,v,v+v^{\prime })\right| \leq C(\tau
-u)^{-d/2}\exp \left( \frac{C\left| v^{\prime }\right| }{\tau -u)}\right) .
\]
For $\int I^{\prime }(u,jh,v,z)(\Re _{1}\otimes _{h}^{\prime }\Phi
)(jh,T,z,y)dz,$ case b), it is enough to estimate for $\left| \nu \right| =3$%
\[
h\int_{u}^{j^{\ast }h}d\tau \int [p(u,\tau ,v,v+v^{\prime })(\sigma
_{lk}(\tau ,v+v^{\prime })-\sigma _{lk}(\tau ,z))]D_{v^{\prime }}^{\nu }%
\frac{\partial ^{2}\widetilde{p}(\tau ,jh,v+v^{\prime },z)}{\partial
v_{l}^{\prime }\partial v_{k}^{\prime }}dv^{\prime }
\]
\begin{equation}
\times (\Re _{1}\otimes _{h}^{\prime }\Phi )(jh,T,z,y)dz.  \label{eq:P23}
\end{equation}
Transferring \ five derivatives from $\widetilde{p}$ \ to $(\Re _{1}\otimes
_{h}^{\prime }\Phi )(jh,T,z,y)$ and using the following estimate for $\left|
\mu \right| =5$%
\[
\left| D_{v^{\prime }}^{\mu }\widetilde{p}(\tau ,jh,v+v^{\prime
},z)+D_{z}^{\mu }\widetilde{p}(\tau ,jh,v+v^{\prime },z)\right|
\]
\[
\leq C(jh-\tau )^{-d/2}\phi _{\sqrt{jh-\tau }}(z-v-v^{\prime })
\]
we obtain that (\ref{eq:P23}) does not exceed
\[
C(j^{\ast }h-\tau )hT^{-7/2}\phi _{\sqrt{T-u}}(y-v)\leq Ch^{1+\delta
}T^{1-\delta }n^{-(1-\delta -7\varkappa /2)}\phi _{\sqrt{T-u}}(y-v)
\]
\begin{equation}
=o(h^{1+\delta }T^{1-\delta })\phi _{\sqrt{T-u}}(y-v).  \label{eq:P24}
\end{equation}
We used that for any \ $0<$ $\delta <1$ \ \ $\varkappa
<\frac{2-2\delta }{7}$ \ (see condition (B2)). \ To estimate $\int
II^{\prime }(u,jh,v,z)(\Re _{1}\otimes _{h}^{\prime }\Phi
)(jh,T,z,y)dz$ we use a decomposition \ (12) in \ Konakov (2006).
The estimate for terms in $\ II^{\prime }$ \ containing \ the
first derivatives \ $\lambda ^{\prime }(ih,z^{\prime })$ \ we use
identity (14) from Konakov (2006) and similar arguments to already
used in estimation of $\int I^{\prime }(u,jh,v,z)(\Re
_{1}\otimes _{h}^{\prime }\Phi )(jh,T,z,y)dz.$ The estimate for terms in $%
II^{\prime }$ containing second derivatives $\lambda ^{^{\prime \prime
}}(ih,z^{\prime })$ \ follows from (22) and (23) in Konakov (2006). At last
for $III^{\prime }$ \ the same estimates hold because of smoothing
properties of the convolution ...$\otimes _{h}^{^{\prime }}\Phi
_{1}(u,jh,v,z).$ This implies (\ref{eq:P16}) and, hence, the expansion (\ref
{eq:P06}).

\textit{Asymptotic replacement of \ }$p_{\Delta }$ \textit{by }$p.$ We shall
compare $hp\otimes \mathcal{F}_{2}[p_{\Delta }](0,T,x,y)$ \ and \ $hp\otimes
\mathcal{F}_{2}[p](0,T,x,y)$ . Note that for \ $\ 2\varkappa <\delta <\frac{2%
}{5},\left| \nu \right| =4$%
\[
\left| h\int_{0}^{h^{\delta }}du\int p(0,u,x,z)\chi _{\nu }(u,z)D_{z}^{\nu
}p(u,T,z,y)dz\right| \leq Ch^{1+\delta }(T-h^{\delta })^{-2}\phi _{\sqrt{T}%
}(y-x)
\]
\begin{equation}
\leq Ch^{1+\delta }\frac{n^{2\varkappa }}{(Tn^{\varkappa }-n^{\varkappa
}h^{\delta })^{2}}\leq Ch^{1+(\delta -2\varkappa )}T^{2\varkappa }\phi _{%
\sqrt{T}}(y-x),  \label{eq:P25}
\end{equation}
and, analogously,
\begin{equation}
\left| h\int_{T-h^{\delta }}^{T}du\int D_{z}^{\nu }[p(0,u,x,z)\chi _{\nu
}(u,z)]p(u,T,z,y)dz\right| \leq Ch^{1+(\delta -2\varkappa )}T^{2\varkappa
}\phi _{\sqrt{T}}(y-x).  \label{eq:P26}
\end{equation}
The same estimates hold for $p_{\Delta }(u,T,z,y).$ Hence, it is
enough to consider $u\in \lbrack h^{\delta },T-h^{\delta }].$ We
consider
\[
h\int_{h^{\delta }}^{T-h^{\delta }}du\int p(0,u,x,z)\chi _{\nu
}(u,z)D_{z}^{\nu }(p-p_{\Delta })(u,T,z,y)dz
\]
\begin{equation}
=h\int_{h^{\delta }}^{T/2}...+h\int_{T/2}^{T-h^{\delta }}...=I+II.
\label{eq:P27}
\end{equation}
By (\ref{eq:P19})
\[
\left| II\right| =\left| h\int_{T/2}^{T-h^{\delta }}du\int D_{z}^{\nu
}[p(0,u,x,z)\chi _{\nu }(u,z)](p-p_{\Delta })(u,T,z,y)dz\right|
\]
\[
\leq Ch^{3/2}n^{\varkappa }\phi _{\sqrt{T}}(y-x)=Ch^{3/2-\varkappa
}T^{\varkappa }\phi _{\sqrt{T}}(y-x)
\]
\begin{equation}
=Ch^{1+\gamma }\phi _{\sqrt{T}}(y-x),\gamma >0  \label{eq:P28}
\end{equation}
For $u\in \lbrack h^{\delta },T/2]$%
\[
\left| I\right| =\left| h\int_{h^{\delta }}^{T/2}du\int D_{z}^{\nu
}[p(0,u,x,z)\chi _{\nu }(u,z)](p-p_{\Delta })(u,T,z,y)dz\right|
\]
\begin{equation}
\leq Ch^{3/2-\delta }\phi _{\sqrt{T}}(y-x),  \label{eq:P29}
\end{equation}
where condition $\ \delta <\frac{2}{5}$ implies $3/2-\delta >1.$ It follows
from (\ref{eq:P25})-(\ref{eq:P29}) that
\begin{equation}
hp\otimes \mathcal{F}_{2}[p_{\Delta }](0,T,x,y)-hp\otimes \mathcal{F}%
_{2}[p](0,T,x,y)=O(h^{1+\gamma }\phi _{\sqrt{T}}(y-x)).  \label{eq:P30}
\end{equation}
To prove that
\[
hp\otimes \mathcal{F}_{1}[p\otimes \mathcal{F}_{1}[p-p_{\Delta
}]]=O(h^{1+\delta }\phi _{\sqrt{T}}(y-x))
\]
we consider a typical summand for fixed $\nu ,\left| \nu \right| =3,$%
\begin{equation}
h\int_{0}^{T}du\int p(0,u,x,z)\chi _{\nu }(u,z)D_{z}^{\nu }\left[
\int_{u}^{T}d\tau \int p(u,\tau ,z,v)\chi _{\nu }(\tau ,v)D_{v}^{\nu
}(p-p_{\Delta })(\tau ,T,v,y)dv\right] dz.  \label{eq:P31}
\end{equation}
As before it's enough to consider $u\in \lbrack h^{\delta },T-h^{\delta }].$
\ The integral in (\ref{eq:P31}) is a sum of \ four integrals
\[
I_{1}=h\int_{h^{\delta }}^{T/2}du\int ...D_{z}^{\nu }\int_{u}^{(T+u)/2}d\tau
\int ...
\]
\[
I_{2}=h\int_{h^{\delta }}^{T/2}du\int ...D_{z}^{\nu }\int_{(T+u)/2}^{T}d\tau
\int ...
\]
\[
I_{3}=h\int_{T/2}^{T-h^{\delta }}du\int ...D_{z}^{\nu
}\int_{u}^{(T+u)/2}d\tau \int ...
\]
\begin{equation}
I_{4}=h\int_{T/2}^{T-h^{\delta }}du\int ...D_{z}^{\nu
}\int_{(T+u)/2}^{T}d\tau \int ....  \label{eq:P32}
\end{equation}
Note \ that \ $\tau -u\geq T/4$ \ in the integrand in $I_{2}$ . Integrating
by parts w.r.t. $v$ \ we obtain by (\ref{eq:P19})
\begin{equation}
\left| I_{2}\right| \leq Ch^{3/2-\varkappa }T^{\varkappa }\phi _{\sqrt{T}%
}(y-x)).  \label{eq:P33}
\end{equation}
Furthermore, $u\geq T/2,\tau -u\geq h^{\delta }/2,T-u\geq
h^{\delta }$ \ in the integrand in $I_{4}.$ Integrating by parts
w.r.t. $z$ \ we obtain
\begin{equation}
\left| I_{4}\right| \leq Ch^{3/2-\delta }T^{-1/2}\phi _{\sqrt{T}}(y-x))\leq
CT^{\varkappa /2}h^{3/2-\varkappa /2-\delta }\phi _{\sqrt{T}}(y-x)),
\label{eq:P34}
\end{equation}
where, by our choice of $\delta ,$ $\ 3/2-\varkappa /2-\delta >1.$ To
estimate $I_{3}$ we use the representation
\[
(p-p_{\Delta })(\tau ,T,v,y)=(p\otimes H-p\otimes _{h}^{\prime }H)(\tau
,T,v,y)
\]
\[
+(p\otimes H-p\otimes _{h}^{\prime }H)\otimes _{h}^{\prime }\Phi _{1}(\tau
,T,v,y)
\]
\[
=\int_{\tau }^{j^{\ast }h}ds\int p(\tau ,s,v,w)H(s,T,w,y)dw
\]
\[
+\frac{h}{2}[p\otimes _{h}^{\prime }(H_{1}+A_{0})](\tau ,T,v,y)
\]
\[
+\frac{1}{2}\sum_{i=j^{\ast
}}^{n-1}\int_{ih}^{(i+1)h}(t-ih)^{2}\int_{0}^{1}(1-\gamma
)\sum_{k=1}^{4}\int p(\tau ,s,v,w)A_{k}(s,T,w,y)\mid _{s=ih+\gamma
(t-ih)}dwd\gamma dt,
\]
\begin{equation}
+(p\otimes H-p\otimes _{h}^{\prime }H)\otimes _{h}^{\prime }\Phi _{1}(\tau
,T,v,y)  \label{eq:P35}
\end{equation}
where $j^{\ast }=j^{\ast }(\tau )=[\tau /h]+1$ (with a convention $[x]=x-1$
\ for integer $x$) , $H_{1}$ and $A_{k},k=0,1,2,3,4,$ are defined in Konakov
(2006) and
\[
\Phi _{1}(ih,i^{\prime }h,z,z^{\prime })=H(ih,i^{\prime }h,z,z^{\prime
})+H\otimes _{h}^{\prime }H(ih,i^{\prime }h,z,z^{\prime })+...
\]
To estimate $D_{v}^{\nu }(p-p_{\Delta })(\tau ,T,v,y)$ \ we note \ that
\[
h\left| D_{v}^{\nu }\int_{\tau }^{j^{\ast }h}ds\int p(\tau
,s,v,w)H(s,T,w,y)dw\right|
\]
\[
=h\left| D_{v}^{\nu }\int_{\tau }^{j^{\ast }h}ds\int p(\tau ,s,v,v+w^{\prime
})H(s,T,v+w^{\prime },y)dw^{\prime }\right|
\]
\begin{equation}
\leq Ch\int_{\tau }^{j^{\ast }h}\frac{ds}{(T-s)^{2}}\phi _{\sqrt{T-\tau }%
}(y-v)\leq Ch^{2-2\delta }\phi _{\sqrt{T-\tau }}(y-v).  \label{eq:P36}
\end{equation}
Furthermore,
\[
\left| \frac{h^{2}}{2}D_{v}^{\nu }[p\otimes _{h}^{\prime }H_{1}](\tau
,T,v,y)\right| =\left| \frac{h^{2}}{2}D_{v}^{\nu }\sum_{\tau \leq jh\leq
T-h}h\int p(\tau ,jh,v,w)H_{1}(jh,T,w,y)dw\right|
\]
\[
\leq \frac{h^{2}}{2}\left| \sum_{\tau \leq jh\leq T-h^{\delta }/2}h\int
D_{v}^{\nu }[p(\tau ,jh,v,v+w^{\prime })H_{1}(jh,T,v+w^{\prime
},y)]dw^{\prime }\right|
\]
\[
+\frac{h^{2}}{2}\left| C\sum_{i,k=1}^{d}\sum_{T-h^{\delta }/2<jh\leq
T-h}h\int D_{w^{\prime }}^{\nu +e_{i}+e_{k}}[p(\tau ,jh,v,v+w^{\prime })]%
\widetilde{p}(jh,T,v+w^{\prime },y)dw^{\prime }\right|
\]
\begin{equation}
\leq Ch^{2-2\delta }\phi _{\sqrt{T-\tau }}(y-v)+Ch^{2-5\delta /2}\phi _{%
\sqrt{T-\tau }}(y-v).  \label{eq:P37}
\end{equation}
Because of a structure of the operator $A_{0}$ (see Konakov (2006)) it is
enough to estimate for fixed $i,l,k$%
\begin{equation}
\frac{h^{2}}{2}D_{v}^{\nu }\sum_{\tau \leq jh\leq T-h}h\int D_{w^{\prime
}}^{e_{k}}p(\tau ,jh,v,v+w^{\prime })\frac{\partial ^{2}\widetilde{p}%
(jh,T,v+w^{\prime },y)}{\partial w_{i}^{\prime }\partial w_{l}^{\prime }}%
dw^{\prime }  \label{eq:P38}
\end{equation}
With the same decomposition as in (\ref{eq:P37}) we obtain that (\ref{eq:P38}%
) \ does not exceed
\begin{equation}
Ch^{2-5\delta /2}\phi _{\sqrt{T-\tau }}(y-v).  \label{eq:P39}
\end{equation}
From (\ref{eq:P37}) and (\ref{eq:P39}) we obtain that
\begin{equation}
\frac{h^{2}}{2}\left| D_{v}^{\nu }[p\otimes _{h}^{\prime
}(H_{1}+A_{0})](\tau ,T,v,y)\right| \leq Ch^{1+\gamma }\phi _{\sqrt{T-\tau }%
}(y-v)  \label{eq:P40}
\end{equation}
for some $\gamma >0.$ It remains to estimate the last summand in (\ref
{eq:P35}). It follows from the structure of the operators $A_{k},k=1,2,3,4,$
(see Konakov (2006)) that it is enough to estimate
\begin{eqnarray}  \label{eq:P41}&&
\frac{1}{2}\sum_{i=j^{\ast
}}^{n-1}\int_{ih}^{(i+1)h}(t-ih)^{2}\int_{0}^{1}(1-\gamma
)\\ \nonumber && \qquad \sum_{k=1}^{4}\int D_{v}^{\nu }\left[ p(\tau ,s,v,v+w^{\prime })\frac{%
\partial ^{4}\widetilde{p}(s,T,v+w^{\prime },y)}{\partial w_{i}^{\prime
}\partial w_{l}^{\prime }\partial w_{p}^{\prime }\partial w_{q}^{\prime }}%
\right] \mid _{s=ih+\gamma (t-ih)}dw^{\prime }d\gamma dt
\end{eqnarray}
for fixed $i,j,p,q.$ As in Konakov (2006), p. 12, we obtain that (\ref
{eq:P41}) does not exceed
\[
Ch^{2}\phi _{\sqrt{T-\tau }}(y-v)\int_{0}^{1}z^{2}\int_{0}^{1}(1-\gamma
)\sum_{i=j^{\ast }}^{n-1}h\frac{1}{[(ih-\tau )+\gamma hz]^{3/2}}\frac{1}{%
[(n-\gamma z)h-ih]^{2}}d\gamma dz
\]
\[
=Ch^{2}\phi _{\sqrt{T-\tau }}(y-v)\int_{0}^{1}z^{2}\int_{0}^{1}(1-\gamma
)\sum_{\{i:j^{\ast }h\leq ih\leq \tau +h^{\delta }/4\}}...
\]
\begin{equation}
+Ch^{2}\phi _{\sqrt{T-\tau }}(y-v)\int_{0}^{1}z^{2}\int_{0}^{1}(1-\gamma
)\sum_{\{i:\tau +h^{\delta }/4<ih\leq T-h\}}...=I^{\prime \prime
}+II^{\prime \prime }.  \label{eq:P42}
\end{equation}
\[
\left| I^{\prime \prime }\right| \leq Ch^{2}h^{-5\delta /2}\phi _{\sqrt{%
T-\tau }}(y-v)\int_{0}^{1}z^{2}\int_{0}^{1}(1-\gamma )\sum_{\{i:j^{\ast
}h\leq ih\leq \tau +h^{\delta }/4\}}h\frac{1}{[(ih-\tau )+\gamma hz]}d\gamma
dz
\]
\[
\leq Ch^{2-\varepsilon }h^{-5\delta /2}\phi _{\sqrt{T-\tau }%
}(y-v)\int_{0}^{1}z^{2-\varepsilon }\int_{0}^{1}\frac{(1-\gamma )}{\gamma
^{\varepsilon }}\sum_{\{i:j^{\ast }h\leq ih\leq \tau +h^{\delta }/4\}}h\frac{%
1}{[(ih-\tau )+\gamma hz]^{1-\varepsilon }}d\gamma dz
\]
\begin{equation}
\leq C(\varepsilon )h^{2-\varepsilon -5\delta /2}\phi _{\sqrt{T-\tau }}(y-v).
\label{eq:P43}
\end{equation}
Using inequality $(h-\gamma z)h-ih=(n-i)h-\gamma zh\geq h(1-\gamma z)\geq
h(1-\gamma )$ \ we obtain that
\[
\left| II^{\prime \prime }\right| \leq Chh^{-5\delta /2}\phi _{\sqrt{T-\tau }%
}(y-v)\int_{0}^{1}z^{2}\int_{0}^{1}d\gamma \sum_{\{i:\tau +h^{\delta
}/4<ih\leq T-h\}}h
\]
\begin{equation}
\leq Ch^{1-5\delta /2}\phi _{\sqrt{T-\tau }}(y-v).  \label{eq:P44}
\end{equation}
Now from (\ref{eq:P35}), (\ref{eq:P36}), (\ref{eq:P40}), (\ref{eq:P41}), (%
\ref{eq:P43}) and (\ref{eq:P43}) \ we obtain that
\begin{equation}
\left| D_{v}^{\nu }(p\otimes H-p\otimes _{h}^{\prime }H)(\tau ,T,v,y)\right|
\leq Ch^{\gamma }\phi _{\sqrt{T-\tau }}(y-v)  \label{eq:P45}
\end{equation}
for some positive $\gamma .$ The last summand in the right hand side of (\ref
{eq:P35}) admits the same estimate (\ref{eq:P45}) because of the smoothing
properties of the operation $\otimes _{h}^{\prime }.$ Hence,
\begin{equation}
\left| D_{v}^{\nu }(p-p_{\Delta })(\tau ,T,v,y)\right| \leq Ch^{\gamma }\phi
_{\sqrt{T-\tau }}(y-v).  \label{eq:P45a}
\end{equation}
Making a chaige of variables $v=z+v^{\prime }$ \ in (\ref{eq:P31}) we get
that the integral w.r.t. $v$ \ is equal to
\begin{equation}
D_{z}^{\nu }\left[ \int_{u}^{T}d\tau \int p(u,\tau ,z,z+v^{\prime })\chi
_{\nu }(\tau ,v)D_{v}^{\nu }(p-p_{\Delta })(\tau ,T,z+v^{\prime
},y)dv^{\prime }\right] .  \label{eq:P46}
\end{equation}
Taking into account (\ref{eq:P45a}) and making an integration by parts in (%
\ref{eq:P46}) we obtain that (\ref{eq:P46}) does not exceed
\begin{equation}
Ch^{\gamma }\int_{u}^{T}\frac{d\tau }{(\tau -u)^{3/2}}\phi _{\sqrt{T-u}%
}(y-z)\leq \frac{Ch^{\gamma }}{\sqrt{T-u}}\phi _{\sqrt{T-u}}(y-z)
\label{eq:P47}
\end{equation}
From (\ref{eq:P31}) and (\ref{eq:P47}) we obtain that
\[
\left| I_{3}\right| \leq Ch^{1+\gamma }\phi _{\sqrt{T}}(y-x).
\]
The estimate for $I_{1}$ \ may be proved absolutely analogously to the
estimate for $I_{3}.$ Thus, we proved that
\begin{equation}
hp\otimes \mathcal{F}_{1}[p\otimes \mathcal{F}_{1}[p-p_{\Delta
}]]=O(h^{1+\delta }\phi _{\sqrt{T}}(y-x))  \label{eq:P48}
\end{equation}
The estimate
\[
h^{1/2}p\otimes \mathcal{F}_{1}[p-p_{\Delta }]=O(h^{1+\delta }\phi _{\sqrt{T}%
}(y-x))
\]
may be proved by using the same decomposition of \ $p-p_{\Delta }.$ This
completes the proof of \ Theorem 1.

\section*{References.}

\begin{description}
\item [Bally,~V.~and~Talay,~D.~\ ](1996a). The law of the Euler
scheme for stochastic differential equations: I. Convergence rate
of the distribution function. \textsl{Probab. Theory Rel. Fields}
\textbf{104}, 43-60. \item [Bally,~V.~and~Talay,~D.~\ ](1996b).
The law of the Euler scheme for stochastic differential equations:
II. Convergence rate of the density. \textsl{Monte Carlo Methods
Appl.} \textbf{2}, 93-128. \item [Bhattacharya,~R.~and~Rao,~R.\
](1976). \textsl{Normal approximations and asymptotic expansions.}
John Wiley \& Sons, New York. \item [Friedman,~A.\ ](1964).
\textsl{Partial differential equations of parabolic type.}
Prentice-Hall, Englewood Cliffs, New Jersey.

\item [Guyon,~J.]\ (2006). Euler scheme and tempered
distributions. \textit{Stoch. Proc. and their Appl.}, to appear.

\item [Jacod,~J.~and Protter,~P.]\ (1998). Asymptotic error
distributions for the Euler method for stochastic differential
equations. \textit{Ann. Prob. } \textbf{26}, 267-307.

\item [Jacod,~J.]\ (2004). The Euler scheme for L\'{e}vy driven
stochastic differential equations: limit theorems. \textit{Ann.
Prob. } \textbf{32}, 1830-1872.

\item [Jacod,~J., Kurtz,~T., M\'{e}l\'{e}ard,~S.~and Protter,~P.]\
(2005). The approximate Euler method for L\'{e}vy driven
stochastic differential equations. \textit{Ann. de l´I.H.P.}
\textbf{41}, 523-558.

 \item [Kloeden,~P.~E.~and~Platen,~E.~\ ](1992).
\textit{Numerical solution of stochastic differential equations.}
Springer, Berlin, Heidelberg. \item
[Konakov,~V.~D.~and~Mammen,~E.]\ (2000). Local limit theorems for
transition densities of Markov chains converging to diffusions.
\textit{Probability Theory and rel. Fields} 117, 551--587. \item
[Konakov,~V.~D.~and~Mammen,~E.]\ (2001). Local approximations of
Markov random walks by diffusions. \textit{Stochastic Processes
and their Applications} 96, 73-98. \item
[Konakov,~V.~D.~and~Mammen,~E.]\ (2002). Edgeworth type expansions
for Euler schemes for stochastic differential equations.
\textit{Monte Carlo Methods and Applications} 8, 271-286.

\item [Konakov,~V.~D.~and~Mammen,~E.]\ (2005). Edgeworth-type
expansions for transition densities of Markov chains converging to
diffusions. \textit{Bernoulli} 4, 591-641.

\item [Konakov,~V.~D.]\ (2006). Small time asymptotics in local
limit theorems for Markov chains converging to diffusions.
\textit{Prepublications LPMA.}
 \item [Konakov,~V.~D.~and~Molchanov,~S.]\textbf{A}.\
(1984). On the convergence of Markov chains to diffusion
processes. \textsl{Teor. Veroyatn. Mat. Stat.}\textbf{31}, 51-64,
(in russian) {[}English translation in Theory Probab. Math. Stat.
(1985) \textbf{31}, 59-73{]}. \item[Lady\u{z}enskaja, O. A.,
Solonnikov,  V. A. and and Ural\'{}ceva (1968).] {\it Linear and
quasi-linear equations of parabolic type.} Amer. Math. Soc.,
Providence, Rhode Island.
 \item
[McKean,~H.~P.~and~Singer,~I.~M.~\ ](1967). Curvature and the
eigenvalues of the Laplacian. \textsl{J. Differential Geometry}
\textbf{1}, 43-69. \item [Mil'shtein,~G.~N.](1974). Approximate
integration of stochastic differential equations. \emph{Theory
Probab. Appl.,} \textbf{\emph{19,}} 557-562. \item [Protter,P.~
and Talay,D.~]\ (1997). The Euler scheme for L\'{e}vy driven
stochastic differential equations. \textit{Ann. Prob. }
\textbf{25}, 393-423.

 \item
[Stroock,~D.~W.,~Varadhan,~S.~R.~](1979). \textsl{Multidimensional
diffusion processes.}\\
 Springer, Berlin, Heidelberg, New York.
\item [Skorohod,~A.~V.\ ](1965). \textsl{Studies in the theory of
random processes.} Addison-Wesley, Reading,
Massachussetts.{[}English translation of Skorohod, A. V.\ (1961).
\textsl{Issledovaniya po teorii sluchainykh processov.} Kiev
University Press{]} \item [Skorohod,~A.~V.~](1987).
\textsl{Asymptotic methods for Stochastic differential equations.}
(in russian). Kiev, Naukova dumka.\end{description}

\end{document}